\documentclass[journal,article,submit,moreauthors]{Definitions/mdpi} 
\usepackage{tikz}
\usetikzlibrary{shapes.geometric, arrows.meta, positioning, shadows}
\usepackage[edges]{forest}

\firstpage{1} 
\pubvolume{1}
\issuenum{1}
\articlenumber{0}
\pubyear{2026}
\copyrightyear{2026}
\datereceived{ } 
\daterevised{ } 
\dateaccepted{ } 
\datepublished{ } 

\usepackage{amsmath,amssymb,amsfonts}
\usepackage{bm}
\usepackage{graphicx}
\usepackage{multirow}
\usepackage{booktabs}
\usepackage{algorithm}
\usepackage{algorithmicx}
\usepackage{algpseudocode}
\usepackage{url}
\usepackage{stmaryrd}
\usepackage{textcomp}
\usepackage{enumitem}

\mathchardef\mhyphen="2D

\newcommand\rnumber{\operatorname{-}}
\newcommand{\acknowledgements}[1]{\section*{Acknowledgements}#1}

\newcommand{\ten}[1]{\underline{\mathbf{#1}}}
\newcommand{\A}{\underline{\mathbf A}}
\newcommand{\T}{\underline{\mathbf T}}
\newcommand{\W}{\underline{\mathbf W}}
\newcommand{\G}{\underline{\mathbf G}}
\newcommand{\D}{\underline{\mathbf D}}
\newcommand{\X}{\underline{\mathbf X}}
\newcommand{\C}{\underline{\mathbf C}}
\newcommand{\R}{\underline{\mathbf R}}
\newcommand{\Y}{\underline{\mathbf Y}}
\newcommand{\I}{\underline{\mathbf I}}
\newcommand{\B}{\underline{\mathbf B}}
\newcommand{\Q}{\underline{\mathbf Q}}
\newcommand{\U}{\underline{\mathbf U}}

\newcommand{\V}{\underline{\mathbf V}}

\newcommand{\Z}{\underline{\mathbf Z}}
\newcommand{\cH}{\underline{\mathbf H}}
\newcommand{\cG}{\underline{\mathbf G}}

\Title{Efficient Techniques for Low-Rank Tensor Approximation and Applications in Robust Object Detection}

\Author{Salman Ahmadi-Asl $^{1,2}$\orcidA{}, Naeim Rezaeian $^{3}$\orcidB{}, Cesar F. Caiafa $^{4}$\orcidC{} and André L. F. de Almeida $^{5}$\orcidD{}} 

\AuthorNames{Salman Ahmadi-Asl, Naeim Rezaeian, Cesar F. Caiafa and André L. F. de Almeida}

\address{%
$^{1}$ \quad Research Center of the Artificial Intelligence Institute, Innopolis University, 420500, Russia;\\
$^{2}$ \quad Lab of Machine Learning and Knowledge Representation, Innopolis University, 420500, Russia;\\
$^{3}$ \quad Department of Applied Artificial Intelligence, Faculty of Artificial Intelligence, RUDN University, 117198, Moscow, Russia;\\
$^{4}$ \quad Instituto Argentino de Radioastronomía—CCT La Plata, CONICET/CIC-PBA/UNLP, Villa Elisa 1894, Argentina;\\
$^{5}$ \quad Department of Teleinformatics Engineering, Federal University of Ceará, Fortaleza, Brazil}

\corres{Correspondence: s.ahmadiasl@innopolis.ru}

\abstract{
This paper introduces efficient randomized fixed-precision and single-pass algorithms for low-tubal-rank approximation of third-order tensors. The proposed fixed-precision algorithms are faster and more efficient than the existing algorithms for approximating the truncated tensor SVD (T-SVD). Besides, unlike existing single-pass methods, which directly extend early, unstable matrix algorithms, the proposed approach adapts enhanced and stabilized matrix techniques to the tensor setting. Through extensive numerical experiments, we identify a critical flaw in current single-pass algorithms: using sketching parameters of equal size often produces ill-conditioned tensor least-squares problems, leading to inaccurate approximations. The proposed algorithms are demonstrably robust to this issue, achieving superior performance under identical conditions. We also evaluate the robustness of existing single-pass methods on real-world data tensors, including images and videos, a topic that has not been thoroughly examined before. Numerical results confirm the effectiveness of the proposed methods. Three applications are presented: image compression, video super-resolution, and deep learning.}


\begin{document}
\section{Introduction}\label{sec:intro}
Tensor decompositions are foundational tools in signal processing, machine learning, and deep learning \cite{sidiropoulos2017tensor,comon2009tensor}. Widely used formats include Tucker \cite{de2000multilinear}, tensor train and tensor ring \cite{oseledets2011tensor,zhao2016tensor}, tensor singular value decomposition (T-SVD) \cite{kilmer2013third}, block term decomposition, and constrained factorizations \cite{confac,stegeman2009}. These techniques underpin a broad range of applications, such as tensor completion \cite{chen2026rank,shu2026robust,yang2026bayesian}, recommender systems \cite{sun2026pcnet,lu2025enabling}, machine learning \cite{miao2026multi,he2025fully,navarro2025low}, and wireless communications \cite{Almeida_Elsevier_2007,Favier2014,deAraujoSAM2020,de2021channel,de2025channel}.

For large-scale tensors, exact decomposition methods often become computationally infeasible. Over the past decade, randomized algorithms have emerged as a practical alternative. Compared to deterministic approaches, they are typically faster, more memory-efficient, and well-suited to large data. Furthermore, they offer controllable accuracy, allowing users to trade off computational cost against approximation quality. 

The rapid growth of data has led to tensors of ever-increasing size. Low-rank approximations are critical in this context, particularly for tensor completion \cite{song2019tensor} and recommender systems \cite{frolov2017tensor}. The difficulty escalates when tensors must be stored out-of-core\footnote{Out-of-core storage means data resides on disk rather than in memory, typically because the dataset exceeds available RAM. This regime introduces challenges such as slower access, higher latency, I/O bottlenecks, and complex data management \cite{warren2015big}.} or distributed across machines. In such environments, communication cost often dominates—and can surpass—computational cost. Minimizing the number of passes over the tensor is therefore essential. Within randomized linear algebra, algorithms designed under this constraint are called \textit{pass-efficient}.

T-SVD generalizes the matrix SVD to third-order tensors via the tubal product (T-product) \cite{kilmer2013third}. It has been successfully applied to image and video completion \cite{zhang2016exact}, face recognition \cite{hao2013facial}, and compression \cite{zeng2020decompositions}. Randomized methods are now standard for computing low-tubal-rank approximations of large tensors \cite{zhang2018randomized,qi2021t,ahmadi2024randomized}. A pass-efficient randomized T-SVD algorithm was proposed in \cite{ahmadi2024randomized} requiring $v\geq 2$ passes. However, many real-world applications permit only a single pass ($v=1$), motivating the need for dedicated single-pass algorithms.

\begin{figure}[htbp]
\centering
\includegraphics[width=0.9\linewidth]{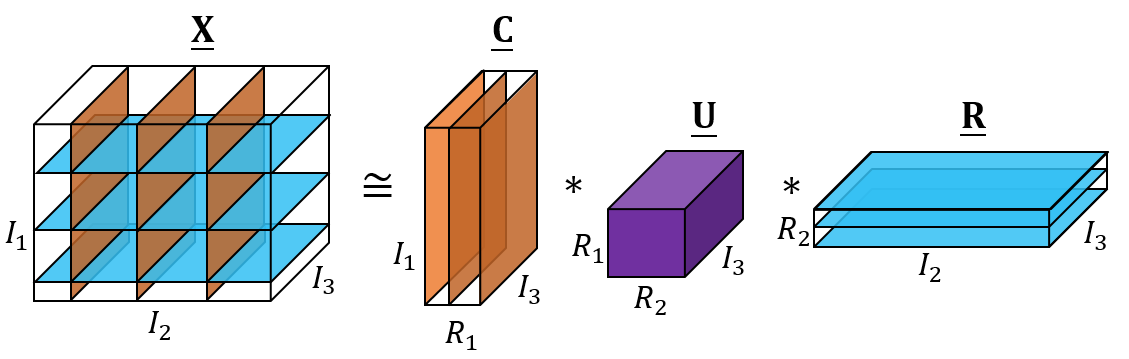}
\caption{Illustration of randomized low-tubal-rank approximation via slice sampling \cite{tarzanagh2018fast}. Lateral slices (brown) and frontal slices (blue) are sampled from the original tensor to construct a low-rank approximation using their intersection as the core tensor. This approach requires only a single pass over the data but suffers from numerical instability when the number of selected lateral and frontal slices is equal.}
\label{Pic2}
\end{figure}

To date, the only single-pass randomized T-SVD algorithms we are aware of are those in \cite{qi2021t}. While seminal, they have notable limitations: they do not exploit the fast T-product and T-QR decomposition from \cite{lu2019tensor}, and they are based on earlier matrix algorithms \cite{tropp2017practical} that have since been improved and stabilized \cite{bjarkason2019pass}. Another related line of work \cite{tarzanagh2018fast} constructs a CUR-like approximation by sampling lateral and frontal slices and using the pseudoinverse of their intersection as the core tensor (see Figure \ref{Pic2}). This approach suffers from poor conditioning and performs especially poorly when the number of selected frontal and lateral slices is equal—a weakness that we will demonstrate numerically.

To address these limitations, this paper proposes several novel robust single-pass algorithms for the T-SVD framework. Through extensive empirical evaluation against established baselines, we demonstrate that the proposed methods yield superior robustness and predictive performance. Furthermore, we showcase the practical utility of these algorithms in image super-resolution and object detection. To the best of our knowledge, this work represents the first application of single-pass randomized tensor decompositions to these computer vision tasks.

In addition, we study randomized fixed-precision algorithms for low-tubal-rank approximation. These are valuable when the tubal rank is unknown a priori: given an error tolerance, the algorithm must automatically select an appropriate rank and produce the corresponding approximation. Building on \cite{feng2023fast,ding2020efficient}, we refine the fixed-precision method proposed in \cite{ahmadi2023efficient}. Extensive simulations show that the refined algorithm is significantly faster and more computationally efficient.

\subsection*{Contributions}
\begin{itemize}
\item Novel efficient and robust randomized fixed-precision algorithms for T-SVD.
\item Three new single-pass algorithms for randomized T-SVD.
\item Theoretical guarantees for all proposed algorithms.
\item Extensive experiments on synthetic and real-world data, including applications to image/video compression, image super-resolution, and deep learning. This work is the first to apply single-pass tensor decomposition algorithms to tensor completion and image super-resolution.
\end{itemize}

\subsection*{Paper Organization}
In Section \ref{sec:related}, we discuss related work on randomized algorithms for T-SVD. Section \ref{Sec:prelim} recalls preliminary concepts and notation. Section \ref{Sec:tSVD} reviews T-SVD and existing algorithms. The new single-pass and fixed-precision algorithms are presented in Sections \ref{Sec:Propo} and \ref{sec:fixed}, respectively, together with theoretical analyses. Section \ref{Sec:Experi} reports experimental results, and the ablation study is discussed in Section \ref{sec:ablation}. Section \ref{Sec:Conclu} concludes the paper.

\section{Related Work}\label{sec:related}
In this section, we review the existing literature on randomized algorithms for tensor decompositions, with a particular focus on single-pass and fixed-precision methods. We identify the key limitations of current approaches and position our contributions within this landscape.

\subsection{Tensor Decompositions and the T-Product}
Tensor decompositions have become essential tools for multidimensional data analysis. The Tucker decomposition \cite{tucker1966some} and the CANDECOMP/PARAFAC (CP) decomposition \cite{carroll1970analysis} are among the most widely used models. However, these decompositions do not fully exploit the spectral structure of third-order tensors. The tensor-train (TT) decomposition \cite{oseledets2011tensor} offers a hierarchical representation but requires careful rank selection. More recently, the T-product framework introduced by Kilmer et al. \cite{kilmer2011factorization, kilmer2013third} provides an elegant algebraic foundation for third-order tensors, treating them as linear operators on matrices. This framework has enabled the development of tensor singular value decomposition (T-SVD) \cite{kilmer2013third}, which generalizes the matrix SVD while preserving many of its properties, including optimal low-rank approximation.

Despite its theoretical elegance, the computational cost of T-SVD is prohibitive for large-scale tensors, requiring $O(I_1 I_2 \min(I_1, I_2) I_3)$ operations. This limitation has motivated a rich line of research on randomized algorithms for T-SVD, which we review below.

\subsection{Randomized Algorithms for T-SVD}
Randomized algorithms have emerged as powerful tools for computing low-rank approximations of large matrices and tensors \cite{halko2011finding, martinsson2020randomized}. The core idea is to project the data onto a lower-dimensional subspace using random sketching matrices, dramatically reducing computational complexity while maintaining high accuracy with high probability.

\subsubsection{Multi-Pass Randomized T-SVD}
Zhang et al. \cite{zhang2018randomized} proposed the first randomized T-SVD algorithm (rT-SVD) based on the T-product framework. Their method computes a sketch $\Y = \X * \underline{\Omega}\in\mathbb{R}^{I_1\times K\times I_3}$, where $\underline{\Omega}\in\mathbb{R}^{I_2\times K\times I_3}$ is a random tensor, and $\underline{\mathbf X}\in\mathbb{R}^{I_1\times I_2\times I_3}$ then constructs an orthonormal basis $\underline{\mathbf{Q}}\in\mathbb{R}^{I_1\times K\times I_3}$ via tensor QR (T-QR) (to be discussed), and finally computes $\underline{\mathbf{X}} \approx \underline{\mathbf{Q}} * (\underline{\mathbf{Q}}^T * \underline{\mathbf{X}})$. This algorithm requires two passes over the data tensor. For a power iteration $q$, we need to compute $(\X*\X^T)*\X*\underline{\Omega}$, which needs $2q+2$ passes. this framework was extended in \cite{ahmadi2024randomized} to support an arbitrary number of passes $v \geq 2$, offering a flexible trade-off between accuracy and computational cost.

\subsubsection{Single-Pass Randomized T-SVD}
Single-pass algorithms are crucial when data cannot be stored in memory or revisited (e.g., streaming data). The first single-pass randomized T-SVD algorithms were introduced by Qi and Yu \cite{qi2021t}. Their method computes two sketches: $\Y = \X * \underline{\Omega}_1\in\mathbb{R}^{I_1\times K\times I_3}$ and $\W = \underline{\Omega}_2 * \X\in\mathbb{R}^{L \times I_2\times I_3}$, where $\underline{\Omega}_1\in\mathbb{R}^{I_2\times K\times I_3}$ and $\underline{\Omega}_2\in\mathbb{R}^{L\times I_1\times I_3}$. Then approximates $\X \approx \Q * (\underline{\Omega}_2 * \Q)^\dagger * \W$, where $\Q$ is an orthonormal basis for $\Y$. This approach generalizes the matrix single-pass algorithm of Tropp et al. \cite{tropp2017practical}.

Independently, Tarzanagh and Michailidis \cite{tarzanagh2018fast} proposed a cross tensor approximation (TCUR) method that samples lateral and horizontal slices and uses their intersection as the core tensor. While this method requires only one pass, it suffers from the limitation that the selected slices may not effectively capture the dominant singular subspaces. Similar TCUR methods were proposed in \cite{ahmadi2024robust}.

\subsection{Limitations of Existing Single-Pass Methods}
Despite their pioneering nature, existing single-pass T-SVD algorithms have several critical limitations that have not been adequately addressed in the literature:

\begin{enumerate}
    \item \textbf{Conditioning issues when $L = K$:} As shown in Bjarkason's work on matrix algorithms \cite{bjarkason2019pass}, single-pass methods become numerically unstable when the sketch sizes satisfy $L = K$ due to ill-conditioned least-squares problems. This issue has been largely ignored in the tensor literature. Both Qi and Yu \cite{qi2021t} and Tarzanagh and Michailidis \cite{tarzanagh2018fast} do not report empirical instability and do not provide a principled solution for this issue. The experiments confirm that when $L=K$, baseline methods produce PSNR values as low as $9.02$ dB (see Table~\ref{tab:kodak_comparison}).

    \item \textbf{Lack of regularization mechanisms:} Existing methods rely on T-QR decomposition for orthonormalization, which does not provide regularization. When sketches are ill-conditioned, the computed pseudoinverse $(\underline{\Omega}_2 * \Q)^\dagger$ becomes numerically unreliable.

    \item \textbf{No theoretical guarantees for the tensor case:} While matrix single-pass algorithms have rigorous error bounds \cite{tropp2017practical, bjarkason2019pass}, these guarantees have not been fully extended to the T-product tensor setting. The proofs in \cite{qi2021t} rely on slice-wise arguments that do not account for the coupling between frontal slices in the Fourier domain.

    \item \textbf{No adaptation for fixed-precision setting:} Existing single-pass methods require the target rank $R$ as an input. In many real-world applications, the optimal rank is unknown a priori, necessitating fixed-precision algorithms that automatically determine the rank from an error tolerance.
\end{enumerate}

\subsection{Recent Developments in Randomized Tensor Decompositions}
Several recent works have advanced randomized tensor computations, though none address the specific limitations identified above:

Che and Wei \cite{che2019randomized} proposed randomized algorithms for TT decomposition. However, these methods focus on higher-order tensors and do not leverage the T-product structure.

For fixed-precision low-rank approximation, Yu et al. \cite{yu2018efficient} proposed adaptive randomized algorithms for matrices that automatically determine the rank from an error tolerance. This idea was extended tensors to tensors via the T-product in \cite{ahmadi2023efficient}, introducing Algorithm \ref{ALgRR}. Feng and Yu \cite{feng2023fast} recently improved the matrix fixed-precision algorithm by replacing the QR decomposition with the LU decomposition and allowing an odd number of passes. However, their tensor extensions remain unexplored.

Zeng et al. \cite{zeng2025efficient} proposed randomized tensor sketch methods for online network security monitoring, demonstrating the practical value of pass-efficient algorithms. An incremental algorithm is proposed in \cite{abdelgawad2025inctsvd} for computing the TSVD of streaming tensors. Chen et al. \cite{chen2026rank} developed rank-revealing Bayesian tensor completion methods, highlighting the ongoing need for adaptive rank estimation.

\subsection{Gaps in the Literature and Our Contributions}
Based on the above review, we identify the following specific gaps that this work addresses:

\begin{table*}[htbp]
\centering
\small
\caption{Summary of limitations in existing single-pass T-SVD methods and how this work addresses each gap.}
\label{tab:gaps}
\begin{tabular}{p{3.2cm}p{3.8cm}p{5.5cm}}
\toprule
\textbf{Gap} & \textbf{Existing Methods} & \textbf{Our Contribution} \\
\midrule
Instability when $L = K$ in single-pass T-SVD & Qi and Yu \cite{qi2021t}, Tarzanagh and Michailidis \cite{tarzanagh2018fast} provide no solution & Algorithms \ref{Single-1}--\ref{Single-3} introduce truncation parameter $H$ as regularizer; Theorem 2 proves $\epsilon$-stability \\
\midrule
No theoretical guarantees for tensor single-pass methods & Matrix theory exists \cite{tropp2017practical,bjarkason2019pass} but not extended to T-product & Theorems 1--5 provide rigorous error bounds via Fourier domain decoupling \\
\midrule
Fixed-precision T-SVD requires even number of passes & Algorithm \ref{ALgRR} \cite{ahmadi2023efficient} requires $2q+2$ passes & Algorithm \ref{fixed-precision} supports arbitrary $q$ (odd or even), reducing pass count \\
\midrule
QR decomposition limits parallel efficiency & Algorithm \ref{ALgRR} \cite{ahmadi2023efficient} uses T-QR for orthonormalization & Algorithm \ref{fixed-precision} replaces T-QR with T-LU, achieving $25$--$30\%$ speedup \\
\midrule
QR decomposition limits parallel efficiency & Algorithm \ref{ALgRR} \cite{ahmadi2023efficient} uses T-QR for orthonormalization & Algorithm \ref{fixed-precis2} replaces T-QR with t-product and inversion of small tensors. \\
\midrule
No single-pass algorithms for tensor completion/super-resolution & Existing single-pass methods tested only on compression & First application to tensor completion and image super-resolution \\
\bottomrule
\end{tabular}
\end{table*}

\subsection{Positioning of Our Work}
This work builds upon but substantially extends the existing literature in the following ways:

\begin{enumerate}
    \item \textbf{Unlike \cite{qi2021t} and \cite{tarzanagh2018fast}:} We provide a principled regularization mechanism (truncated T-SVD with parameter $H$) that eliminates the instability when $L = K$. We address these contributions theoretically and numerically, with an ablation study quantifying the resulting improvements.
    
    \item \textbf{Unlike \cite{ahmadi2023efficient} (Algorithm \ref{ALgRR}):} We reduce the number of passes from $2q+2$ to any arbitrary $q$ (Algorithm \ref{ALgRR}), making the method more practical for streaming applications. We also replace T-QR with T-LU for better parallelization (Algorithm \ref{fixed-precision}).

    \item \textbf{Unlike all prior work:} We are the first to apply single-pass randomized T-SVD to tensor completion and image super-resolution, demonstrating practical utility beyond compression.

    \item \textbf{Theoretical contribution:} We provide the first error bounds for regularized single-pass T-SVD that explicitly account for the conditioning improvement via truncation.
\end{enumerate}

In summary, while our algorithms are inspired by existing matrix methods \cite{bjarkason2019pass, ding2020efficient, feng2023fast}, their extension to the T-product tensor setting is non-trivial due to the coupling between frontal slices in the Fourier domain. Our theoretical analysis, regularization strategy, and empirical validation collectively represent a significant advancement over the state of the art.

\section{Background}\label{Sec:prelim}
This section provides the necessary background on tensors and establishes the notation used throughout the paper. Tensors, matrices, and vectors are denoted by underlined bold uppercase letters (e.g., $\underline{\mathbf{X}}$), bold uppercase letters (e.g., $\mathbf{X}$), and bold lowercase letters (e.g., $\mathbf{x}$), respectively. For a third-order tensor $\X$, the slices $\X(:,:,k)$, $\X(:,j,:)$, and $\X(i,:,:)$ are referred to as frontal, lateral, and horizontal slices, respectively. We denote the $i$th frontal slice of $\X$ by $\X^{(i)}$. The conjugate transpose of a matrix $\mathbf{X}$ is denoted by ${\mathbf X}^*$. Additionally, the sub-tensor $\X(i,j,:)$ is called a tube. The Frobenius norm of a tensor is denoted by $\|\cdot\|_F$. The notation ``$\mathrm{conj}$'' indicates the complex conjugate of a complex number or the component-wise complex conjugate of a matrix. The ceiling function $\lceil n \rceil$ represents the smallest integer greater than or equal to $n$. While this paper focuses on real-valued tensors, all results can be straightforwardly extended to complex-valued tensors.

\begin{table}[htbp]
\centering
\caption{Notation used in Section 3 (preliminaries).}
\label{tab:notations}
\begin{tabular}{|c|c|}
\hline
\textbf{Notation} & \textbf{Description} \\
\hline
$\underline{\mathbf{X}}$ & Third-order tensor (underlined bold uppercase) \\
\hline
$\mathbf{X}$ & Matrix (bold uppercase) \\
\hline
$\mathbf{x}$ & Vector (bold lowercase) \\
\hline
$\X(:,:,k)$ & Frontal slices of tensor $\X$ \\
\hline
$\X(:,j,:)$ & Lateral slices of tensor $\X$ \\
\hline
$\X(i,:,:)$ & Horizontal slices of tensor $\X$ \\
\hline
$\X_1$ & First frontal slice of tensor $\X$ \\
\hline
$\X(i,j,:)$ & Tube (fiber along third dimension) of tensor $\X$ \\
\hline
$\|\cdot\|_F$ & Frobenius norm of a tensor \\
\hline
$\mathrm{conj}$ & Complex conjugate (component-wise for matrices) \\
\hline
$\lceil n \rceil$ & Ceiling function \\
\hline
$\boxplus_1$ & Concatenation along mode 1 (first dimension) \\
\hline
$\boxplus_2$ & Concatenation along mode 2 (second dimension) \\
\hline
$\begin{bmatrix} \A\\ \B \end{bmatrix}$ & Alternative notation for concatenation along mode 1 \\
\hline
$[\A,\B]$ & Alternative notation for concatenation along mode 2 \\
\hline
$\X * \Y$ & T-product (tubal product) of two tensors \\
\hline
$\mathrm{fold}(\cdot)$ & Operator reshaping unfolded tensor to original form \\
\hline
$\mathrm{unfold}(\cdot)$ & Operator stacking frontal slices vertically \\
\hline
$\X^T$ & Transpose of tensor $\X$ \\
\hline
$\I$ & Identity tensor \\
\hline
$\X^{\dag}$ & Moore–Penrose pseudoinverse of tensor $\X$ \\
\hline
$\X^{-1}$ & Inverse of tensor $\X$ \\
\hline
$\underline{\Omega}$ & Random tensor (first frontal slice Gaussian, others zero) \\
\hline
\end{tabular}
\end{table}

Two tensors can be concatenated along the first or second modes. A concatenation along mode 1 of tensors $\underline{\mathbf{A}}\in\mathbb{R}^{I_1\times I_2\times I_3}$ and $\underline{\mathbf{B}}\in\mathbb{R}^{J_1\times J_2\times J_3}$ is denoted as $\underline{\mathbf{C}} = \underline{\mathbf{A}} \boxplus_1 \underline{\mathbf{B}} \in {\mathbb{R}^{(I_1+J_1)\times I_2\times I_3}},$ where $I_2=J_2,\,I_3=J_3$, and the same definition can be stated for concatenation along the second mode. Alternative notations for the concatenation along the first and second modes are $\A\boxplus_1 \B \equiv \begin{bmatrix} \A\\ \B \end{bmatrix}$ and $\A\boxplus_2 \B \equiv [\A,\B]$.

\begin{Definition}[Tube and Circular Convolution]
Let $\underline{\mathbf{X}} \in \mathbb{R}^{I_1 \times I_2 \times I_3}$ be a third-order tensor. The \emph{tubes} of $\underline{\mathbf{X}}$ are the vectors $\underline{\mathbf{X}}(i,j,:) \in \mathbb{R}^{I_3}$ along the third dimension. For two tubes $\mathbf{a}, \mathbf{b} \in \mathbb{R}^{I_3}$, their \emph{circular convolution} is the tube $\mathbf{c} = \mathbf{a} * \mathbf{b} \in \mathbb{R}^{I_3}$ defined componentwise as
\[
\mathbf{c}(k) = \sum_{\ell=1}^{I_3} \mathbf{a}(\ell) \, \mathbf{b}(k - \ell + 1 \bmod I_3),
\]
where indices are taken modulo $I_3$ (with $1$ replacing $I_3+1$).
\end{Definition}

\begin{Definition}[T-Product]
Let $\underline{\mathbf{X}} \in \mathbb{R}^{I_1 \times I_2 \times I_3}$ and $\underline{\mathbf{Y}} \in \mathbb{R}^{I_2 \times I_4 \times I_3}$. The T-product $\underline{\mathbf{C}} = \underline{\mathbf{X}} * \underline{\mathbf{Y}} \in \mathbb{R}^{I_1 \times I_4 \times I_3}$ is defined such that for each $(i,j)$, the $(i,j)$-th tube of $\underline{\mathbf{C}}$ is the sum of circular convolutions of tubes from $\underline{\mathbf{X}}$ and $\underline{\mathbf{Y}}$:
\[
\underline{\mathbf{C}}(i,j,:) = \sum_{k=1}^{I_2} \underline{\mathbf{X}}(i,k,:) * \underline{\mathbf{Y}}(k,j,:),
\]
where $*$ denotes circular convolution along the third dimension. Equivalently, using the block circulant matrix formulation,
\[
\underline{\mathbf{C}} = \underline{\mathbf{X}} * \underline{\mathbf{Y}} = \mathrm{fold}\left( \mathrm{bcirc}\left( \underline{\mathbf{X}} \right) \mathrm{unfold}\left( \underline{\mathbf{Y}} \right) \right),
\]
where
\[
\mathrm{bcirc} \left(\underline{\mathbf{X}}\right)
=
\begin{bmatrix}
\underline{\mathbf{X}}(:,:,1) & \underline{\mathbf{X}}(:,:,I_3) & \cdots & \underline{\mathbf{X}}(:,:,2)\\
\underline{\mathbf{X}}(:,:,2) & \underline{\mathbf{X}}(:,:,1) & \cdots & \underline{\mathbf{X}}(:,:,3)\\
 \vdots & \vdots & \ddots &  \vdots \\
 \underline{\mathbf{X}}(:,:,I_3) & \underline{\mathbf{X}}(:,:,I_3-1) & \cdots & \underline{\mathbf{X}}(:,:,1)
\end{bmatrix},
\]
and
\[
\mathrm{unfold}(\underline{\mathbf{Y}})=
\begin{bmatrix}
\underline{\mathbf{Y}}(:,:,1)\\
\underline{\mathbf{Y}}(:,:,2)\\
\vdots\\
\underline{\mathbf{Y}}(:,:,I_3)
\end{bmatrix},\quad
\underline{\mathbf{Y}}=\mathrm{fold} \left(\mathrm{unfold}\left(\underline{\mathbf{Y}}\right)\right).
\]
\end{Definition}

\begin{Remark}
The equivalence between the tube-wise circular convolution definition and the block circulant formulation follows from the fact that block circulant matrices diagonalize under the Fourier transform. Indeed, for $\widehat{\underline{\mathbf{X}}} = \mathrm{fft}(\underline{\mathbf{X}},[],3)$, we have
\[
\widehat{\underline{\mathbf{C}}}(:,:,k) = \widehat{\underline{\mathbf{X}}}(:,:,k) \, \widehat{\underline{\mathbf{Y}}}(:,:,k), \quad k=1,2,\dots,I_3,
\]
so the T-product reduces to standard matrix multiplication in the Fourier domain. Algorithm \ref{ALG:Tp} summarizes the computational process of the T-product.
\end{Remark}

It can be proved that for a tensor $\underline{\mathbf{X}} \in \mathbb{R}^{I_1 \times I_2 \times I_3}$, we have
\begin{eqnarray}\label{eq_fou}
\|\underline{\mathbf{X}}\|^2_F=\frac{1}{I_3}\sum_{i=1}^{I_3}\|\widehat{\underline{\mathbf{X}}}(:,:,i)\|_F^2,
\end{eqnarray}
where $\widehat{\underline{\mathbf{X}}}(:,:,i)$ is the $i$-th frontal slice of the tensor $\widehat{\underline{\mathbf{X}}}=\mathrm{fft}(\underline{\mathbf{X}},[],3)$, which calculates the fast Fourier transform of all tubes of the tensor $\underline{\mathbf{X}}$ \cite{lu2019tensor,zhang2018randomized}.  Algorithm \ref{ALG:Tp} summarizes the process of the T-product of two tensors.

\begin{Definition} (Transpose)
Let $\underline{\mathbf{X}}\in\mathbb{R}^{I_1\times I_2\times I_3}$ be a given tensor. Then the transpose of the tensor $\underline{\mathbf{X}}$ is denoted by $\underline{\mathbf{X}}^{T}\in\mathbb{R}^{I_2\times I_1\times I_3}$, which is constructed by transposing all its frontal slices and then reversing the order of transposed frontal slices $2$ through $I_3$. A tensor $\underline{\mathbf{X}}\in\mathbb{R}^{I\times I\times K}$ is called symmetric if $\underline{\mathbf{X}}^T=\underline{\mathbf{X}}$.
\end{Definition}

\begin{Definition} (Identity tensor)
The tensor $\underline{\mathbf{I}}\in\mathbb{R}^{I_1\times I_1\times I_3}$ is called identity if its first frontal slice is an identity matrix of size $I_1\times I_1$ and all other frontal slices are zero. It is easy to show ${\I}*{\X}={ \X}$ and ${\X}*{\I} ={\X}$ for all tensors of conforming sizes.
\end{Definition}

\begin{Definition} (Orthogonal tensor)
A tensor $\underline{\mathbf{X}}\in\mathbb{R}^{I_1\times I_1\times I_3}$ is orthogonal if ${\underline{\mathbf{X}}^T} * \underline{\mathbf{X}} = \underline{\mathbf{X}} * {\underline{\mathbf{X}}^ T} = \underline{\mathbf{I}}$.
\end{Definition}

\begin{Definition} (Moore–Penrose pseudoinverse of a tensor)
Let $\X\in\mathbb{R}^{I_1\times I_2\times I_3}$ be given. The Moore-Penrose (MP) pseudoinverse of the tensor $\X$ is denoted by $\X^{\dag}\in\mathbb{R}^{I_2\times I_1\times I_3}$ and is a unique tensor satisfying the following four equations:
\begin{eqnarray*}
\X^{\dag}*\X*\X^{\dag}=\X^{\dag},\quad \X*\X^{\dag}*\X=\X,\\
(\X*\X^{\dag})^T=\X*\X^{\dag},\quad (\X^{\dag}*\X)^T=\X^{\dag}*\X.
\end{eqnarray*}
The MP pseudoinverse of a tensor can also be computed in the Fourier domain, and this is shown in Algorithm \ref{ALG:TMP}. The inverse of a tensor $\X\in\mathbb{R}^{I_1\times I_1\times I_3}$ is denoted by $\X^{-1}\in\mathbb{R}^{I_1\times I_1\times I_3}$ and is a special case of the MP for which we have ${\underline{\mathbf{X}}^{-1}} * \underline{\mathbf{X}} = \underline{\mathbf{X}} * {\underline{\mathbf{X}}^{-1}} = \underline{\mathbf{I}}$.
\end{Definition}

\begin{Definition} (f-diagonal tensor)
If all frontal slices of a tensor are diagonal, then the tensor is called an f-diagonal tensor.
\end{Definition}

\begin{Definition} (Random tensor)
A tensor $\underline{\Omega}$ is random if its first frontal slice $\underline{\Omega}(:,:,1)$ is a standard Gaussian matrix, while the other frontal slices are zero.
\end{Definition}

\begin{Lemma}
Let $\underline{\mathbf{B}} \in \mathbb{R}^{n \times n \times n_3}$ be a third-order tensor, and let $\underline{\mathbf{B}}^T$ denote its transpose under the t-product. Define $\underline{\mathbf{G}} = \underline{\mathbf{B}}^T * \underline{\mathbf{B}}$, where $*$ denotes the t-product (circular convolution along the third dimension). Then the following identity holds:
\[
\|\underline{\mathbf{B}}\|_F^2 = \operatorname{trace}\bigl(\underline{\mathbf{G}}(:,:,1)\bigr),
\]
where $\underline{\mathbf{G}}(:,:,1)$ is the first frontal slice of $\underline{\mathbf{G}}$.
\end{Lemma}

\begin{proof}
For clarity, denote the entries of $\underline{\mathbf{B}}$ by $\underline{\mathbf{B}}(i,j,k)$, where $i,j \in \{1,2,\dots,n\}$ and $k \in \{1,2,\dots,n_3\}$. By definition of the t-product, for any indices $p,q \in \{1,2,\dots,n\}$, the tube $\underline{\mathbf{G}}(p,q,:)$ is given by the circular convolution
\[
\underline{\mathbf{G}}(p,q,:) = \sum_{r=1}^{n} \underline{\mathbf{B}}^T(p,r,:) * \underline{\mathbf{B}}(r,q,:),
\]
where $*$ denotes circular convolution along the third dimension. Equivalently, in terms of entries, for each $k \in \{1,2,\dots,n_3\}$,
\[
\underline{\mathbf{G}}(p,q,k) = \sum_{r=1}^{n} \sum_{\ell=1}^{n_3} \underline{\mathbf{B}}^T(p,r,\ell) \, \underline{\mathbf{B}}(r,q, k - \ell + 1 \bmod n_3),
\]
with indices taken modulo $n_3$. Now take the trace of the first frontal slice. By definition,
\[
\operatorname{trace}\bigl(\underline{\mathbf{G}}(:,:,1)\bigr) = \sum_{i=1}^{n} \underline{\mathbf{G}}(i,i,1).
\]
Using the entrywise formula for the t-product at index $k=1$, we have
\[
\underline{\mathbf{G}}(i,i,1) = \sum_{r=1}^{n} \sum_{\ell=1}^{n_3} \underline{\mathbf{B}}^T(i,r,\ell) \, \underline{\mathbf{B}}(r,i, 1 - \ell + 1 \bmod n_3).
\]
Since $\underline{\mathbf{B}}^T$ is the tensor transpose, its entries satisfy
\[
\underline{\mathbf{B}}^T(i,r,\ell) = \underline{\mathbf{B}}(r,i,\ell)
\]
(up to conjugation; over $\mathbb{R}$ this is simply transposition). Also, the circular convolution at index $1$ effectively sums over all $\ell$ because the shift $1 - \ell + 1 \bmod n_3$ runs through all indices as $\ell$ varies. Thus,
\[
\underline{\mathbf{G}}(i,i,1) = \sum_{r=1}^{n} \sum_{\ell=1}^{n_3} \underline{\mathbf{B}}(r,i,\ell) \, \underline{\mathbf{B}}(r,i,\ell)
= \sum_{r=1}^{n} \sum_{\ell=1}^{n_3} |\underline{\mathbf{B}}(r,i,\ell)|^2.
\]
Summing this equality over $i = 1,2, \dots, n$ gives
\[
\operatorname{trace}\bigl(\underline{\mathbf{G}}(:,:,1)\bigr)
= \sum_{i=1}^{n} \sum_{r=1}^{n} \sum_{\ell=1}^{n_3} |\underline{\mathbf{B}}(r,i,\ell)|^2.
\]
But the right-hand side is precisely the definition of the squared Frobenius norm:
\[
\|\underline{\mathbf{B}}\|_F^2 = \sum_{i=1}^{n} \sum_{j=1}^{n} \sum_{k=1}^{n_3} |\underline{\mathbf{B}}(i,j,k)|^2.
\]
Renaming the summation indices (interchanging $i$ and $r$) confirms that
\[
\operatorname{trace}\bigl(\underline{\mathbf{G}}(:,:,1)\bigr) = \|\underline{\mathbf{B}}\|_F^2.
\]
Hence the lemma is proved.
\end{proof}

\begin{algorithm}
\caption{T-product in the Fourier domain}\label{ALG:Tp}
\textbf{Input:} Two data tensors $\underline{\mathbf{X}} \in {\mathbb{R}^{{I_1} \times {I_2} \times {I_3}}}$ and $\underline{\mathbf{Y}} \in {\mathbb{R}^{{I_2} \times {I_4} \times {I_3}}}$;\\
\textbf{Output:} The T-product $\underline{\mathbf{C}} \in {\mathbb{R}^{{I_1} \times {I_4} \times {I_3}}}$;
\begin{algorithmic}[1]
\State $\widehat{\underline{\mathbf{X}}} = \mathrm{fft}\left( {\underline{\mathbf{X}},[],3} \right)$;
\State $\widehat{\underline{\mathbf{Y}}} = \mathrm{fft}\left( {\underline{\mathbf{Y}},[],3} \right)$;
\For{$i=1,2,\ldots,\lceil \frac{I_3+1}{2}\rceil$}
\State $\widehat{\underline{\mathbf{C}}}\left( {:,:,i} \right) = \widehat{\underline{\mathbf{X}}}\left( {:,:,i} \right)\,\widehat{\underline{\mathbf{Y}}}\left( {:,:,i} \right)$;
\EndFor
\For{$i=\lceil\frac{I_3+1}{2}\rceil+1,\ldots,I_3$}
\State $\widehat{\underline{\mathbf{C}}}\left( {:,:,i} \right)=\mathrm{conj}(\widehat{\underline{\mathbf{C}}}\left( {:,:,I_3-i+2} \right))$;
\EndFor
\State $\underline{\mathbf{C}} = \mathrm{ifft}\left( {\widehat{\underline{\mathbf{C}}},[],3} \right)$;
\end{algorithmic}
\end{algorithm}

\begin{algorithm}
\caption{Fast Moore-Penrose pseudoinverse computation of the tensor $\underline{\mathbf{X}}$}\label{ALG:TMP}
\textbf{Input:} A data tensor $\underline{\mathbf{X}} \in {\mathbb{R}^{{I_1} \times {I_2} \times {I_3}}}$;\\
\textbf{Output:} The MP pseudoinverse $\underline{\mathbf{X}}^{\dagger} \in {\mathbb{R}^{{I_2} \times {I_1} \times {I_3}}}$;
\begin{algorithmic}[1]
\State $\widehat{\underline{\mathbf{X}}} = \mathrm{fft}\left( {\underline{\mathbf{X}},[],3} \right)$;
\For{$i=1,2,\ldots,\lceil \frac{I_3+1}{2}\rceil$}
\State $\widehat{\underline{\mathbf{C}}}\left( {:,:,i} \right) = \mathrm{pinv}\,(\widehat{\underline{\mathbf{X}}}(:,:,i))$;
\EndFor
\For{$i=\lceil\frac{I_3+1}{2}\rceil+1,\ldots,I_3$}
\State $\widehat{\underline{\mathbf{C}}}\left( {:,:,i} \right)=\mathrm{conj}(\widehat{\underline{\mathbf{C}}}\left( {:,:,I_3-i+2} \right))$;
\EndFor
\State $\underline{\mathbf{X}}^{\dag} = \mathrm{ifft}\left( {\widehat{\underline{\mathbf{C}}},[],3} \right)$;
\end{algorithmic}
\end{algorithm}

\section{Generalization of Standard Matrix Decompositions to Tensors through the T-Product}\label{Sec:tSVD}
The T-product provides a straightforward way to generalize standard matrix decompositions---such as QR, LU, eigenvalue decomposition, and SVD---to tensors. For a tensor $\X\in\mathbb{R}^{I_1\times I_2\times I_3}$, the tensor QR (T-QR) decomposition expresses $\X$ as $\X=\Q*\R$ and can be computed using Algorithm \ref{ALG:TQR}. In line 3 of the algorithm, the MATLAB command \texttt{qr(X,0)} performs the economic QR decomposition of the matrix $\mathbf{X}$. The operator ``orth'' provides an orthonormal basis for a tensor by means of the T-QR decomposition. That is, $\Q = \text{orth}(\X)$ signifies $[\Q, \sim] = \text{T-QR}(\X)$. In the special case where $\X$ is a matrix, this operator coincides with the standard QR-based orthonormalization of its columns.

Tensor LU (T-LU), tensor eigenvalue decomposition (T-EIG), and tensor SVD (T-SVD) can be obtained by slightly modifying Algorithm \ref{ALG:TQR}. Specifically, we replace the QR decomposition in line 4 of Algorithm \ref{ALG:TQR} with the SVD, LU decomposition, or eigenvalue decomposition of the frontal slices $\widehat{\X}(:,:,i)$ for $i=1,2,\ldots,I_3$. Note that Algorithm \ref{ALG:TQR} only requires the thin QR of the first $\lceil (I_3+1)/2\rceil$ slices, whereas the original T-QR algorithm introduced in \cite{kilmer2013third} computes the QR of all frontal slices. This approach is recommended to eliminate redundant computations. For example, the T-EIG decomposition of a symmetric tensor $\X\in\mathbb{R}^{I_1\times I_1\times I_3}$ guarantees the following decomposition:
\begin{eqnarray}\label{t_eig}
    \X=\V*\D*\V^T,
\end{eqnarray}
where $\D$ is f-diagonal. The eigenvalue decomposition \eqref{t_eig} can also be expressed as:
\begin{eqnarray}
     \X=(\V*\D^{1/2})*(\D^{1/2}*\V^T),
\end{eqnarray}
where $\D^{1/2}\equiv\mathrm{sqrt}(\D)$ is an f-diagonal tensor whose tubes are computed by performing the following operations on each tube of $\D$, denoted $\D(i,i,:)$:

\begin{algorithm}
\begin{algorithmic}[1]
\State $\widehat{\D}(i,i,:)=\mathrm{fft}(\D(i,i,:),[],3)$;
\For{$j=1:I_3$}
\State $\widehat{\D}(i,i,j)=\mathrm{sqrt}(\widehat{\D}(i,i,j))$;
\EndFor
\State $\D^{1/2}(i,i,:)=\mathrm{ifft}( \widehat{\D}(i,i,:),[],3)$;
\end{algorithmic}
\end{algorithm}

Notably, the T-SVD provides a useful decomposition by representing a tensor as a T-product of three tensors. The first and last tensors are orthogonal, while the middle tensor is f-diagonal. Let $\underline{\mathbf{X}}\in\mathbb{R}^{I_1\times I_2\times I_3}$ be a tensor; then the T-SVD of this tensor is given by:
\begin{eqnarray}
\underline{\mathbf{X}}=\underline{\mathbf{U}} * \underline{\mathbf{S}}* \underline{\mathbf{V}}^T,
\end{eqnarray}
where $\underline{\mathbf{U}}\in\mathbb{R}^{I_1\times I_1\times I_3}$ and $\underline{\mathbf{V}}\in\mathbb{R}^{I_2\times I_2\times I_3}$ are orthogonal, and $\underline{\mathbf{S}}\in\mathbb{R}^{I_1\times I_2\times I_3}$ is f-diagonal. Note that, unlike the matrix SVD where the middle matrix is diagonal and symmetric, the middle $f$-diagonal tensor $\underline{\mathbf{S}}$ in the tensor SVD is not necessarily symmetric. Here, $\underline{\mathbf{U}}$ and $\underline{\mathbf{V}}$ are called the left and right singular tensors. The number of non-zero diagonal tubes of $\underline{\mathbf{S}}\in\mathbb{R}^{I_1\times I_2\times I_3}$ is called the tubal rank. The formal definition of the tubal rank is as follows.

\begin{Definition}\label{def:tubal_rank}
The tubal rank of a tensor $\X \in \mathbb{R}^{I_1 \times I_2 \times I_3}$ is defined as the maximum rank among all frontal slices in the Fourier domain:
\begin{equation}
\mathrm{rank}_t(\X) = \max_{i=1,\ldots,I_3} \mathrm{rank}(\widehat{\X}(:,:,i)).
\end{equation}
\end{Definition}

\begin{figure*}[htbp]
\begin{center}
\includegraphics[width=9.5cm,height=4.8cm]{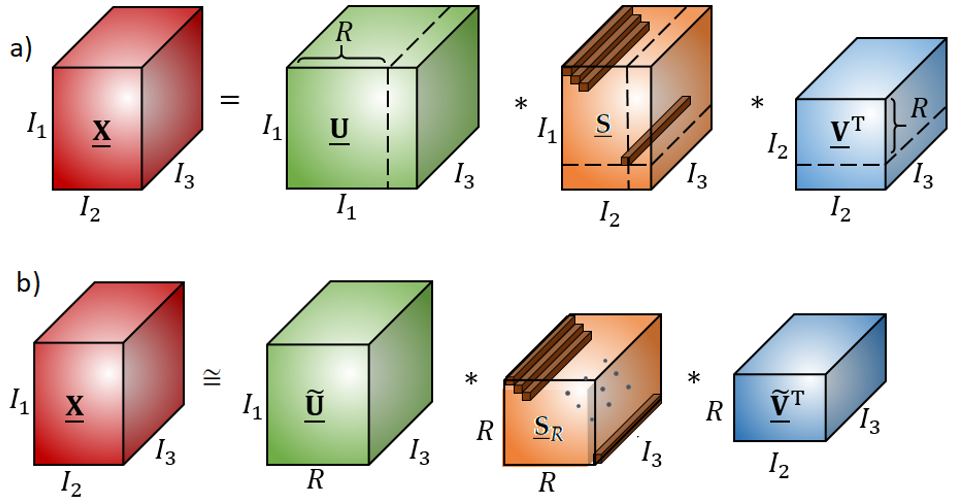}
\caption{Structure of tensor SVD (T-SVD) and its truncated version for a third-order tensor. (a) Full T-SVD: $\underline{\mathbf{X}} = \underline{\mathbf{U}} * \underline{\mathbf{S}} * \underline{\mathbf{V}}^T$ with orthogonal $\underline{\mathbf{U}}$, $\underline{\mathbf{V}}$ and f-diagonal $\underline{\mathbf{S}}$. (b) Truncated T-SVD of tubal rank $R$: only the first $R$ tubes of $\underline{\mathbf{S}}$ and corresponding columns of $\underline{\mathbf{U}}$, $\underline{\mathbf{V}}$ are retained. The tubal rank is the maximum rank of all frontal slices in the Fourier domain.}
\label{Pic7}
\end{center}
\end{figure*}

Thus, in the Fourier domain, the T-SVD decouples into $I_3$ independent matrix SVD problems. The truncated T-SVD is defined by truncating the tensors $\underline{\mathbf{U}}$, $\underline{\mathbf{S}}$, and $\underline{\mathbf{V}}$. For example, a truncated T-SVD of tubal rank $R$ for the tensor $\underline{\mathbf{X}}$ is:
\begin{eqnarray}
\underline{\mathbf{X}}\approx\underline{\mathbf{U}}_R * \underline{\mathbf{S}}_R* \underline{\mathbf{V}}_R^T,
\end{eqnarray}
where $\underline{\mathbf{U}}_R=\underline{\mathbf{U}}(:,1:R,:)\in\mathbb{R}^{I_1\times R\times I_3}$, $\underline{\mathbf{V}}_R=\underline{\mathbf{V}}(:,1:R,:)\in\mathbb{R}^{R\times I_2\times I_3}$, and $\underline{\mathbf{S}}_R=\underline{\mathbf{S}}(1:R,1:R,:)\in\mathbb{R}^{R \times R\times I_3}$. The tensors $\underline{\mathbf{U}}_R$ and $\underline{\mathbf{V}}_R$ are the left and right leading singular tensors. Figure \ref{Pic7} illustrates the structure of this decomposition and its truncated version.

Similar to the economic SVD, the economic T-SVD is defined. More precisely, the economic T-SVD yields $\X = \U * \underline{\bf S} * \V^T$. When $I_1 \leq I_2$ and $\X$ is of full tubal rank, then $\U \in \mathbb{R}^{I_1 \times I_1 \times I_3}$ and $\V \in \mathbb{R}^{I_2 \times I_1 \times I_3}$; when $I_1 \geq I_2$, then $\U \in \mathbb{R}^{I_1 \times I_2 \times I_3}$ and $\V \in \mathbb{R}^{I_2 \times I_2 \times I_3}$.

\begin{algorithm}
\caption{Fast T-QR decomposition of the tensor $\underline{\mathbf{X}}$}\label{ALG:TQR}
\textbf{Input:} The data tensor $\underline{\mathbf{X}} \in {\mathbb{R}^{{I_1} \times {I_2} \times {I_3}}}$;\\
\textbf{Output:} The T-QR decomposition $\underline{\mathbf{X}}=\underline{\mathbf{Q}}*\underline{\mathbf{R}}$;
\begin{algorithmic}[1]
\State $\widehat{\underline{\mathbf{X}}} = \mathrm{fft}\left( {\underline{\mathbf{X}},[],3} \right)$;
\For{$i=1,2,\ldots,\lceil \frac{I_3+1}{2}\rceil$}
\State $[\widehat{\underline{\mathbf{Q}}}\left( {:,:,i} \right),\widehat{\underline{\mathbf{R}}}(:,:,i)] = \mathrm{qr}\,(\widehat{\underline{\mathbf{X}}}(:,:,i),0)$;
\EndFor
\For{$i=\lceil\frac{I_3+1}{2}\rceil+1,\ldots,I_3$}
\State $\widehat{\underline{\mathbf{Q}}}\left( {:,:,i} \right)=\mathrm{conj}(\widehat{\underline{\mathbf{Q}}}\left( {:,:,I_3-i+2} \right))$;
\State $\widehat{\underline{\mathbf{R}}}\left( {:,:,i} \right)=\mathrm{conj}(\widehat{\underline{\mathbf{R}}}\left( {:,:,I_3-i+2} \right))$;
\EndFor
\State $\underline{\mathbf{Q}}= \mathrm{ifft}\left( {\widehat{\underline{\mathbf{Q}}},[],3} \right)$;
\State $\underline{\mathbf{R}}= \mathrm{ifft}\left( {\widehat{\underline{\mathbf{R}}},[],3} \right)$;
\end{algorithmic}
\end{algorithm}

\section{Proposed Randomized Single-Pass Algorithms}\label{Sec:Propo}
Randomized single-pass algorithms process a stream of data in a single pass and do not need to store the entire dataset in memory, which is crucial for large datasets that cannot fit in memory. They play a vital role in processing massive datasets and providing efficient solutions in various applications. Their ability to process data in a single pass, reduce memory usage, and handle worst-case scenarios makes them essential for real-time data analysis, data streaming, and machine learning tasks. For these reasons, developing single-pass and, more generally, pass-efficient algorithms for fast low-rank tensor approximation has been a hot topic over the last decade.

Let us first briefly introduce the standard randomized algorithm with oversampling and power iteration for computing a low tubal rank approximation of a tensor, as it will be needed later in this paper. Algorithm \ref{ALG:RandomizedTSVD} outlines the basic randomized T-SVD with power iteration for computing a low-rank approximation of a third-order tensor $\underline{\mathbf{X}} \in \mathbb{R}^{I_1 \times I_2 \times I_3}$. The algorithm first generates a random tensor $\underline{\boldsymbol{\Omega}}$ and forms the sample tensor $\underline{\mathbf{X}} * \underline{\boldsymbol{\Omega}}$, whose range approximates that of $\underline{\mathbf{X}}$. An orthonormal basis $\underline{\mathbf{Q}}$ for this sample tensor is then computed via the T-QR decomposition. To improve the accuracy of the approximation, the power iteration scheme is applied $p$ times, alternating between the tensors $\underline{\mathbf{X}}^T$ and $\underline{\mathbf{X}}$, which enhances the decay of the singular values. Subsequently, the reduced tensor $\underline{\mathbf{B}} = \underline{\mathbf{Q}}^T * \underline{\mathbf{X}}$ is formed, and its economic T-SVD is computed. Finally, the rank-$k$ factors $\underline{\mathbf{U}}$, $\underline{\mathbf{S}}$, and $\underline{\mathbf{V}}$ are obtained by truncating the resulting T-SVD components to the first $k$ frontal slices.

\begin{algorithm}
\caption{Basic Randomized T-SVD with Power Iteration}\label{ALG:RandomizedTSVD}
\textbf{Input:} $\underline{\mathbf{X}} \in \mathbb{R}^{I_1 \times I_2 \times I_3}$, rank parameter $k$, oversampling parameter $s$, power parameter $p$\\
\textbf{Output:} $\underline{\mathbf{U}} \in \mathbb{R}^{I_1 \times k \times I_3}$, $\underline{\mathbf{S}} \in \mathbb{R}^{k \times k \times I_3}$, $\underline{\mathbf{V}} \in \mathbb{R}^{I_2 \times k \times I_3}$
\begin{algorithmic}[1]
\State $l \leftarrow k + s$, $\underline{\boldsymbol{\Omega}} \leftarrow \text{randn}(I_2, l, I_3)$;
\State $\underline{\mathbf{Q}} \leftarrow \text{orth}(\underline{\mathbf{X}} * \underline{\boldsymbol{\Omega}})$; \Comment{$\underline{\mathbf{Q}} \in \mathbb{R}^{I_1 \times l \times I_3}$};
\For{$j \leftarrow 1, 2, \ldots, p$};
\State $\underline{\mathbf{G}} \leftarrow \text{orth}(\underline{\mathbf{X}}^T * \underline{\mathbf{Q}})$; \Comment{$\underline{\mathbf{G}} \in \mathbb{R}^{I_2 \times l \times I_3}$};
\State $\underline{\mathbf{Q}} \leftarrow \text{orth}(\underline{\mathbf{X}} * \underline{\mathbf{G}})$; \Comment{$\underline{\mathbf{Q}} \in \mathbb{R}^{I_1 \times l \times I_3}$};
\EndFor
\State $\underline{\mathbf{B}} \leftarrow \underline{\mathbf{Q}}^T * \underline{\mathbf{XA}}$; \Comment{$\underline{\mathbf{B}} \in \mathbb{R}^{l \times I_2 \times I_3}$}
\State $[\underline{\mathbf{U}}, \underline{\mathbf{S}}, \underline{\mathbf{V}}] \leftarrow \text{Economic T-SVD}(\underline{\mathbf{B}})$;
\State $\underline{\mathbf{U}} \leftarrow \underline{\mathbf{Q}} * \underline{\mathbf{U}}(:, 1:k, :)$;
\State $\underline{\mathbf{S}} \leftarrow \underline{\mathbf{S}}(1:k, 1:k, :)$;
\State $\underline{\mathbf{V}} \leftarrow \underline{\mathbf{V}}(:, 1:k, :)$;
\end{algorithmic}
\end{algorithm}

\begin{Definition}\label{def:stable_approximation}
A randomized algorithm provides an $\epsilon$-stable approximation to $\X$ if, with probability at least $1-\delta$, the output $\widetilde{\X}$ satisfies:
\begin{equation}
\|\X - \widetilde{\X}\|_F \leq (1 + \epsilon)\|\X - \X_R\|_F,
\end{equation}
where $\X_R$ is the best tubal rank-$R$ approximation of $\X$, and $\delta$ is a failure probability that depends on $\epsilon$ and other parameters of the algorithm.
\end{Definition}

\subsection{Existing Single-Pass Methods and Their Limitations}
In this paper, we focus on the T-SVD and propose new single-pass algorithms. To the best of our knowledge, the only single-pass algorithms proposed for low tubal rank approximation of tensors were proposed in \cite{tarzanagh2018fast} and \cite{qi2021t}. Let us discuss these single-pass algorithms below:

\subsection{Cross Tensor Approximation}
The proposed approach in \cite{tarzanagh2018fast} is based on cross tensor approximation. Specifically, for a given tensor $\X\in\mathbb{R}^{I_1\times I_2\times I_3}$, a low tubal rank approximation is computed by sampling some lateral and horizontal slices of the original data tensor $\X$, as illustrated in Figure~\ref{Pic2}. Let $\C\in\mathbb{R}^{I_1\times L\times I_3}$ and $\R\in\mathbb{R}^{K\times I_2\times I_3}$ be selected lateral and horizontal slices\footnote{Note that the lateral and frontal slices do not necessarily need to have the same size.} of the data tensor $\X$. Then, a low tubal rank approximation is computed as follows:
\begin{eqnarray}
\X \approx \C * \U * \R,
\end{eqnarray}
where $\U = \C^\dag * \X * \R^\dag$. The middle tensor $\U$ is the optimal candidate as it minimizes the residual $\|\X - \C * \U * \R\|_F$. This is known as cross tensor approximation or tensor CUR (TCUR) approximation. The tensor CY (TCY) and tensor YR (TYR) approximations are defined as $\X \approx \C * \Y$ and $\X \approx \Y * \R$ for given lateral and horizontal slices $\C$ and $\R$, respectively.

Notably, this approach requires only one pass over the data to compute the middle tensor $\U$. Our experimental results confirm that this method is stable with respect to the number of selected lateral and horizontal slices. However, it has three main drawbacks. 

First, its accuracy is generally lower than that of other techniques because the factor tensors $\C$ and $\R$ are not very accurate estimates of the left and right singular tensors. Second, it typically requires a higher tubal rank to achieve the desired accuracy. Third, the computation of $\C^\dag * \X * \R^\dag$ is more expensive than multiplication with random tensors, as one can use structured random tensors to accelerate the computations. 

In \cite{tarzanagh2018fast}, it is suggested to use $\U$ as the intersection tensor of some selected lateral and horizontal slices, as shown in Figure~\ref{Pic2}. This modified method is quite fast but suffers from the same instability problem when $L = K$, i.e., when the number of lateral slices equals the number of horizontal slices. This algorithm is summarized in Algorithm~\ref{crosstensor}. The following theorem establishes an upper bound on the TCUR approximation error in terms of the optimal CX and YR approximations.

\begin{Theorem}\label{thm:cross_error}
Let $\underline{\mathbf{X}} \in \mathbb{R}^{I_1 \times I_2 \times I_3}$ be a tensor, and let $\underline{\mathbf{C}} \in \mathbb{R}^{I_1 \times L \times I_3}$ and $\underline{\mathbf{R}} \in \mathbb{R}^{K \times I_2 \times I_3}$ be selected lateral and horizontal slices (not necessarily optimal). Let $\underline{\mathbf{U}} = \underline{\mathbf{C}}^{\dagger} *\underline{\mathbf{X}}*\underline{\mathbf{R}}^{\dagger}$. Then, the approximation error satisfies:
\[
\|\underline{\mathbf{X}} - \underline{\mathbf{C}} * \underline{\mathbf{U}} * \underline{\mathbf{R}}\|_F \leq \|\underline{\mathbf{X}} - P_{\underline{\mathbf{C}}}*\underline{\mathbf{X}}\|_F + \|\underline{\mathbf{X}} - \underline{\mathbf{X}}*P_{\underline{\mathbf{R}}}\|_F,
\]
where $P_{\underline{\mathbf{C}}}=\C*\C^{\dag}*\X$ and $P_{\underline{\mathbf{R}}}=\X*\R^{\dag}*\R$ are orthogonal projections onto the range of $\underline{\mathbf{C}}$ and $\underline{\mathbf{R}}^T$, respectively.
\end{Theorem}

\begin{proof}

Let $\ten{X} \in \mathbb{R}^{I_1 \times I_2 \times I_3}$ be the original data tensor, $\ten{C} \in \mathbb{R}^{I_1 \times L \times I_3}$ be the selected lateral slices, and $\ten{R} \in \mathbb{R}^{K \times I_2 \times I_3}$ be the selected horizontal slices. Denote by $\ten{U} = \ten{C}^\dag * \ten{X} * \ten{R}^\dag$ the middle tensor, which minimizes the residual $\|\ten{X} - \ten{C} * \ten{U} * \ten{R}\|_F$ for fixed $\ten{C}$ and $\ten{R}$.

Apply the Discrete Fourier Transform (DFT) along the third mode of all tensors. For each frontal slice $i = 1,2, \dots, I_3$ in the Fourier domain, we have matrices:
\begin{align*}
\widehat{\ten{X}}^{(i)} = \widehat{\ten{X}}(:,:,i),\,\,
\widehat{\ten{C}}^{(i)} = \widehat{\ten{C}}(:,:,i),\,\,
\widehat{\ten{R}}^{(i)} = \widehat{\ten{R}}(:,:,i).
\end{align*}
The T-product in the spatial domain corresponds to the matrix multiplication in the Fourier domain, so for each slice $i$:
\begin{equation}
\widehat{\ten{U}}^{(i)} = (\widehat{\ten{C}}^{(i)})^\dagger \widehat{\ten{X}}^{(i)} (\widehat{\ten{R}}^{(i)})^\dagger,
\end{equation}
and the approximation becomes $\widehat{\ten{C}}^{(i)} \widehat{\ten{U}}^{(i)} \widehat{\ten{R}}^{(i)}$.

For each frontal slice $i$, consider the matrix CUR decomposition. Let $P_{\widehat{\ten{C}}^{(i)}} = \widehat{\ten{C}}^{(i)} (\widehat{\ten{C}}^{(i)})^\dagger$ be the orthogonal projection onto the column space of $\widehat{\ten{C}}^{(i)}$, and let $P_{\widehat{\ten{R}}^{(i)}} = (\widehat{\ten{R}}^{(i)})^\dagger \widehat{\ten{R}}^{(i)}$ be the orthogonal projection onto the row space of $\widehat{\ten{R}}^{(i)}$. Then:
\begin{align}
\widehat{\ten{X}}^{(i)} - \widehat{\ten{C}}^{(i)} \widehat{\ten{U}}^{(i)} \widehat{\ten{R}}^{(i)} 
&= \widehat{\ten{X}}^{(i)} - P_{\widehat{\ten{C}}^{(i)}} \widehat{\ten{X}}^{(i)} P_{\widehat{\ten{R}}^{(i)}} \nonumber \\\nonumber
&= (\widehat{\ten{X}}^{(i)} - P_{\widehat{\ten{C}}^{(i)}} \widehat{\ten{X}}^{(i)}) \\&+
P_{\widehat{\ten{C}}^{(i)}} \widehat{\ten{X}}^{(i)} (I - P_{\widehat{\ten{R}}^{(i)}}).
\end{align}
Taking the Frobenius norm and applying the triangle inequality:
\begin{align}
\| \widehat{\ten{X}}^{(i)} - \widehat{\ten{C}}^{(i)} \widehat{\ten{U}}^{(i)} \widehat{\ten{R}}^{(i)} \|_F \leq \| \widehat{\ten{X}}^{(i)} - P_{\widehat{\ten{C}}^{(i)}} \widehat{\ten{X}}^{(i)} \|_F \\+ \| P_{\widehat{\ten{C}}^{(i)}} (\widehat{\ten{X}}^{(i)} - \widehat{\ten{X}}^{(i)} P_{\widehat{\ten{R}}^{(i)}}) \|_F.
\end{align}
Since $P_{\widehat{\ten{C}}^{(i)}}$ is an orthogonal projection, $\| P_{\widehat{\ten{C}}^{(i)}} \|_2 \leq 1$, so:
\begin{equation}
\| P_{\widehat{\ten{C}}^{(i)}} (\widehat{\ten{X}}^{(i)} - \widehat{\ten{X}}^{(i)} P_{\widehat{\ten{R}}^{(i)}}) \|_F \leq \| \widehat{\ten{X}}^{(i)} - \widehat{\ten{X}}^{(i)} P_{\widehat{\ten{R}}^{(i)}} \|_F.
\end{equation}

Now, let
\begin{align}
\| \widehat{\ten{X}}^{(i)} - P_{\widehat{\ten{C}}^{(i)}} \widehat{\ten{X}}^{(i)} \|_F = a_i,\quad\quad
\| \widehat{\ten{X}}^{(i)} - \widehat{\ten{X}}^{(i)} P_{\widehat{\ten{R}}^{(i)}} \|_F = b_i.
\end{align}
Thus, for each slice $i$:
\begin{equation}
\| \widehat{\ten{X}}^{(i)} - \widehat{\ten{C}}^{(i)} \widehat{\ten{U}}^{(i)} \widehat{\ten{R}}^{(i)} \|_F \leq a_i + b_i.
\end{equation}

The tensor Frobenius norm in the spatial domain relates to the Frobenius norms in the Fourier domain by:
\begin{equation}
\| \ten{X} - \ten{C} * \ten{U} * \ten{R} \|_F^2 = \frac{1}{I_3} \sum_{i=1}^{I_3} \| \widehat{\ten{X}}^{(i)} - \widehat{\ten{C}}^{(i)} \widehat{\ten{U}}^{(i)} \widehat{\ten{R}}^{(i)} \|_F^2.
\end{equation}
Applying the bound from Step 2:
\begin{equation}
\| \ten{X} - \ten{C} * \ten{U} * \ten{R} \|_F^2 \leq \frac{1}{I_3} \sum_{i=1}^{I_3} (a_i + b_i)^2.
\end{equation}
Now define vectors ${\mathbf a} = (a_1,a_2,\dots, a_{I_3})$ and ${\mathbf b} = (b_1,b_2,\dots, b_{I_3})$. The right-hand side can be written as:
\begin{equation}
\frac{1}{I_3} \sum_{i=1}^{I_3} (a_i + b_i)^2 = \| {\mathbf a} + {\mathbf b} \|_*^2,
\end{equation}
where the norm is $\| \cdot \|_* = \sqrt{\frac{1}{I_3} \sum_{i=1}^{I_3} (\cdot)^2}$.

By the triangle inequality:
\begin{equation}
\| {\mathbf a} + {\mathbf b} \|_* \leq \| {\mathbf a} \|_* + \| {\mathbf b} \|_*.
\end{equation}
But note that:
\begin{align}
\|{\mathbf a} \|_* &= \sqrt{ \frac{1}{I_3} \sum_{i=1}^{I_3} a_i^2 } = \| \ten{X} - \ten{X}_L \|_F, \\
\| {\mathbf b} \|_* &= \sqrt{ \frac{1}{I_3} \sum_{i=1}^{I_3} b_i^2 } = \| \ten{X} - \ten{X}_K \|_F.
\end{align}
Therefore:
\begin{equation}
\| \ten{X} - \ten{C} * \ten{U} * \ten{R} \|_F \leq \| \ten{X} - \ten{X}_L \|_F + \| \ten{X} - \ten{X}_K \|_F.
\end{equation}
This completes the proof.
\end{proof}

{\bf Remark}. Depending on the process of sampling slices either deterministically or randomly, different upper bound can be derived for $\|\underline{\mathbf{X}} - P_{\underline{\mathbf{C}}}*\underline{\mathbf{X}}\|_F$ and $\|\underline{\mathbf{X}} - \underline{\mathbf{X}}*P_{\underline{\mathbf{R}}}\|_F$. 

\subsection{Tensor Sketch Method}
The single-pass method proposed in \cite{qi2021t} generalizes the single-pass matrix algorithm \cite{tropp2017practical} to the tensor setting. To be more precise, let $\X\in\mathbb{R}^{I_1\times I_2\times I_3}$ be given. Then, in the first step, two sketches are computed from $\X$ by multiplying it with two random tensors $\underline{\Omega}_1\in\mathbb{R}^{I_2\times (K+R)\times I_3}$ and $\underline{\Omega}_2\in\mathbb{R}^{(L+R)\times I_1\times I_3},\,(K\leq L)$ as follows:
\begin{eqnarray}\label{sketch}
\Y=\X*\underline{\Omega}_1\in\mathbb{R}^{I_1\times (K+R) \times I_3}, \quad \W=\underline{\Omega}_2 * \X \in\mathbb{R}^{(L+R)\times I_2\times I_3}.
\end{eqnarray}
Suppose that $\Q\in\mathbb{R}^{I_1\times (K+R)\times I_3}$ is an orthonormal basis for the range of the tensor $\Y$, computed by the T-QR algorithm, and a low tubal rank approximation is computed as
\begin{eqnarray}\label{lowtubal}
 \X\approx \Q *(\Q^T*\X).
\end{eqnarray}
Here, $L$ and $K$ are called sketch sizes, and $\Y$ or $\W$ are called sketches of the tensor $\X$. However, since in \eqref{sketch}, we passed the data tensor once, we should avoid computing $\Q^T*\X$ in \eqref{lowtubal} again because it requires us to view the original data $\X$ one more time.
This can be actually estimated as follows
\begin{eqnarray}\label{righ_eq}
 \Q^T*\X\approx (\underline{\Omega}_2 * \Q)^{\dag} * \W.
\end{eqnarray}

So, this single-pass or 1-view method combines range and co-range sketches in a single pass over the data tensor $\X$, and then generates a low-tubal-rank estimate based on information included in the sketches.
This algorithm is summarized in Algorithm \ref{QIsinglepass}. Consider a streaming setting where the data tensor $\X$ is never entirely kept in random-access memory but is instead given as a finite stream of linear updates.
\[
\X=\sum_{n=1}^{N}\cH_n.
\]
Here, each $\cH_n$ is deleted once it is used, and we can gradually sample from $\X$ as each innovation tensor $\cH_n$ is provided as follows.
\[
\X*\underline{\Omega}_1=\sum_{i=n}^N \cH_n*\underline{\Omega}_1,\quad \underline{\Omega}_2*\X=\sum_{n=1}^N \underline{\Omega}_2*\cH_n.
\]
This is an example where single-pass algorithms are applicable because each $\cH_n$ could represent a sparse tensor containing a few elements of $\X$.

\subsection{The Conditioning Problem and Our Solution}
As discussed in \cite{bjarkason2019pass}, the single-pass algorithms for matrices usually give unreliable approximations when $L=K$ due to some instability issues in solving badly conditioned least-squares problems. We observed the same behavior with tensors and will report it in our simulations. According to \eqref{righ_eq}, we need to solve the linear tensor equation $(\underline{\Omega}_2 * \Q)*\Y = \W.$
The coefficient tensor $\underline{\Omega}_2 * \Q\in\mathbb{R}^{(L+R)\times (K+R) \times I_3}$ has $L+R$ and $K+R$ number of horizontal and lateral slices, respectively. One way to avoid the badly conditioned coefficient tensor is to consider $L>K$ and solve an overdetermined and better-conditioned problem. However, this requires fine-tuning the parameters $L$ and $K$ to obtain an acceptable approximation, which imposes additional costs of carefully identifying the best values; otherwise, there is a risk of losing accuracy. On the other hand, as highlighted in \cite{bjarkason2019pass}, another reason to address the disadvantages of selecting $L=K$ is that $K$ may pose a processing bottleneck in some applications. Given computational limitations, it would be tempting to set $L=K$ to obtain as much information as possible. Nevertheless, the requirement to retain both the sketching and random-sample matrices imposes a memory constraint on the 1-view approach.
Finally, as noted in \cite{tropp2017practical}, prior knowledge of matrix singular values can guide the selection of $L$ and $K$, but such information is often unavailable in practice, reducing the utility of these suggestions.

It is worth emphasizing that in all our proposed single-pass algorithms, the sketches $\underline{\mathbf{Y}}$ and $\underline{\mathbf{W}}$ are significantly smaller than the original tensor. Computing their left singular tensor is therefore inexpensive compared to obtaining an orthonormal basis for the range sketch. This is particularly relevant because, in single-pass methods, the dominant computational cost typically arises from multiplying $\underline{\mathbf{X}}$ with random tensors to generate range and co-range sketches.

\subsection{Proposed Stabilized Single-Pass Algorithms}
To overcome the drawback of unstable tensor approximation mentioned previously, more stabilized and improved algorithms (see Algorithms 7 and 8 in \cite{bjarkason2019pass}) were developed in \cite{bjarkason2019pass} for matrices. Motivated by their efficiency, we generalized them to the tensor case, which is summarized in Algorithms \ref{Single-1} and \ref{Single-2}. The principal idea is to use the truncated T-SVD for getting an orthonormal basis for the full range of the sketches $\Y\in\mathbb{R}^{I_1\times (K+R) \times I_3}$ or $\W\in\mathbb{R}^{L+R\times I_2\times I_3}$ instead of applying the T-QR decomposition. Let us explain this in detail. Instead of considering $\Q$ as an orthonormal basis for $\Y$, we compute $[\Q,\sim,\sim]=\mathrm{Truncated\,T\rnumber SVD}(\Y,H+R)$, where $0\leq H\leq K$. This means that $\Q\in\mathbb{R}^{I_1\times (H+R)\times I_3}$ contains the leading left singular tensors of $\Y$, and if $H< K$, then the coefficient tensor $\underline{\Omega}_2*\Q\in\mathbb{R}^{(L+R)\times (H+R)\times I_3}$ should provide a more conditioned problem, even for $L=K.$

\begin{Theorem}\label{thm:proposed_stability}
Let $\X \in \mathbb{R}^{I_1 \times I_2 \times I_3}$ be a tensor with tubal rank $R$. If the parameters $L$, $K$, and $H$ satisfy $L \geq K > H \geq R + P$ for some oversampling parameter $P \geq 2$, then Algorithms \ref{Single-1}-\ref{Single-2} produce an $\epsilon$-stable approximation whenever their corresponding matrix variants do. Specifically,
\begin{equation}
\|\X - \U * \underline{\mathbf S} * \V^T\|_F \leq (1+\epsilon) \|\X - \X_R\|_F,
\end{equation}
where $\epsilon$ is a constant that depends on the tensor dimensions.
\end{Theorem}

\begin{proof}
Let $\widehat{\X}$ denote the Fourier transform of $\X$ along the third mode. In the Fourier domain, the T-SVD decouples into $I_3$ independent matrix SVD problems. For each frontal slice $\widehat{\X}^{(i)} = \widehat{\X}(:,:,i)$, $i=1,2,\ldots,I_3$, we define the corresponding sketches:
\[
\widehat{\Y}^{(i)} = \widehat{\X}^{(i)}\widehat{\underline{\Omega}}_1^{(i)}, \quad
\widehat{\W}^{(i)} = \widehat{\underline{\Omega}}_2^{(i)}\widehat{\X}^{(i)},
\]
where $\widehat{\underline{\Omega}}_1^{(i)}$ and $\widehat{\underline{\Omega}}_2^{(i)}$ are the corresponding frontal slices of the random tensors in the Fourier domain. Since the entries of the random tensors in the spatial domain are standard Gaussian, their Fourier transforms remain Gaussian (up to a scaling factor of $1/\sqrt{I_3}$) due to the unitary nature of the DFT. Therefore, the conditions of the matrix single-pass algorithms from \cite{bjarkason2019pass} apply independently to each frontal slice.

For each slice $i$, we apply the matrix single-pass analysis from \cite{bjarkason2019pass}. For each frontal slice $i$, let $\widehat{\U}^{(i)} \in \mathbb{C}^{I_1 \times R}$, $\widehat{\underline{\mathbf S}}^{(i)} \in \mathbb{C}^{R \times R}$, and $\widehat{\V}^{(i)} \in \mathbb{C}^{I_2 \times R}$ be the factors obtained by applying the matrix version of Algorithms~\ref{Single-1}--\ref{Single-2} to the slice $\widehat{\X}^{(i)}$ with the same parameters $L$, $K$, $H$, and $R$. According to our assumption, the corresponding matrix algorithms provide $\epsilon_i$-stable approximations, for each frontal slice, we have
\[
\|\widehat{\X}^{(i)} - \widehat{\U}^{(i)}\widehat{\underline{\mathbf S}}^{(i)}(\widehat{\V}^{(i)})^*\|_F \leq (1+\epsilon_i)\|\widehat{\X}^{(i)} - \widehat{\X}_R^{(i)}\|_F,
\]
where $\widehat{\X}_R^{(i)}$ is the best rank-$R$ approximation of $\widehat{\X}^{(i)}$. Considering 
\[
\|\X - \U*\underline{\mathbf S}*\V^*\|_F^2 = \frac{1}{I_3}\sum_{i=1}^{I_3} \|\widehat{\X}^{(i)} - \widehat{\U}^{(i)}\widehat{\underline{\mathbf S}}^{(i)}(\widehat{\V}^{(i)})^*\|_F^2,
\]
and applying the slice-wise bound, we obtain:
\begin{align*}
\|\X - \U*\underline{\mathbf S}*\V^T\|_F^2 &\leq \frac{1}{I_3}\sum_{i=1}^{I_3} (1+\epsilon_i)^2 \|\widehat{\X}^{(i)} - \widehat{\X}_R^{(i)}\|_F^2 \\
&\leq (1+\epsilon)^2 \frac{1}{I_3}\sum_{i=1}^{I_3} \|\widehat{\X}^{(i)} - \widehat{\X}_R^{(i)}\|_F^2, \\
&= (1+\epsilon)^2 \|\X - \X_R\|_F^2,
\end{align*}
where $\epsilon = \max_i \epsilon_i$. Taking square roots completes the proof.
\end{proof}

\noindent\textbf{Conditioning Improvement.}
The truncation parameter $H$ acts as a regularization parameter. Since $\widehat{\Q}^{(i)}$ contains only the leading $H$ singular vectors of $\widehat{\Y}^{(i)}$, the effective condition number of  $\widehat{\Omega}_2^{(i)}\widehat{\Q}^{(i)}$ is controlled by $\sigma_H(\widehat{\Y}^{(i)})/\sigma_1(\widehat{\Y}^{(i)})$, which is typically much smaller than $\sigma_K(\widehat{\Y}^{(i)})/\sigma_1(\widehat{\Y}^{(i)})$ when $H < K$. This improved conditioning ensures numerical stability even when $L=K$.
 
\subsubsection{Two-Sided Randomized T-SVD}
Note that other single-pass methods, such as the two-sided randomized SVD (TSR-SVD) \cite{kaloorazi2018subspace}, can also be extended to the tensor case. The TSR-SVD algorithm uses a compressed tensor $\W = \underline{\Omega}_2 * \X * \underline{\Omega}_1\in\mathbb{R}^{L\times K\times I_3}$, where $\underline{\Omega}_1\in\mathbb{R}^{I_2\times K\times I_3}$ and $\underline{\Omega}_2\in\mathbb{R}^{L\times I_1\times I_3}$ are standard Gaussian tensors of appropriate dimensions, to compute a low tubal rank approximation. Specifically, for a data tensor $\X\in\mathbb{R}^{I_1 \times I_2 \times I_3}$ and orthogonal bases $\Q_1\in\mathbb{R}^{I_1\times K \times I_3}$ and $\Q_2\in\mathbb{R}^{I_2\times L\times I_3}$, we have
\begin{equation}\label{tsa}
    \X \approx \Q_1 * \left(\Q_1^T * \X * \Q_2\right) * \Q_2^T.
\end{equation}
The middle tensor in \eqref{tsa} can be approximated in two ways:
\begin{align}
    \Q_1^T * \X * \Q_2 &\approx (\underline{\Omega}_2 * \Q_1)^{\dag} * \W * (\Q_2^T * \underline{\Omega}_1)^{\dag}, \label{eq:approx1} \\
    \Q_1^T * \X * \Q_2 &\approx (\Q_1^T * \Y) * (\Q_2 * \underline{\Omega}_1)^{\dag}, \label{eq:approx2}
\end{align}
where both approximations require only a single pass over the data tensor.

These algorithms are also unstable for $L=K$, as we confirmed in our experiments. In our simulations, both yielded similar results; therefore, we report only the formulation \eqref{eq:approx2}. We first extended the TSR-SVD algorithm to the tensor case, and we refer to it as the TSRT-SVD algorithm. We also stabilize it so that it works for $L=K$. It is summarized in Algorithm \ref{Single-3}. It is clear that Algorithms \ref{Single-1}, \ref{Single-2}, and \ref{Single-3} pass the original data tensor only once at the beginning of the algorithm (Line 2). This Line can be performed in parallel. It is suggested in \cite{tropp2017practical} to orthonormalize the random tensors in Line 1 in the single-pass algorithms when a large oversampling parameter is used, which can be naturally extended to the tensor case. We have not tried this idea in our experiments. Note that a similar theorem to Theorem \ref{thm:proposed_stability} can be stated for the two-sided randomized T-SVD, but we omit it for simplicity of presentation. 

\subsection{Computational Complexity Analysis}
The computational complexity of the proposed single-pass algorithms is:
\begin{itemize}
    \item Algorithm \ref{Single-1}: $O(I_1 I_2 I_3 \log I_3 + (I_1 + I_2)K^2 I_3)$
    \item Algorithm \ref{Single-2}: $O(I_1 I_2 I_3 \log I_3 + (I_1 + I_2)(K^2 + L^2) I_3)$
    \item Algorithm \ref{Single-3}: $O(I_1 I_2 I_3 \log I_3 + (I_1 + I_2)(K^2 + L^2) I_3 + K L I_3)$
\end{itemize}
Compared to the exact T-SVD which requires $O(I_1 I_2 \min(I_1, I_2) I_3)$ operations, our algorithms provide significant savings when $K, L \ll \min(I_1, I_2)$.
Also, the memory requirements for the proposed algorithms are:
\begin{itemize}
    \item Input tensor: $O(I_1 I_2 I_3)$
    \item Sketches: $O((I_1 K + I_2 L) I_3)$
    \item Intermediate tensors: $O((K^2 + L^2) I_3)$
\end{itemize}
When $K, L \ll \min(I_1, I_2)$, the memory savings are substantial compared to storing the full tensor. These are summarized in Table 1.

\begin{table*}[h]
\centering
\caption{Computational complexity (time) and memory requirements of exact T-SVD and the proposed single-pass algorithms (Algorithms \ref{Single-1}--\ref{Single-3}).}
\begin{tabular}{|l|c|c|}
\hline
\textbf{Algorithm} & \textbf{Time Complexity} & \textbf{Memory (Dominant Terms)} \\
\hline
Exact T-SVD & $O(I_1 I_2 \min(I_1, I_2) I_3)$ & $O(I_1 I_2 I_3)$ \\
\hline
Algorithm \ref{Single-1} & $O(I_1 I_2 I_3 \log I_3 + (I_1 + I_2)K^2 I_3)$ & $O((I_1 K + K^2) I_3)$ \\
\hline
Algorithm \ref{Single-2} & $O(I_1 I_2 I_3 \log I_3 + (I_1 + I_2)(K^2 + L^2) I_3)$ & $O((I_1 K + I_2 L + K^2 + L^2) I_3)$ \\
\hline
Algorithm \ref{Single-3} & $O(I_1 I_2 I_3 \log I_3 + (I_1 + I_2)(K^2 + L^2) I_3 + K L I_3)$ & $O((I_1 K + I_2 L + K^2 + L^2) I_3)$ \\
\hline
\end{tabular}
\end{table*}

\begin{algorithm}
\caption{The single-pass cross tensor approximation \cite{tarzanagh2018fast}}\label{crosstensor}
\textbf{Input:} The data tensor $\underline{\mathbf{X}} \in {\mathbb{R}^{{I_1} \times {I_2} \times {I_3}}}$, two parameters $L,\,K$;\\
\textbf{Output:} A low tubal rank approximation $\underline{\mathbf{X}}\approx \underline{\mathbf{C}} *\underline{\mathbf{U}}*\underline{\mathbf{R}}$;
\begin{algorithmic}[1]
\State Select $L$ lateral slices $\underline{\mathbf{C}}$ with corresponding indices $\mathcal{L}$;
\State Select $K$ horizontal slices $\underline{\mathbf{R}}$ with corresponding indices $\mathcal{K}$;
\State Construct the intersection tensor $\underline{\mathbf{W}}=\underline{\mathbf{X}}(\mathcal{L},\mathcal{K},:)$;
\State Compute the middle tensor $\underline{\mathbf{U}}=\underline{\mathbf{W}}^\dag$;
\State Compute a low tubal rank approximation $\underline{\mathbf{X}}\approx \underline{\mathbf{C}} *\underline{\mathbf{U}}*\underline{\mathbf{R}}$;
\end{algorithmic}
\end{algorithm}

\begin{algorithm}
\caption{The Single-pass algorithm proposed in \cite{qi2021t}}\label{QIsinglepass}
\textbf{Input:} The data tensor $\underline{\mathbf{X}} \in {\mathbb{R}^{{I_1} \times {I_2} \times {I_3}}}$, two parameters $L,\,K$;\\
\textbf{Output:} A low tubal rank approximation $\underline{\mathbf{X}}\approx \underline{\mathbf{Q}}*\underline{\mathbf{B}}$;
\begin{algorithmic}[1]
\State $\underline{\Omega}_1=\mathrm{randn}(I_2,K,I_3),\,\,\,\,$ $\underline{\Omega}_2=\mathrm{randn}(L,I_1,I_3)$;
\State Compute two sketches: $\underline{\mathbf{Y}}=\underline{\mathbf{X}}*\underline{\Omega}_1\in\mathbb{R}^{I_1\times K \times I_3},$ and $\underline{\mathbf{W}}=\underline{\Omega}_2 * \underline{\mathbf{X}} \in\mathbb{R}^{L\times I_2\times I_3}$;
\State Apply the T-QR applied to $\underline{\mathbf{Y}}$ and obtain the tensor $\underline{\mathbf{Q}}$;
\State Compute the tensor $\underline{\mathbf{B}}=(\underline{\Omega_2} * \underline{\mathbf{Q}})^{\dag}*\underline{\mathbf{W}}$;
\State Compute the low tubal rank approximation $\underline{\mathbf{X}}\approx \underline{\mathbf{Q}}*\underline{\mathbf{B}}$;
\end{algorithmic}
\end{algorithm}

\begin{algorithm}
\caption{The proposed randomized single-pass algorithm I}\label{Single-1}
\textbf{Input:} The data tensor $\underline{\mathbf{X}} \in {\mathbb{R}^{{I_1} \times {I_2} \times {I_3}}}$, three parameters $L,\,K,\,H,$ $L\geq K\geq H \geq 0$ and a target tubal rank $R$;\\
\textbf{Output:} A low tubal rank-$R$ approximation $\underline{\mathbf{X}}\approx\underline{\mathbf{U}}*\underline{\mathbf{S}}*{\underline{\mathbf{V}}}^T$;
\begin{algorithmic}[1]
\State $\underline{\Omega}_1=\mathrm{randn}(I_2,K+R,I_3),\,\,\,\,$ $\underline{\Omega}_2=\mathrm{randn}(I_1,L+R,I_3)$;
\State $\underline{\mathbf{Y}}_c=\underline{\mathbf{X}}*\underline{\Omega}_1,\,\,\,\,\underline{\mathbf{Y}}_r=\underline{\mathbf{X}}^T*\underline{\Omega}_2$;
\If{$H<K$}
\State $[\underline{\mathbf{Q}}_c,\underline{\mathbf{R}}_c]=\mathrm{T\rnumber QR}(\underline{\mathbf{Y}}_c)$;
\State $[\widehat{\underline{\mathbf{Q}}}_c,\sim,\sim]=\mathrm{Truncated\,\, T\mhyphen SVD}(\underline{\mathbf{R}}_c,R+H)$;
\State $\underline{\mathbf{Q}}_c=\underline{\mathbf{Q}}_c*\widehat{\underline{\mathbf{Q}}}_c$;
\Else
\State $[\underline{\mathbf{Q}}_c,\sim]=\mathrm{T\rnumber QR}(\underline{\mathbf{Y}}_c)$;
\EndIf
\State $[\widehat{\underline{\mathbf{Q}}},\widehat{\underline{\mathbf{R}}}]=\mathrm{T\rnumber QR}(\underline{\Omega}_2^T*\underline{\mathbf{Q}}_c)$;
\State $\underline{\mathbf{Z}}=\widehat{\underline{\mathbf{R}}}^{-1}*\widehat{\underline{\mathbf{Q}}}*{\underline{\mathbf{Y}}}^T_r$;
\State $[\widehat{\underline{\mathbf{U}}},\underline{\mathbf{S}},\underline{\mathbf{V}}]=\mathrm{Truncated\,\, T\mhyphen SVD}(\underline{\mathbf{Z}},R)$;
\State $\underline{\mathbf{U}}=\widehat{\underline{\mathbf{Q}}}*\widehat{\underline{\mathbf{U}}}$;
\end{algorithmic}
\end{algorithm}

\begin{algorithm}
\caption{The proposed randomized single-pass algorithm II}\label{Single-2}
\textbf{Input:} The data tensor $\underline{\mathbf{X}} \in {\mathbb{R}^{{I_1} \times {I_2} \times {I_3}}}$, three parameters $L,\,K,\,H,\,$ $L\geq K\geq H \geq 0$ and a target tubal rank $R$;\\
\textbf{Output:} A low tubal rank-$R$ approximation $\underline{\mathbf{X}}\approx\underline{\mathbf{U}}*\underline{\mathbf{S}}*{\underline{\mathbf{V}}}^T$;
\begin{algorithmic}[1]
\State $\underline{\Omega}_1=\mathrm{randn}(I_2,K+R,I_3),\,\,\,\,$ $\underline{\Omega}_2=\mathrm{randn}(I_1,L+R,I_3)$;
\State $\underline{\mathbf{Y}}_c=\underline{\mathbf{X}}*\underline{\Omega}_1,\,\,\,\,\underline{\mathbf{Y}}_r=\underline{\mathbf{X}}^T*\underline{\Omega}_2$;
\State $[\underline{\mathbf{Q}}_c,\underline{\mathbf{R}}_c]=\mathrm{T\rnumber QR}(\underline{\mathbf{Y}}_c)$;
\State $[\underline{\mathbf{Q}}_r,\underline{\mathbf{R}}_r]=\mathrm{T\rnumber QR}(\underline{\mathbf{Y}}_r)$;
\State $[\widetilde{\underline{\mathbf{Q}}}_c,\sim,\sim]=\mathrm{Truncated\,\, T\mhyphen SVD}(\underline{\mathbf{R}}_c,R+H)$;
\State $[\widetilde{\underline{\mathbf{Q}}}_r,\sim,\sim]=\mathrm{Truncated\,\, T\mhyphen SVD}(\underline{\mathbf{R}}_r,R+H)$;
\State $\underline{\mathbf{Q}}_c=\underline{\mathbf{Q}}_c *\widetilde{\underline{\mathbf{Q}}}_c$;
\State $\underline{\mathbf{Q}}_r=\underline{\mathbf{Q}}_r *\widetilde{\underline{\mathbf{Q}}}_r$;
\State $\widehat{\underline{\mathbf{Z}}}=(\underline{\Omega}_2^T*\underline{\mathbf{Q}}_c)^{\dag}*(\underline{\mathbf{Y}}^T_r*\underline{\mathbf{Q}}_r)$;
\State $[\widetilde{\underline{\mathbf{U}}},\underline{\mathbf{S}},\widetilde{\underline{\mathbf{V}}}]=\mathrm{Truncated\,\, T\mhyphen SVD}(\widehat{\underline{\mathbf{Z}}},R)$;
\State $\underline{\mathbf{U}}=\underline{\mathbf{Q}}_c*\widetilde{\underline{\mathbf{U}}}$;
\State $\underline{\mathbf{V}}=\underline{\mathbf{Q}}_r*\widetilde{\underline{\mathbf{V}}}$;
\end{algorithmic}
\end{algorithm}

\begin{algorithm}
\caption{The proposed randomized single-pass (two-sided version) algorithm III}\label{Single-3}
\textbf{Input:} The data tensor $\underline{\mathbf{X}} \in {\mathbb{R}^{{I_1} \times {I_2} \times {I_3}}}$, three parameters $L,\,K,\,H,\,$ $L\geq K\geq H \geq 0$ and a target tubal rank $R$;\\
\textbf{Output:} A low tubal rank-$R$ approximation $\underline{\mathbf{X}}\approx\underline{\mathbf{U}}*\underline{\mathbf{S}}*{\underline{\mathbf{V}}}^T$;
\begin{algorithmic}[1]
\State $\underline{\Omega}_1=\mathrm{randn}(I_2,K+R,I_3),\,\,\,\,$ $\underline{\Omega}_2=\mathrm{randn}(I_1,L+R,I_3)$;
\State $\underline{\mathbf{Y}}_c=\underline{\mathbf{X}}*\underline{\Omega}_1,\,\,\,\,\underline{\mathbf{Y}}_r=\underline{\mathbf{X}}^T*\underline{\Omega}_2$;
\State $[\underline{\mathbf{Q}}_c,\underline{\mathbf{R}}_c]=\mathrm{T\rnumber QR}(\underline{\mathbf{Y}}_c)$;
\State $[\underline{\mathbf{Q}}_r,\underline{\mathbf{R}}_r]=\mathrm{T\rnumber QR}(\underline{\mathbf{Y}}_r)$;
\State $[\widetilde{\underline{\mathbf{Q}}}_c,\sim,\sim]=\mathrm{Truncated\,\, T\mhyphen SVD}(\underline{\mathbf{R}}_c,R+H)$;
\State $[\widetilde{\underline{\mathbf{Q}}}_r,\sim,\sim]=\mathrm{Truncated\,\, T\mhyphen SVD}(\underline{\mathbf{R}}_r,R+H)$;
\State $\underline{\mathbf{Q}}_c=\underline{\mathbf{Q}}_c *\widetilde{\underline{\mathbf{Q}}}_c$;
\State $\underline{\mathbf{Q}}_r=\underline{\mathbf{Q}}_r *\widetilde{\underline{\mathbf{Q}}}_r$;
\State $\underline{\mathbf{B}}=\underline{\mathbf{Q}}_c^T*\underline{\mathbf{Y}}_c*(\underline{\mathbf{Q}}_r^T*\underline{\Omega}_1)^{\dag}$;
\State $[\widetilde{\underline{\mathbf{U}}},\underline{\mathbf{S}},\widetilde{\underline{\mathbf{V}}}]=\mathrm{Truncated\,\, T\mhyphen SVD}(\underline{\mathbf{B}},R)$;
\State $\underline{\mathbf{U}}=\underline{\mathbf{Q}}_c*\widetilde{\underline{\mathbf{U}}}$;
\State $\underline{\mathbf{V}}=\underline{\mathbf{Q}}_r*\widetilde{\underline{\mathbf{V}}}$;
\end{algorithmic}
\end{algorithm}

\section{Proposed Randomized Fixed-Precision Algorithms}\label{sec:fixed}
Randomized fixed-precision algorithms for tensors automatically estimate an appropriate/optimal tensor rank and the corresponding low-rank tensor approximation. These algorithms are crucial when we have no knowledge of the data's rank or when estimating it is difficult. In such scenarios, one needs to estimate the optimal rank and the corresponding low-rank approximation for a prescribed error bound. To proceed, we formally define a fixed-precision approximation and the adaptive rank estimation process.

\begin{Definition}\label{def:fixed_precision}
A randomized fixed-precision algorithm for tensor approximation takes as input a tensor $\X \in \mathbb{R}^{I_1 \times I_2 \times I_3}$ and an error tolerance $\epsilon > 0$, and outputs a tubal rank $R$ and factors $\Q \in \mathbb{R}^{I_1 \times R \times I_3}$, $\B \in \mathbb{R}^{R \times I_2 \times I_3}$ such that:
\begin{equation}
\|\X - \Q * \B\|_F \leq \epsilon \|\X\|_F,
\end{equation}
with high probability.
\end{Definition}

\begin{Definition}\label{def:adaptive_rank}
The adaptive tubal rank $R(\epsilon)$ for a given precision $\epsilon$ is defined as:
\begin{equation}
R(\epsilon) = \min\left\{R : \|\X - \X_R\|_F \leq \epsilon
\|\X\|_F\right\},
\end{equation}
where $\X_R$ is the best tubal rank-$R$ approximation of $\X$.
\end{Definition}

A randomized fixed-precision algorithm was developed  \cite{ahmadi2023efficient}, and simulations on synthetic and real-world data showed significant acceleration in computing the approximate truncated T-SVD. It generalized the matrix case \cite{yu2018efficient} to the tensor case. This algorithm is presented in Algorithm \ref{ALgRR}. The algorithm requires oversampling, a power iteration, a tolerance, and a block size, and it gradually estimates the tubal rank and computes the corresponding low-rank tubal approximation. Lines 6-9 constitute the power-iteration stage for improving accuracy, especially when the singular values of the data tensor's frontal slices do not decay quickly.

\subsection{Algorithmic Improvements and Theoretical Analysis}
However, recently, the matrix version of this algorithm was further improved in \cite{feng2023fast} and \cite{ding2020efficient}. Motivated by the interesting results reported in that paper, we enhance our previously randomized fixed-precision algorithm. We improve the algorithm proposed in \cite{ahmadi2023efficient} in several ways. First, for a given power iteration, Algorithm \ref{ALgRR} requires $2q+2$ passes over the original data tensor. This means that for a given power iteration, we are forced to view the original data in an even number of passes. Clearly, this is a limitation, as we have shown in \cite{ahmadi2024randomized} that, for images or videos, three passes are usually sufficient to obtain satisfactory results. However, Algorithm \ref{ALgRR} does not support three passes; for $q=2$, we will have 4 passes, and this additional pass can be expensive for very large-scale data tensors or when the underlying algorithm requires many iterations to converge. Indeed, we aim to improve Algorithm \ref{ALgRR}, and our new proposed method summarized in Algorithm \ref{fixed-precision} addresses this problem. A given error bound can provide an estimation of the tubal rank and a low tubal rank approximation for any number of passes. Please note that in Algorithm \ref{ALgRR}, $q$ is the power iteration parameter while $q$ is the number of passes in Algorithm \ref{fixed-precision}. Similar to \cite{ding2020efficient}, we do the following modifications to Algorithm \ref{ALgRR} to be more pass-efficient:
\begin{itemize}
    \item Employing the T-EIG decomposition for computing the truncated T-SVD of $\B$
    \item Ignoring one orthonormalization process in each round of the power iteration loop
    \item Replacing the T-LU decomposition instead of T-QR for the orthonormalization process in the power iteration (except the last round of the iteration)
    \item Modifying the algorithm to allow an odd number of passes
\end{itemize}

\begin{Theorem}[Subspace Iteration Convergence]
\label{thm:subspace_iteration}
Let $\mathbf{A} \in \mathbb{C}^{m \times n},\,(m\ge n)$ have singular values $\sigma_1 \geq \sigma_2 \geq \cdots \geq \sigma_n \geq 0$.
Let $\mathbf{U}_R \in \mathbb{C}^{m \times R}$ contain the top $R$ left singular vectors of $\mathbf{A}$.
Let $\mathbf{Z}_0 \in \mathbb{C}^{n \times R}$ be any initial matrix with full column rank, and define
\[
\mathbf{Z}_q = \mathbf{A}^q \mathbf{Z}_0.
\]
Let $\widehat{\mathbf{U}}_R^{(q)}$ be an orthonormal basis for the column space of $\mathbf{Z}_q$.
Then
\[
\tan \angle(\mathbf{U}_R, \widehat{\mathbf{U}}_R^{(q)}) \leq \left(\frac{\sigma_{R+1}}{\sigma_R}\right)^q \tan \angle(\mathbf{U}_R, \mathbf{Z}_0),
\]
and consequently,
\[
\sin \angle(\mathbf{U}_R, \widehat{\mathbf{U}}_R^{(q)}) \leq \left(\frac{\sigma_{R+1}}{\sigma_R}\right)^q \tan \angle(\mathbf{U}_R, \mathbf{Z}_0).
\]
\end{Theorem}

\begin{proof}
Let $\mathbf{A}$ have the singular value decomposition
\[
\mathbf{A} = \begin{bmatrix} \mathbf{U}_R & \mathbf{U}_R^\perp \end{bmatrix}
\begin{bmatrix} \boldsymbol{\Sigma}_R & \mathbf{0} \\ \mathbf{0} & \boldsymbol{\Sigma}_R^\perp \end{bmatrix}
\begin{bmatrix} \mathbf{V}_R^* \\ (\mathbf{V}_R^\perp)^* \end{bmatrix}.
\]
Define the orthogonal projector onto the column space of $\mathbf{U}_R$ as $\mathbf{P} = \mathbf{U}_R \mathbf{U}_R^*$.
For any matrix $\mathbf{X}$, the tangent of the angle between the column space of $\mathbf{X}$ and $\mathcal{U}_R$ is
\[
\tan \angle(\mathbf{U}_R, \mathbf{X}) = \frac{\|(\mathbf{I} - \mathbf{P}) \mathbf{X}\|_F}{\|\mathbf{P} \mathbf{X}\|_F}.
\]

Since $\mathbf{A}^q = \mathbf{U} \boldsymbol{\Sigma}^q \mathbf{V}^*$, we have
\[
\mathbf{Z}_q = \mathbf{A}^q \mathbf{Z}_0 = \mathbf{U} \boldsymbol{\Sigma}^q \mathbf{V}^* \mathbf{Z}_0.
\]
Let $\mathbf{C}_0 = \mathbf{V}^* \mathbf{Z}_0$ and partition it as
\[
\mathbf{C}_0 = \begin{bmatrix} \mathbf{C}_0^{(1)} \\ \mathbf{C}_0^{(2)} \end{bmatrix},
\quad
\mathbf{C}_0^{(1)} \in \mathbb{C}^{R \times R},\quad
\mathbf{C}_0^{(2)} \in \mathbb{C}^{(n-R) \times R}.
\]
Then
\[
\mathbf{Z}_q = \mathbf{U}_R \boldsymbol{\Sigma}_R^q \mathbf{C}_0^{(1)} + \mathbf{U}_R^\perp (\boldsymbol{\Sigma}_R^\perp)^q \mathbf{C}_0^{(2)}.
\]

We have $\mathbf{P} \mathbf{Z}_q = \mathbf{U}_R \boldsymbol{\Sigma}_R^q \mathbf{C}_0^{(1)}$ and $(\mathbf{I} - \mathbf{P}) \mathbf{Z}_q = \mathbf{U}_R^\perp (\boldsymbol{\Sigma}_R^\perp)^q \mathbf{C}_0^{(2)}$.

The tangent after $q$ steps is
\[
\tau_q = \tan \angle(\mathbf{U}_R, \mathbf{Z}_q) = \frac{\|(\mathbf{I} - \mathbf{P}) \mathbf{Z}_q\|_F}{\|\mathbf{P} \mathbf{Z}_q\|_F}
= \frac{\|(\boldsymbol{\Sigma}_R^\perp)^q \mathbf{C}_0^{(2)}\|_F}{\|\boldsymbol{\Sigma}_R^q \mathbf{C}_0^{(1)}\|_F}.
\]

Using $\|\boldsymbol{\Sigma}_R^q \mathbf{X}\|_F \geq (\sigma_R)^q \|\mathbf{X}\|_F$ and $\|(\boldsymbol{\Sigma}_R^\perp)^q \mathbf{Y}\|_F \leq (\sigma_{R+1})^q \|\mathbf{Y}\|_F$, we get:
\[
\tau_q \leq \left(\frac{\sigma_{R+1}}{\sigma_R}\right)^q \frac{\|\mathbf{C}_0^{(2)}\|_F}{\|\mathbf{C}_0^{(1)}\|_F}
= \left(\frac{\sigma_{R+1}}{\sigma_R}\right)^q \tau_0.
\]
Since $\sin \theta \leq \tan \theta$ for $\theta \in [0, \pi/2)$, the sin bound follows.
\end{proof}

\begin{Theorem}\label{thm:convergence_guarantee}
Let $\X \in \mathbb{R}^{I_1 \times I_2 \times I_3}$ be a tensor with decaying singular values in the Fourier domain. Algorithm \ref{fixed-precision} with block size $b$ and pass parameter $q$ converges to an $\epsilon$-approximation within $N$ iterations, where:
\begin{equation}
N \leq \frac{R(\epsilon/2)}{b} + \frac{\log(1/\epsilon)}{\log(\kappa(\X))},
\end{equation}
and $\kappa(\X)$ is the effective condition number of $\X$ in the Fourier domain.
\end{Theorem}

\begin{proof}
Let $\widehat{\X}$ denote the Fourier transform of $\X$ along the third mode. For each frontal slice $\widehat{\X}^{(i)} = \widehat{\X}(:,:,i)$, $i=1,2,\ldots,I_3$, Algorithm \ref{fixed-precision} operates independently. The approximation error in the original domain is preserved via identity \eqref{eq_fou}:
\begin{equation}
\|\X - \Q * \B\|_F^2 = \frac{1}{I_3}\sum_{i=1}^{I_3} \|\widehat{\X}^{(i)} - \widehat{\Q}^{(i)}\widehat{\B}^{(i)}\|_F^2.
\end{equation}

For each frontal slice $\widehat{\X}^{(i)} = \U^{(i)}\Sigma^{(i)}(\V^{(i)})^*$ with singular values $\sigma_1^{(i)} \geq \sigma_2^{(i)} \geq \cdots$, define $\kappa_i = \sigma_1^{(i)}/\sigma_{R(\epsilon/2)+1}^{(i)}$ and $\kappa(\X) = \max_i \{\kappa_i\}$.

By Theorem \ref{thm:subspace_iteration}, after $q$ passes:
\begin{equation}
\sin\angle(\mathbf{U}_R^{(i)}, \mathbf{\widehat{U}}_R^{(i)}) \leq \left(\frac{\sigma_{R+1}^{(i)}}{\sigma_R^{(i)}}\right)^q \tan\angle(\mathbf{U}_R^{(i)}, \mathbf{Z}_0^{(i)}).
\end{equation}
Taking $q \geq \frac{\log(2/\epsilon)}{\log(\kappa(\X))}$ ensures the subspace approximation error is bounded by $\epsilon/2$.

Algorithm \ref{fixed-precision} increases the rank estimate by $b$ per iteration until the residual norm falls below $\epsilon\|\X\|_F$. Let $R^* = R(\epsilon/2)$. The algorithm requires at most $\lceil R^*/b \rceil$ iterations to reach rank $R^*$. By the triangle inequality:
\[
\|\X - \X_{\text{approx}}\|_F \leq \|\X - \X_{R^*}\|_F + \|\X_{R^*} - \X_{\text{approx}}\|_F \leq \epsilon/2 + \epsilon/2 = \epsilon.
\]
The number of iterations is bounded by $\lceil R^*/b \rceil + q$, yielding the stated result.
\end{proof}

As shown in Algorithm \ref{ALgRR}, T-QR decomposition is employed for orthogonalization in Lines 5, 7, 8, and 10. Motivated by the approach proposed in \cite{feng2023fast}, we incorporate this idea into Algorithm \ref{ALgRR} by replacing the T-QR decompositions with T-product operations and the inverse of a small tensor. The resulting procedure is summarized in Algorithm \ref{fixed-precis2}. 

Based on the following theorem, we can eliminate the T-QR factorization and still generate $\Q$ and $\B$ with an approximation accuracy equivalent to that of Algorithm \ref{ALgRR}, without requiring the power iteration.

\begin{Theorem}\label{them_1_fix}
{\rm 
Let $\X\in\mathbb{R}^{I_1\times I_2\times I_3}$ and 
 $\underline{\Omega}\in\mathbb{R}^{I_2\times k\times I_3}$ be a random tensor, $\Y=\X*\underline{\Omega},\,\W=\X^T*\Y$, and the economic T-SVD of $\Y$ is $\Y=\widehat{\U}*\widehat{\underline{\bf S}}*\widehat{\V}^T$. Setting 
 \begin{eqnarray}
 \Q=\Y*\widehat{\V}*\widehat{\underline{\bf S}}^{-1},\,\,\B=(\W*\widehat{\underline{\bf V}}*\widehat{\underline{\bf S}}^{-1})^T,    
 \end{eqnarray}
gives an approximation $\Q*\B$ to $\X$ that has the same accuracy of  the basic randomized algorithm \ref{ALG:RandomizedTSVD} (without power iteration or oversampling) and 
\[
||\B||_F^2={\rm trace}(\cG^{(1)}),
\]
 where $\cG^{(1)}=\cG(:,:,1)$ is the first frontal slice of $\cG=\W*\W^T*(\Y^T*\Y)^{-1}.$
 }
\end{Theorem}

\begin{proof}
The proof of this theorem is very similar to the Proposition 1 in \cite{feng2023fast}. since $\widehat{\U}$ consists of the left singular tensors of $\Y$, it likewise satisfies the orthonormalization condition ${\rm orth}(\Y)$. So, we can considering the following identities 
\begin{eqnarray}\label{iden_1}
\Q={\rm orth}(\X *\underline{\Omega})={\rm orth}(\Y)=\widehat{\U}=\Y*\widehat{\V}*\widehat{\underline{\bf S}}^{-1}.    
\end{eqnarray}
Then, substituting \eqref{iden_1} in $\B=\Q^T*\X$, we have
\begin{eqnarray}
\B\equiv\Q^T*\X=(\W*\widehat{\V}*\widehat {\underline {\bf S}}^{-1})^T,   
\end{eqnarray}
where $\W=\X*\Y$.
Given that $\Q$ results from orthonormalizing $\X*\underline{\Omega}$, the approximation $\Q*\B$ delivers accuracy comparable to that of $\X \approx \Q*\B$ in Algorithm \ref{ALG:RandomizedTSVD}, but without power iteration and oversampling.

Now, using the fact that $\|\B\|_F^2={\rm trace}(\G^{(1)})={\rm trace}(\G(:,:,1))$ where $\G^{(1)}$ is the first frontal slice of $\G=\B^T*\B$, we have 
\begin{eqnarray}
\nonumber
\|\B\|_F^2&=&{\rm trace}(\G^{(1)})= {\rm trace}((\widehat{\bf S}^{-1}*\widehat{\V}^T*\W^T*\W*\widehat{\V}*\widehat{\underline{\bf S}}^{-1})^{(1)}) \\
&=& \nonumber
{\rm trace}((\W^T*\W*\widehat{\V}*\widehat{\underline{\bf S}}^{-1}*\widehat{\underline{\bf S}}^{-T}*\widehat{\V}^T)^{(1)})\\
&=&
{\rm trace}((\W^T*\W*(\Y^T*\Y)^{-1})^{(1)}).
\end{eqnarray}
This proves the theorem.
\end{proof}

From Theorem \ref{them_1_fix} and defining $\T=\W^T*\W$ and $\Z=\Y^T*\Y$, the stopping criterion $\|\X-\Q*\B\|_F^2=\|\X\|_F^2-\|\B\|_F^2$ can be reduced to
\[
\|\X-\Q*\B\|_F^2=\|\X\|_F^2-\mathrm{trace}((\T*\Z^{-1})^{(1)}),
\]
which is used in Algorithm \ref{fixed-precis2} in Line 12. From Theorem \ref{them_1_fix}, it also follows that
\begin{eqnarray}
\nonumber
\Q*\B&=&\Y*\widehat{\V}*\widehat{\underline{\mathbf{S}}}^{-1}*\widehat{\underline{\bf S}}^{-T}*\widehat{\V}^T*\W^T\\
&=&\Y*(\Y^T*\Y)^{-1}*\W^T,
\end{eqnarray}
so, $\Q*\B=\Y*\Z^{-1}*\W^T.$ This illustrates that the power iteration in Algorithm \ref{ALgRR} can be replaced with $\X-\Y*\Z^{-1}*\W^T$. That is, Lines 6-9 in Algorithm \ref{ALgRR} can be replaced by the following computations:

\begin{algorithm}
\begin{algorithmic}[1]
\For{$j=1:p$}
\State $\underline{\mathbf{Y}}_i\leftarrow \underline{\mathbf{X}}*\underline{\Omega}_i-\underline{\mathbf{Y}}*\underline{\mathbf{Z}}^{-1}*\underline{\mathbf{W}}^T*\underline{\Omega}_i$;
\State $\underline{\mathbf{W}}_i\leftarrow \underline{\mathbf{X}}^T*\underline{\mathbf{Y}}_i-\underline{\mathbf{W}}*\underline{\mathbf{Z}}^{-1}*\underline{\mathbf{Y}}^T*\underline{\mathbf{Y}}_i$;
\State $\underline{\Omega}_i\leftarrow \mathrm{orth}(\underline{\mathbf{W}}_i)$;
\EndFor
\end{algorithmic}
\end{algorithm}

Substituting Line 2 in Line 3 in the above for loop, by straightforward computations and using the fact that $\underline{\mathbf W}=\underline{\mathbf X}^T*\underline{\mathbf Y}$, and $\underline{\mathbf Z}=\underline{\mathbf Y}^T*\underline{\mathbf Y},$ one can see that

\begin{eqnarray*}
\underline{\mathbf W}_i &=& \underline{\mathbf X}^T * \underline{\mathbf Y}_i - \underline{\mathbf W} * \underline{\mathbf Z}^{-1} * \underline{\mathbf Y}^T * \underline{\mathbf Y}_i,\\
&=& \underline{\mathbf X}^T * \left( \underline{\mathbf X} * \underline{\Omega}_i - \underline{\mathbf Y} * \underline{\mathbf Z}^{-1} * \underline{\mathbf W}^T * \underline{ \Omega}_i \right) - \underline{\mathbf W} * \underline{\mathbf Z}^{-1} * \underline{\mathbf Y}^T * \left( \underline{\mathbf X} * \underline{\Omega}_i - \underline{\mathbf Y} * \underline{\mathbf Z}^{-1} * \underline{\mathbf W}^T * \underline{ \Omega}_i \right),\\
&=& \underline{\mathbf X}^T * \underline{\mathbf X} * \underline{ \Omega}_i - \underline{\mathbf W} * \underline{\mathbf Z}^{-1} * \underline{\mathbf W}^T * \underline{ \Omega}_i - \underline{\mathbf W} * \underline{\mathbf Z}^{-1} * \underline{\mathbf W}^T * \underline{ \Omega}_i + \underline{\mathbf W} * \underline{\mathbf Z}^{-1} * \underline{\mathbf W}^T * \underline{ \Omega}_i,
\\
&=& \underline{\mathbf X}^T * \underline{\mathbf X} *\underline{ \Omega}_i - \underline{\mathbf W} * \underline{\mathbf Z}^{-1} * \underline{\mathbf W}^T *\underline{ \Omega}_i.
\end{eqnarray*}
So, using these changes, we can modify Algorithm \ref{ALgRR} to Algorithm \ref{fixed-precis2}. Following \cite{feng2023fast}, one step of the orthonormalization has been skipped to reduce the computing time with a slight accuracy drop.

\subsection{Complexity and Performance Analysis}
Let $\X \in \mathbb{R}^{I_1 \times I_2 \times I_3}$ with an adaptive rank $R(\epsilon)$. The computational complexity of the proposed fixed-precision algorithms is summarized in Table \ref{tab:complexities}, where $q$ is the number of power iterations/passes and $b$ is the block size.

\begin{table}[ht]
\centering
\caption{Computational complexity of fixed-precision algorithms.}
\begin{tabular}{|l|c|}
\hline
Algorithm & Computational Complexity \\
\hline
Algorithm \ref{ALgRR} & $O\big( I_1 I_2 I_3  R(\epsilon)(\log I_3 + q) \big)$ \\
Algorithm \ref{fixed-precision} & $O\big( I_1 I_2 I_3 R(\epsilon) (\log I_3 + \lfloor q/2 \rfloor) \big)$ \\
Algorithm \ref{fixed-precis2} & $O\big( I_1 I_2 I_3 (R(\epsilon)\log I_3 + q b) \big)$ \\
\hline
\end{tabular}
\label{tab:complexities}
\end{table}

\begin{algorithm}
\caption{Randomized fixed-precision algorithm \cite{ahmadi2023efficient}}\label{ALgRR}
\textbf{Input:} A tensor $\underline{\mathbf{X}}\in\mathbb{R}^{I_1\times I_2\times I_3}$; an error bound $\epsilon$; a block size $b$ and a power iteration $q$;\\
\textbf{Output:} $\underline{\mathbf{Q}}=[\underline{\mathbf{Q}}^{(1)},\underline{\mathbf{Q}}^{(2)},\ldots,\underline{\mathbf{Q}}^{(i)}],\,\underline{\mathbf{B}}=\begin{bmatrix} \underline{\mathbf{B}}^{(1)}\\ \underline{\mathbf{B}}^{(2)}\\ \vdots\\ \underline{\mathbf{B}}^{(i)} \end{bmatrix}$ such that ${\left\| {\underline{\mathbf{X}} - \underline{\mathbf{Q}}*\underline{\mathbf{B}}} \right\|_F} < \varepsilon$ and corresponding optimal tubal rank $R$;
\begin{algorithmic}[1]
\State $\underline{\mathbf{Q}}=[],\,\,\underline{\mathbf{B}}=[~{}]$;
\State $E = \left\| \underline{\mathbf{X}} \right\|_F^2$;
\For{$i=1,2,\ldots$}
\State ${\underline{\Omega}^{(i)}} = \mathrm{randn}\left( {I_2,b,I_3} \right)$;
\State ${\underline{\mathbf{Q}}^{(i)}} = \mathrm{orth}\left( {\underline{\mathbf{X}}*{\underline{\Omega}^{(i)}} - \underline{\mathbf{Q}}*\left( {\underline{\mathbf{B}}*{\underline{\Omega}^{(i)}}} \right)} \right)$;
\For{$j=1,2,\ldots, q$}
\State ${\underline{\mathbf{Q}}^{(i)}} = \mathrm{orth}\left({{\underline{\mathbf{X}}}^T*{\underline{\mathbf{Q}}^{(i)}}} -\underline{\mathbf{B}}^T*\underline{\mathbf{Q}}^T*\underline{\mathbf{Q}}^{(i)}\right)$;
\State ${\underline{\mathbf{Q}}^{(i)}} = \mathrm{orth}\left( {\underline{\mathbf{X}}*{\underline{\mathbf{Q}}^{(i)}}} -\underline{\mathbf{Q}}*\underline{\mathbf{B}}*\underline{\mathbf{Q}}^{(i)}\right)$;
\EndFor
\State ${\underline{\mathbf{Q}}^{(i)}} = \mathrm{orth}\left( {{\underline{\mathbf{Q}}^{(i)}} - \underline{\mathbf{Q}}*\left( {{\underline{\mathbf{Q}}^T}*{\underline{\mathbf{Q}}^{(i)}}} \right)} \right)$;
\State ${\underline{\mathbf{B}}^{(i)}} = \underline{\mathbf{Q}}^{(i)T}*\underline{\mathbf{X}}$;
\State $\underline{\mathbf{Q}}=\underline{\mathbf{Q}} \boxplus_2{\underline{\mathbf{Q}}^{(i)}}$;
\State $\underline{\mathbf{B}}=\underline{\mathbf{B}} \boxplus_1 {\underline{\mathbf{B}}^{(i)}}$;
\State $E = E - \left\| {{\underline{\mathbf{B}}^{(i)}}} \right\|_F^2$;
\If{$E < {\varepsilon}^2$}
\State \textbf{break}
\EndIf
\EndFor
\end{algorithmic}
\end{algorithm}

\begin{algorithm}
\caption{The proposed randomized fixed-precision algorithm I}\label{fixed-precision}
\textbf{Input:} The data tensor $\underline{\mathbf{X}} \in {\mathbb{R}^{{I_1} \times {I_2} \times {I_3}}}\,(I_1\leq I_2$), an approximation error bound $\epsilon$, the pass number $q>2$ and block size $b$;\\
\textbf{Output:} A low tubal rank approximation of the tensor $||\underline{\mathbf{X}}-\underline{\mathbf{U}}*\underline{\mathbf{S}}*\underline{\mathbf{V}}^T||\leq \epsilon$ and $\underline{\mathbf{U}}\in\mathbb{R}^{I_1\times k\times I_3},\,$ $\underline{\mathbf{S}}\in\mathbb{R}^{k\times k\times I_3}$ and $\underline{\mathbf{V}}\in\mathbb{R}^{I_2\times k\times I_3}$;
\begin{algorithmic}[1]
\State $\underline{\mathbf{Q}}=[\,];\quad\underline{\mathbf{B}}=[\,]$;
\For{$l=1,2,\ldots$}
\If{$q$ is an even number}
\State $\underline{\Omega} = \mathrm{randn}\left( {I_2,b,I_3} \right)$;
\State $\underline{\mathbf{Y}}=\underline{\mathbf{X}}*\underline{\Omega}-\underline{\mathbf{Q}}*(\underline{\mathbf{B}}*\underline{\Omega})$;
\State $[\underline{\mathbf{Q}}_l,\sim]=\mathrm{T\rnumber LU}(\underline{\mathbf{Y}})$;
\Else
\State Set ${\underline{\mathbf{Q}}}_l$ as a random tensor of size $I_1\times b\times I_3$;
\EndIf
\For{$t=1,2,3,\ldots,\lfloor\frac{q-1}{2}\rfloor$}
\If{$t==\lfloor \frac{q-1}{2}\rfloor$}
\State $\underline{\mathbf{R}}=\underline{\mathbf{X}}^T*\underline{\mathbf{Q}}_l$;
\State $\underline{\mathbf{Q}}_l=\mathrm{orth}(\underline{\mathbf{X}}*\underline{\mathbf{R}}-\underline{\mathbf{Q}}*(\underline{\mathbf{B}}*\underline{\mathbf{R}}))$;
\Else
\State $[\underline{\mathbf{Q}}_l,\sim]=\mathrm{T\rnumber LU}(\underline{\mathbf{X}}*(\underline{\mathbf{X}}^T*\underline{\mathbf{Q}}_l))$;
\EndIf
\EndFor
\State $\underline{\mathbf{Q}}_l=\mathrm{orth}(\underline{\mathbf{Q}}_l-\underline{\mathbf{Q}}*(\underline{\mathbf{Q}}^T*\underline{\mathbf{Q}}_l))$;
\State $\underline{\mathbf{B}}_i=\underline{\mathbf{Q}}_l^T*\underline{\mathbf{X}}$;
\State $\underline{\mathbf{Q}}=\underline{\mathbf{Q}} \boxplus_2 \underline{\mathbf{Q}}_l$;
\State $\underline{\mathbf{B}}=\underline{\mathbf{B}} \boxplus_1 \underline{\mathbf{B}}_l$;
\If{$\|\underline{\mathbf{X}}-\underline{\mathbf{Q}}*\underline{\mathbf{B}}\|_F < \epsilon$}
\State Set the tubal rank as $k$ and \textbf{break}
\EndIf
\EndFor
\State $[\widehat{\underline{\mathbf{U}}},\widehat{\underline{\mathbf{S}}},\widehat{\underline{\mathbf{V}}}]=\mathrm{T\rnumber SVD}(\underline{\mathbf{B}}^T)$;
\State $ind=lb:-1:lb-k+1$;
\State $\underline{\mathbf{U}}=\underline{\mathbf{Q}}*\widehat{\underline{\mathbf{V}}}(:,ind,:),\,\, \underline{\mathbf{S}}=\widehat{\underline{\mathbf{S}}}(ind,ind,:),\,\, \underline{\mathbf{V}}=\widehat{\underline{\mathbf{U}}}(:,ind,:)$;
\end{algorithmic}
\end{algorithm}

\begin{algorithm}
\caption{The proposed randomized fixed-precision algorithm II}\label{fixed-precis2}
\textbf{Input:} The data tensor $\underline{\mathbf{X}} \in {\mathbb{R}^{{I_1} \times {I_2} \times {I_3}}}$, a block size $b$, a power iteration $q$ and an approximation error bound $\epsilon$;\\
\textbf{Output:} The QB approximation of the tensor $\|\underline{\mathbf{X}}- \underline{\mathbf{Q}}*\underline{\mathbf{B}}\|_F\leq\epsilon$;
\begin{algorithmic}[1]
\State $\underline{\mathbf{Y}}=[\,],\quad\underline{\mathbf{W}}=[\,]$;
\State $E=\|\underline{\mathbf{X}}\|_F^2,\,{tol}=\epsilon^2$;
\For{$i=1,2,\ldots$}
\State Generate a random tensor $\underline{\Omega}_i$ of size $I_2\times b\times I_3$;
\For{$j=1,2,\ldots,q$}
\State $\underline{\mathbf{W}}_i=\underline{\mathbf{X}}^T*\underline{\mathbf{X}}*\underline{\Omega}_i-\underline{\mathbf{W}}*\underline{\mathbf{Z}}^{-1}*\underline{\mathbf{W}}^T*\underline{\Omega}_i$;
\State $\underline{\Omega}_i=\mathrm{orth}(\underline{\mathbf{W}}_i)$;
\EndFor
\State $\underline{\mathbf{Y}}_i=\underline{\mathbf{X}}*\underline{\Omega}_i,\quad \underline{\mathbf{W}}_i=\underline{\mathbf{X}}^T*\underline{\mathbf{Y}}_i$;
\State $\underline{\mathbf{Y}}= \underline{\mathbf{Y}}\boxplus_2\underline{\mathbf{Y}}_i,\quad \underline{\mathbf{W}}=\underline{\mathbf{W}}\boxplus_2\underline{\mathbf{W}}_i$;
\State $\underline{\mathbf{Z}}=\underline{\mathbf{Y}}^T*\underline{\mathbf{Y}},\quad \underline{\mathbf{T}}=\underline{\mathbf{W}}^T*\underline{\mathbf{W}}$;
\If{$E-(\mathrm{trace}((\underline{\mathbf{T}}*\underline{\mathbf{Z}}^{-1})^{(1)}))<tol$}
\State \textbf{break}
\EndIf
\EndFor
\State $[\widehat{\underline{\mathbf{V}}},\widehat{\underline{\mathbf{D}}}]=\mathrm{T\rnumber EIG}(\underline{\mathbf{Z}})$;
\State $\underline{\mathbf{Q}}=\underline{\mathbf{Y}}*\widehat{\underline{\mathbf{V}}}*\mathrm{sqrt}(\widehat{\underline{\mathbf{D}}})^{-1}$;
\State $\underline{\mathbf{B}}=(\underline{\mathbf{W}}*\widehat{\underline{\mathbf{V}}}*\mathrm{sqrt}(\widehat{\underline{\mathbf{D}}})^{-1})^T$;
\end{algorithmic}
\end{algorithm}

\subsection{Practical Decision Framework for Algorithm Selection}
The proposed algorithms address different operational scenarios in large-scale tensor computations. The choice of algorithm depends on three key factors: (1) data accessibility (can data be revisited?), (2) prior knowledge of the target rank, and (3) computational constraints (time, memory, and parallel resources). So, the decision workflow can be presented as follows:

\begin{enumerate}
    \item \textbf{Can the data be stored and revisited?}
    \begin{itemize}
        \item \textbf{No} (streaming data, limited memory): Use single-pass algorithms (Algorithms \ref{Single-1}, \ref{Single-2}, or \ref{Single-3}).
        \item \textbf{Yes}: Proceed to step 2.
    \end{itemize}
    
    \item \textbf{Is the target tubal rank known a priori?}
    \begin{itemize}
        \item \textbf{No}: Use fixed-precision algorithms (Algorithms \ref{fixed-precision} or \ref{fixed-precis2}).
        \item \textbf{Yes}: Use multi-pass randomized T-SVD (rT-SVD) for higher accuracy, or single-pass algorithms if data cannot be revisited.
    \end{itemize}
    
    \item \textbf{What are the computational constraints?}
    \begin{itemize}
        \item \textbf{Limited time}: Use T-LU based algorithms (Algorithm \ref{fixed-precision}) for 25--30\% speedup.
        \item \textbf{Limited memory}: Use single-pass algorithms (Algorithms \ref{Single-1}-\ref{Single-3}).
        \item \textbf{Parallel resources available}: Use T-LU for better parallelization.
    \end{itemize}
\end{enumerate}

The typical use cases of the algorithms are described below:

\begin{itemize}
    \item \textbf{Video surveillance (streaming):} Algorithm \ref{Single-2} (single-pass, stabilized) is ideal for real-time video processing where data cannot be stored.
    \item \textbf{Image compression (batch processing):} Algorithm \ref{fixed-precision} or \ref{fixed-precis2} (fixed-precision) for optimal rank selection.
    \item \textbf{Scientific computing (known rank):} Multi-pass rT-SVD for maximum accuracy, or Algorithm \ref{Single-1} for single-pass approximation.
\end{itemize}

\begin{figure*}[htbp]
\begin{center}
\includegraphics[width=1.05\linewidth]{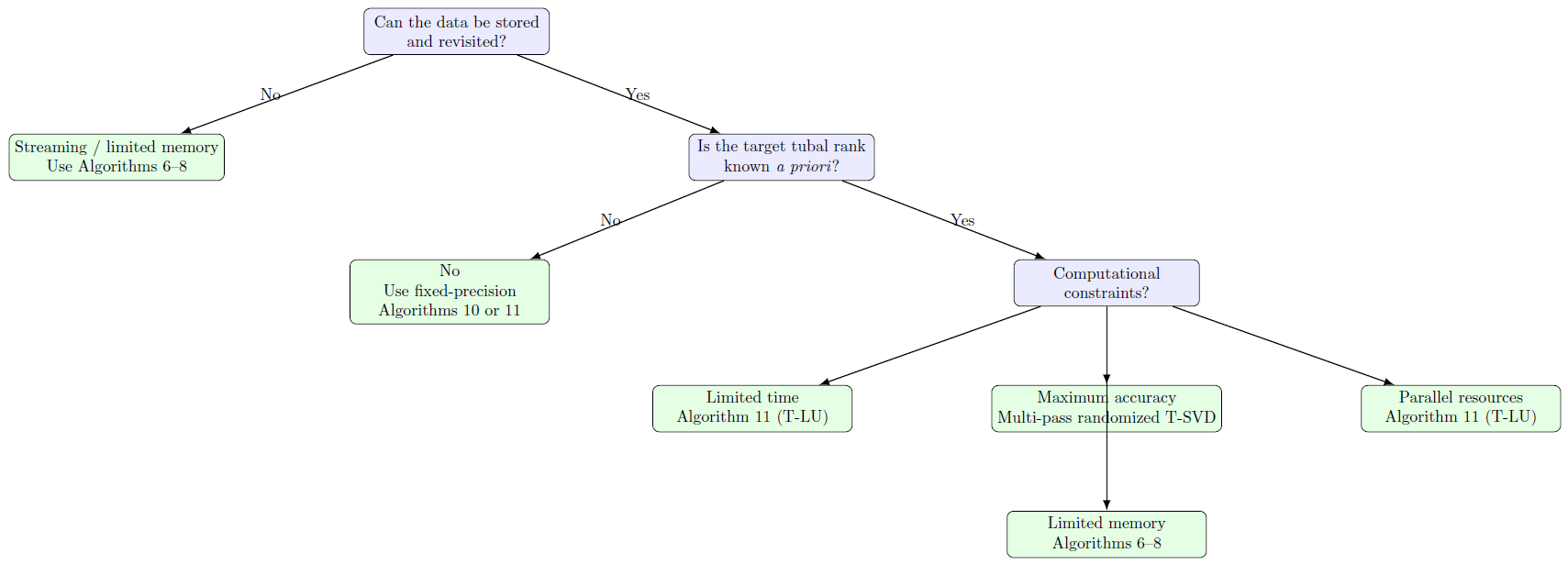}
\caption{Decision workflow for selecting randomized tensor decomposition algorithms based on data accessibility, prior knowledge of the tubal rank, and computational constraints..}
\label{Workflow}
\end{center}
\end{figure*}

Figure \ref{Workflow} provides a visual guide to the algorithm selection process outlined above.

\section{Experimental Results}\label{Sec:Experi}
In this section, we evaluate the performance of the proposed algorithms on both synthetic and real-world data tensors. All experiments were conducted using MATLAB on a computer equipped with a 2.60 GHz Intel(R) Core(TM) i7-5600U processor and 8 GB of memory.

The first experiment uses synthetic data. The second and third experiments focus on image and video compression tasks, respectively. The final two experiments demonstrate the application of the proposed approaches to image super-resolution and deep learning tasks.

The PSNR of two images $\underline{\mathbf{X}}$ and $\underline{\mathbf{Y}}$ is defined as
\[
\mathrm{PSNR}=10\log _{10}\left({\frac{255^2}{\mathrm{MSE}}}\right)\,\mathrm{dB},
\]
where $\mathrm{MSE}=\frac{\|\underline{\mathbf{X}}-\underline{\mathbf{Y}}\|_F^2}{\mathrm{num}(\underline{\mathbf{X}})}$ and ``num($\underline{\mathbf{X}}$)'' stands for the number of elements of the data tensor $\underline{\mathbf{X}}$.
The relative error is also defined as
\[
e(\widetilde{\X})=\frac{\|\X-\widetilde{\X}\|_F}{\|\X\|_F},
\]
where $\X$ is the original tensor and $\widetilde{\X}$ is the approximated tensor.

\begin{Example}\label{Exa_1}
({\bf Synthetic data tensors}) 
In this example, we compare the efficiency of the proposed algorithms with baseline methods on synthetic tensor data. We generate a random tensor using a low-rank approximation. First, we construct a noiseless tensor $\X\in\mathbb{R}^{I_1\times I_2\times I_3}$ of tubal rank $R$ as follows:
\begin{eqnarray}\label{ranten}
\X_{\mathrm{clean}}=\mathrm{randn}(I_1,R,I_3)*\mathrm{randn}(R,I_2,I_3).
\end{eqnarray}
We then add noise to obtain a perturbed tensor $\X_{\mathrm{perturb}}=\X_{\mathrm{clean}}+\delta\frac{\Y}{\|\Y\|_F}\|\X_{\mathrm{clean}}\|_F$, where $\Y$ is a standard Gaussian tensor of the same size as $\X$. In all simulations, we set $R=50$, $\delta=10^{-3}$, and assume $I_1=I_2=I_3$. Consequently, the tubal rank of $\X_{\mathrm{perturb}}$ is approximately $50$. We then evaluate how accurately the proposed fixed-precision algorithms estimate this tubal rank.

For an error tolerance of $\epsilon=10^{-5}$ and a block size of $b=100$, we apply the proposed fixed-precision algorithms, along with the baseline methods (including Algorithm \ref{ALgRR} and the truncated T-SVD), to the perturbed tensor $\X_{\mathrm{perturb}}$. Note that Algorithms \ref{ALgRR}--\ref{fixed-precis2} automatically estimate the tubal rank, whereas the truncated T-SVD requires the tubal rank as an input. To address this issue, we first estimate the tubal rank using Algorithms \ref{ALgRR}--\ref{fixed-precis2} and then use the estimated rank as input to the truncated T-SVD. Our first observation is that while all fixed-precision algorithms correctly estimate the tubal rank, the proposed algorithms require significantly less CPU time than the baseline Algorithm \ref{ALgRR}.

Additionally, for tubal rank $R=50$, we applied the truncated T-SVD to the tensor $\X_{\mathrm{perturb}}$. The running time of the truncated T-SVD was significantly higher than that of the proposed fixed-precision algorithms, while its accuracy was nearly identical to that of the proposed randomized algorithms. We also observe that as the tensor size increases, the running times of the proposed algorithms become much shorter than those of the baselines. These observations confirm that the proposed fixed-precision algorithms are faster and more efficient than the baseline methods.

Table \ref{tab:resnet_inference} reports the simulation results, including the relative error and computation time, for the proposed randomized fixed-precision algorithms and the baseline methods (Algorithm \ref{ALgRR} and the truncated T-SVD) with $I_1=200,300,400$ and tubal rank $R=50$. The results clearly show that the proposed randomized fixed-precision algorithms achieve more accurate approximations than the baselines in significantly less computation time.

\begin{table*}[htbp]
\centering
\caption{Comparison of running time (seconds) and relative error for fixed-precision algorithms (\ref{ALgRR},\ref{fixed-precision},\ref{fixed-precis2}) and deterministic truncated T-SVD on synthetic tensors of increasing size $n = 200$--$500$ with tubal rank $50$ (Example \ref{Exa_1}).}
\label{tab:resnet_inference}
\begin{tabular}{|l | c | c | c | c|}
\hline
Tensor size & Algorithm \ref{ALgRR} &  Algorithm \ref{fixed-precision} &  Algorithm \ref{fixed-precis2} & Truncated T-SVD\\\hline
$n=200$ & (2.94,\,2.95e-10) & (2.71,\,2.58e-09) & (1.18,\,4.72e-09) &  (11.43,\,1.34e-09)   \\\hline
$n=300$ & (8.74,\,6.47e-10) & (4.21,\,6.95e-09) & (3.07,\,9.20e-09) &  (36.81,\,1.17e-09) \\\hline
$n=400$ & (20.86,\,1.15e-09) & (7.98,\,1.27e-08) & (6.65,\,1.63e-08) &  (81.83,\,2.11e-09) \\\hline
$n=500$ & (45.89,\,1.90e-09) & (20.59,\,2.34e-08) &  (19.32,\,1.61e-08) & (195.25,\,3.43e-09)\\\hline
\end{tabular}
\end{table*}

Furthermore, to assess the performance of the proposed randomized single-pass Algorithms~\ref{crosstensor}--\ref{Single-3}, we applied them to the tensor in \eqref{ranten} with $I_1=300$ and $R=50$, and compared them with the baseline Algorithms~\ref{crosstensor}--\ref{Single-3}. For Algorithms~\ref{Single-1}--\ref{Single-3}, we used the sketch parameters $L=K=50$, $H=45$, and $R=40$, while for Algorithms~\ref{crosstensor}--\ref{QIsinglepass}, we used $L=K=40$. Thus, all algorithms compute a low-rank approximation of tubal rank $R=40$. The results are shown in Table~\ref{tab:sing_pass}.

As presented in Table~\ref{tab:sing_pass}, the running time of the proposed single-pass algorithms is slightly higher than that of the baselines. However, they are more stable in selecting sketch sizes. We will see in Example~\ref{im_com} that the sensitivity of Algorithms~\ref{crosstensor}--\ref{QIsinglepass} to the sketch parameters $L=K$ is more serious, and the proposed algorithms are more reliable.

\begin{table}[htbp]
\centering
\caption{Computing time (seconds) and relative error for Algorithms \ref{crosstensor}--\ref{Single-3} on a $300\times300\times300$ synthetic tensor with tubal rank $40$ (Example~1).}
\label{tab:sing_pass}
\begin{tabular}{| c| c| c| c|c|}
\hline
Algorithm \ref{crosstensor} & Algorithm \ref{QIsinglepass} & Algorithm \ref{Single-1} & Algorithm \ref{Single-2} & Algorithm \ref{Single-3} \\\hline
(2.96,\,5.75) & (4.97,\,8.10) & (9.60,\,0.26) & (9.01,\,0.26) & (12.18,\,0.26)\\\hline
\end{tabular}
\end{table}

To further evaluate the performance of the proposed algorithms, we consider three new synthetic data tensors defined as follows:

\begin{itemize}
    \item \textbf{Tensor Case I:} $\underline{\mathbf{X}}(i,j,k)=\dfrac{1}{\sqrt{i^2+j^2+k^2}}$
    \item \textbf{Tensor Case II:} $\underline{\mathbf{Y}}(i,j,k)=\dfrac{1}{(i^3+j^3+k^3)^{1/3}}$
    \item \textbf{Tensor Case III:} $\underline{\mathbf{Z}}(i,j,k)=\dfrac{1}{\sin(i)+\tanh(j+k)}$
\end{itemize}

We applied the proposed single-pass algorithms and the baseline Algorithms~\ref{crosstensor}--\ref{QIsinglepass} to the tensors described above, each of size $300 \times 300 \times 300$, with tubal rank $R=40$. The numerical results, reported in Table~\ref{tab:synthe_}, demonstrate the robustness of the proposed algorithms compared to the baselines. These simulations confirm that the proposed algorithm is more efficient than the other techniques.

\begin{table*}[htbp]
\centering
\caption{Computing time and relative error for Algorithms \ref{crosstensor}--\ref{Single-3} on three synthetic data tensors with non-polynomial decaying entries (Cases I--III).}
\label{tab:synthe_}
\begin{tabular}{|c| c| c| c|c|}
\hline
\multicolumn{5}{c}{Data Tensor Case I}\\\hline
Algorithm \ref{crosstensor} & Algorithm \ref{QIsinglepass} & Algorithm \ref{Single-1} & Algorithm \ref{Single-2} & Algorithm \ref{Single-3} \\\hline
(1.6,\,0.07) &   (3.6,\,2.81e-14) &  (7.7,\,1.91e-14) &  (7.6,\,1.26e-14) & (6.3,\,3.04e-14)\\\hline
\multicolumn{5}{c}{Data Tensor Case II}\\\hline
Algorithm \ref{crosstensor} & Algorithm \ref{QIsinglepass} & Algorithm \ref{Single-1} & Algorithm \ref{Single-2} & Algorithm \ref{Single-3} \\\hline
(1.5,\,0.014) &  (3.10,\,3.62e-12) &  (7.00,\,2.80e-14) &  (7.45,\,5.79e-14) & (9.32,\,2.92e-14)\\\hline
\multicolumn{5}{c}{Data Tensor Case III}\\\hline
Algorithm \ref{crosstensor} & Algorithm \ref{QIsinglepass} & Algorithm \ref{Single-1} & Algorithm \ref{Single-2} & Algorithm \ref{Single-3} \\\hline
(1.62,\,2.60) &   (3.53,\,8.50e-15) &  (7.35,\,1.91e-14) &  (7.54,\,2.36e-14) & (6.23,\,5.19e-14)\\\hline
\end{tabular}
\end{table*}
\end{Example}

\begin{Example}\label{im_com}
({\bf Image Compression})
In this experiment, we evaluate the effectiveness of the proposed randomized single-pass algorithms for image compression. We used the Kodak dataset and considered four images in our simulations: "Kodim15", "Kodim17", "Kodim18", and "Kodim23". The first two images are of size $512\times 768 \times 3$, while the latter two are of size $768\times 512\times 3$. We applied our proposed algorithms and compared them with the baseline methods, namely single-pass T-CUR \cite{tarzanagh2018fast} and tensor sketch \cite{qi2021t}. In our simulations, we set the sketch parameters to $L=350$, $K=350$, $H=100$, and $R=30$. The reconstructed images, along with their corresponding PSNR and relative error values, are shown in Figure~\ref{Res_1} and Table~\ref{tab:kodak_comparison}. The results demonstrate that our proposed randomized single-pass algorithms outperform both single-pass T-CUR and the tensor sketch algorithms.

We observed experimentally that when $K=L$, all algorithms except our proposed randomized single-pass Algorithms~\ref{Single-1}--\ref{Single-3} were unstable and yielded poor results. In this regard, our Algorithm~\ref{Single-2} proves to be more reliable. Better results are achieved when the sketch parameters are different, i.e., $L \neq K$. Another important numerical finding is that for $L < K$, the results were consistently poor across all experiments. This instability arises from the ill-conditioning of the underlying least-squares problem. More precisely, the Moore-Penrose pseudoinverse of the coefficient tensor in the least-squares formulation is poorly computed, leading to unsatisfactory results. In contrast, the proposed randomized single-pass Algorithms~\ref{Single-1}--\ref{Single-3} are robust because the truncated T-SVD acts as a regularizer, improving numerical stability. This illustration clearly demonstrates that the proposed approach is more reliable and capable of producing superior results.

\begin{figure*}[htbp]
\begin{center}
\includegraphics[width=0.6\linewidth]{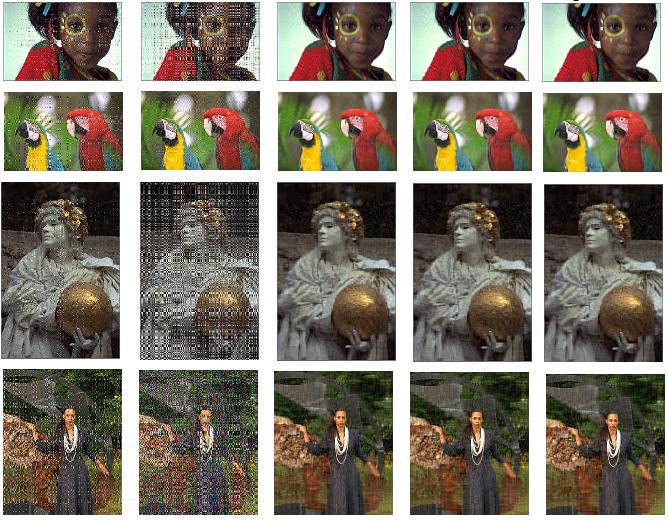}
\caption{Reconstructed images from the Kodak dataset (Kodim15, Kodim17, Kodim18, Kodim23) using different single-pass algorithms with sketch sizes $L = 350$, $K = 350$, truncation parameter $H = 100$ and target tubal rank $R = 30$. Columns from left to right: Algorithm \ref{crosstensor} (cross tensor approximation), Algorithm \ref{QIsinglepass} (tensor sketch), Algorithm \ref{Single-1} (proposed single-pass I), Algorithm \ref{Single-2} (proposed single-pass II), Algorithm \ref{Single-3} (proposed two-sided single-pass).}
\label{Res_1}
\end{center}
\end{figure*}

\begin{table*}[htbp]
\centering
\small
\caption{PSNR (dB) and relative error for Kodak images using single-pass algorithms with $L = K = 350$, $H = 100$, $R = 30$.}
\label{tab:kodak_comparison}
\begin{tabular}{lccccc}
\toprule
\textbf{Image} & \textbf{Alg. 5} & \textbf{Alg. 6} & \textbf{Alg. 7} & \textbf{Alg. 8} & \textbf{Alg. 9} \\
\midrule
Kodim 15 & 16.62 / 0.278 & 13.82 / 0.384 & 27.21 / 0.082 & 27.04 / 0.084 & 27.21 / 0.082 \\
Kodim 17 & 13.40 / 0.465 & 15.89 / 0.389 & 29.62 / 0.072 & 29.53 / 0.073 & 29.62 / 0.072 \\
Kodim 18 & 12.36 / 0.679 & 13.09 / 0.625 & 26.55 / 0.138 & 26.37 / 0.135 & 26.53 / 0.133 \\
Kodim 23 & 9.02 / 0.615 & 12.03 / 0.571 & 23.50 / 0.216 & 23.26 / 0.222 & 23.49 / 0.216 \\
\bottomrule
\end{tabular}
\par\medskip
\footnotesize{Note: Values shown as PSNR (dB) / Relative Error. Alg. = Algorithm.}
\end{table*}

\end{Example}
 
\begin{Example}
({\bf Video compression}) In this example, we study the performance of the proposed randomized single-pass algorithms for the video compression task. We primarily use the video datasets ``Foreman'' and ``News''. The videos are third-order tensors of size $144\times 176\times 300$. 

First, we evaluated the efficiency of the single-pass algorithms for low-rank tubal approximations on the aforementioned videos, using sketch parameters $K=90$, $L=90$, $H=20$, and $R=20$. The PSNR values for all frames of the Foreman and News videos achieved by the proposed single-pass algorithms and the baselines are reported in Figures~\ref{psnr_video_2} and~\ref{psnr_video}, respectively. Additionally, reconstructed images from selected frames of these videos are shown in Figures~\ref{video_recon_2} and~\ref{video_recon_1}, respectively. 

Our results demonstrate that the proposed Algorithms~\ref{Single-1}--\ref{Single-3} are robust with respect to sketch size selection. This example illustrates the reliability of the randomized single-pass algorithms for video compression tasks.

\begin{figure*}[htbp]
\begin{center}
\includegraphics[width=1\linewidth]{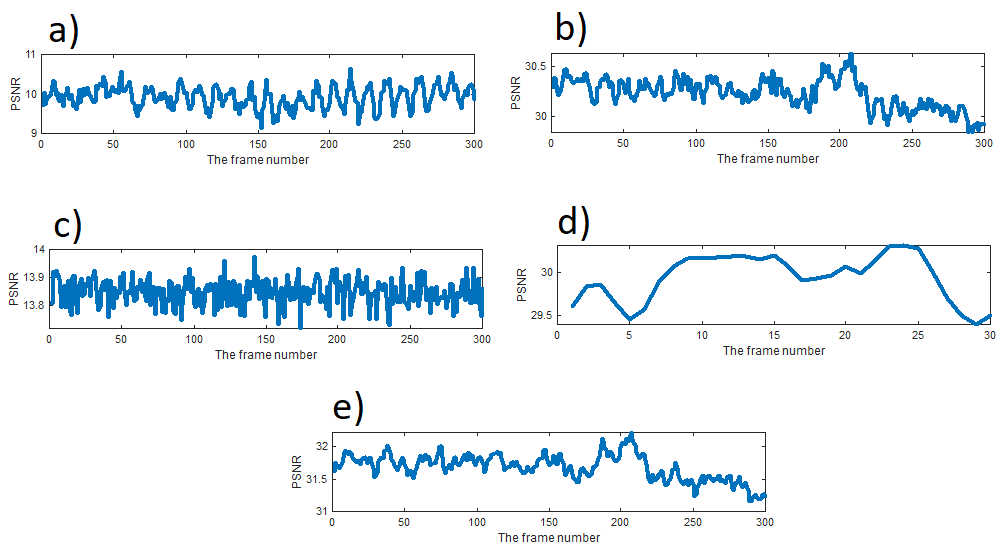}
\caption{PSNR (dB) of all 300 frames of the Foreman video (size $144\times176\times300$) achieved by different single-pass algorithms with $L = 90$, $K = 90$, $H = 20$, $R = 50$. (a) Algorithm \ref{crosstensor}, (b) Algorithm \ref{Single-1}, (c) Algorithm \ref{QIsinglepass}, (d) Algorithm \ref{Single-2}, (e) Algorithm \ref{Single-3}.}
\label{psnr_video_2}
\end{center}
\end{figure*}

\begin{figure*}[htbp]
\begin{center}
\includegraphics[width=1\linewidth]{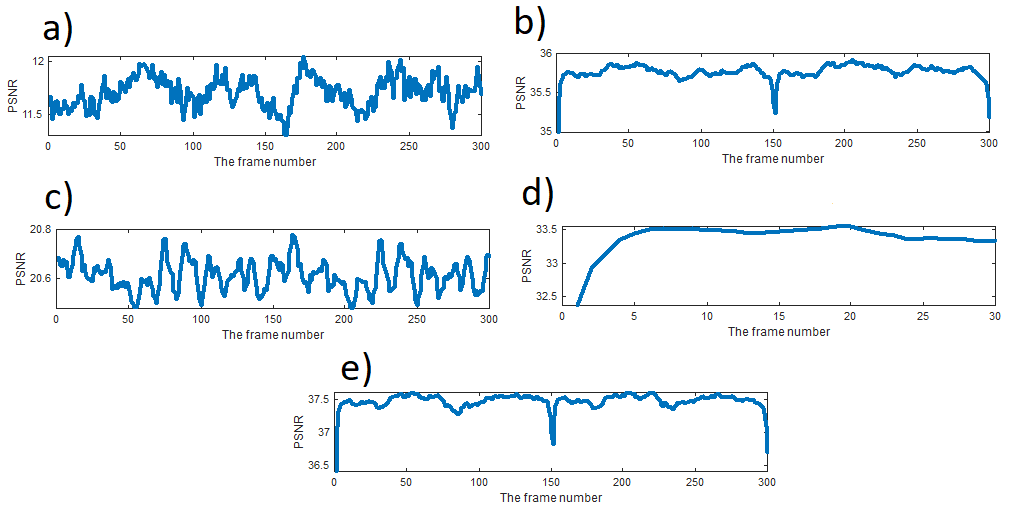}
\caption{PSNR (dB) of all frames of the News video (size $144\times176\times300$) using the same single-pass algorithms and parameters as in Figure~\ref{psnr_video_2}. (a) Algorithm \ref{crosstensor}, (b) Algorithm \ref{Single-1}, (c) Algorithm \ref{QIsinglepass}, (d) Algorithm \ref{Single-2}, (e) Algorithm \ref{Single-3}}
\label{psnr_video}
\end{center}
\end{figure*}

\begin{figure*}[htbp]
\begin{center}
\includegraphics[width=0.6\linewidth]{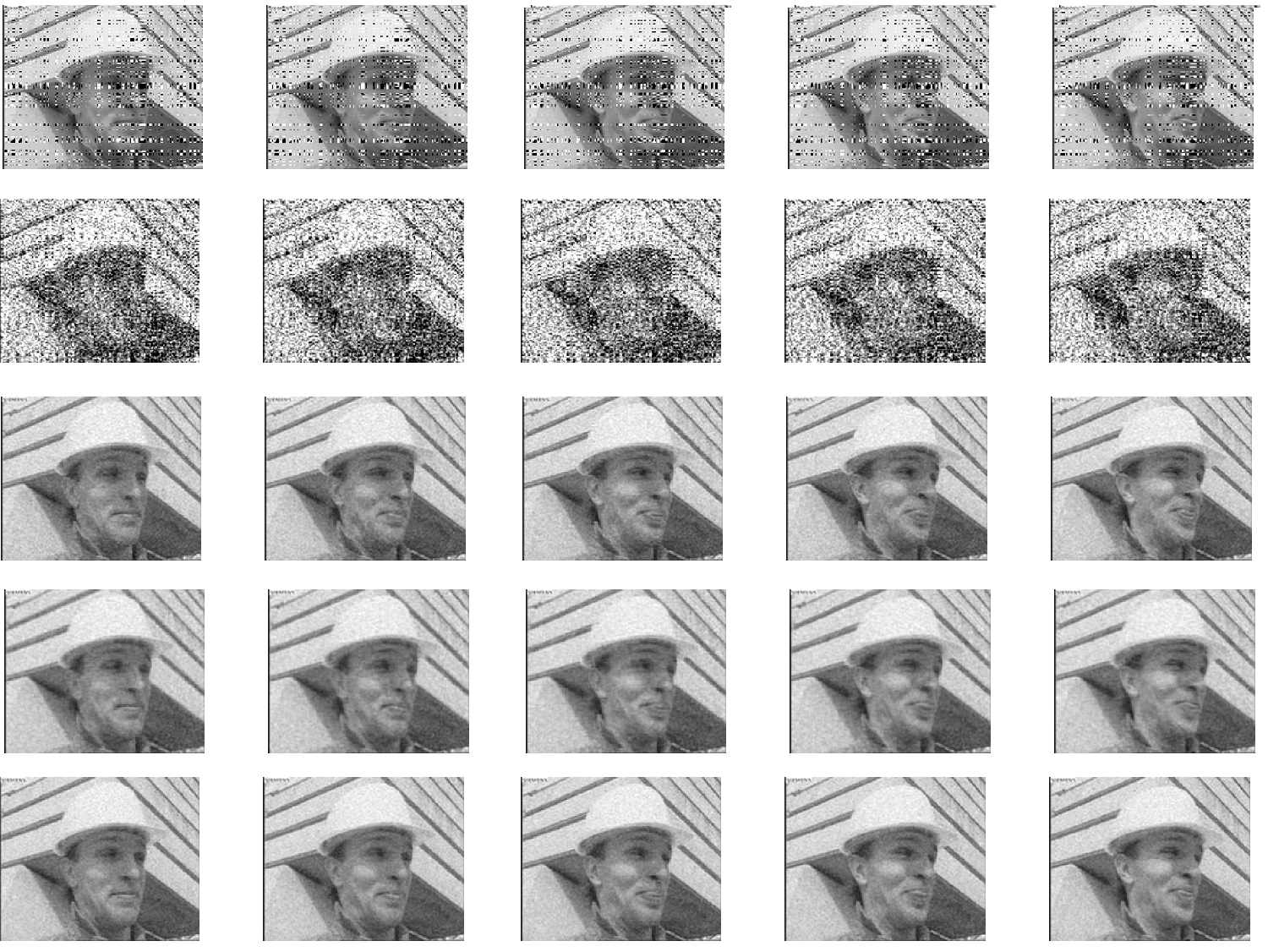}
\caption{Reconstruction of randomly selected frames from the Foreman video. Rows from top to bottom correspond to Algorithms \ref{crosstensor}, \ref{QIsinglepass}, \ref{Single-1}, \ref{Single-2}, \ref{Single-3} (same parameters as Figure~\ref{psnr_video_2}).}
\label{video_recon_2}
\end{center}
\end{figure*}

\begin{figure*}[htbp]
\begin{center}
\includegraphics[width=0.6\linewidth]{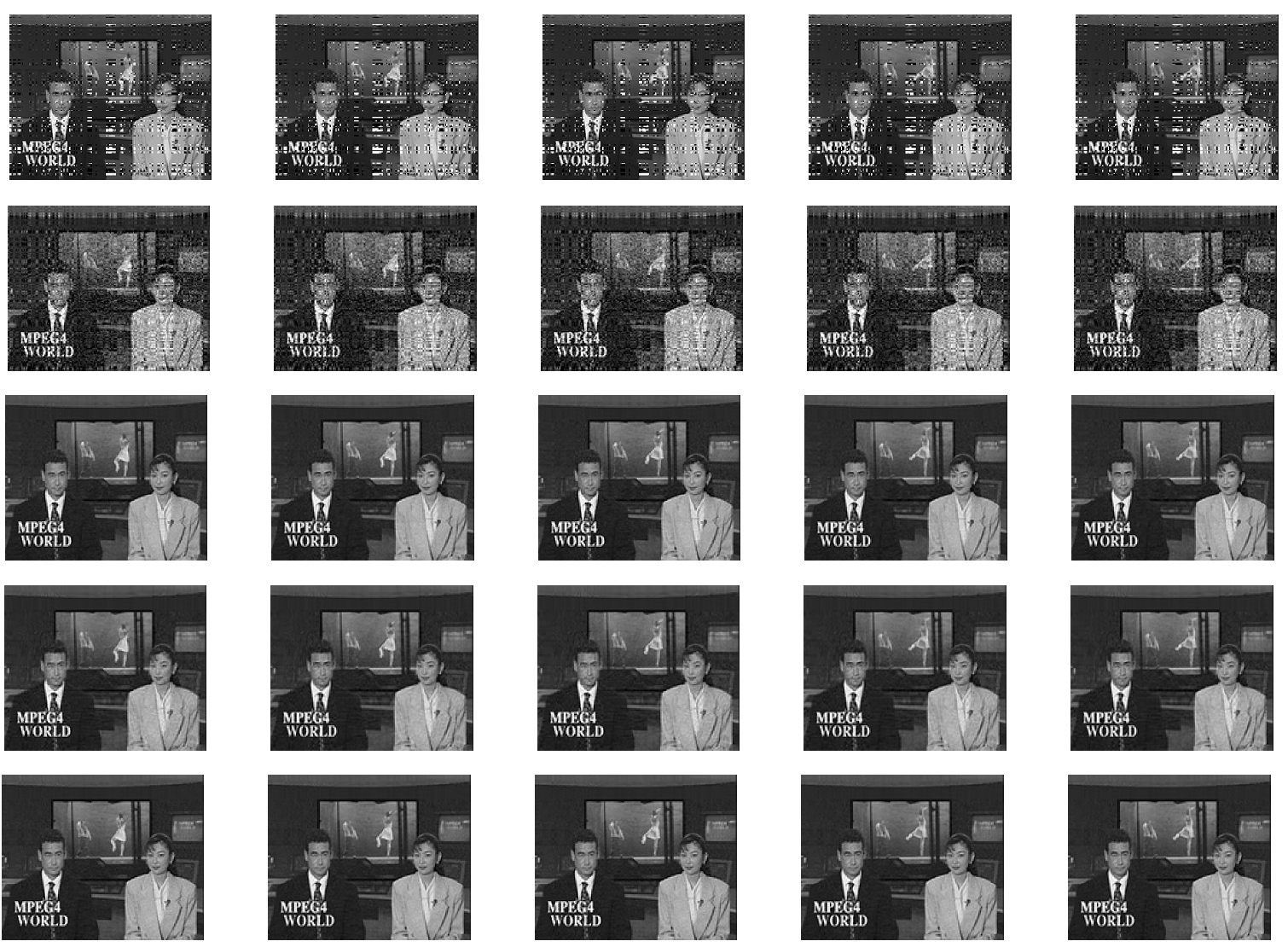}
\caption{Reconstruction of randomly selected frames from the News video. Rows from top to bottom correspond to Algorithms \ref{crosstensor}, \ref{QIsinglepass}, \ref{Single-1}, \ref{Single-2}, \ref{Single-3} (same parameters as Figure ~\ref{psnr_video_2}).}
\label{video_recon_1}
\end{center}
\end{figure*}
\end{Example}
 
\begin{Example}
({\bf Image Super-resolution}) \label{img_super} In this experiment, we study the efficiency and applicability of the proposed approaches to the super-resolution task. Super-resolution is an important computer vision task that involves reconstructing a higher-resolution image from a lower-resolution input. The goal is to produce a more aesthetically pleasing and detailed image while reducing artifacts and noise. Image super-resolution can be addressed using tensor completion. More precisely, a given small image is up-sampled along the first and second dimensions to obtain an incomplete image. Tensor completion is then applied to this incomplete image to recover the missing pixels. We adopt the tensor completion method proposed in \cite{ahmadi2023fast} for this task.

To the best of our knowledge, this is the first paper to use single-pass algorithms for tensor completion and image super-resolution. The following iterative steps are performed for $n=0,1,2,\ldots$:
\begin{eqnarray}\label{Step1}
\underline{\mathbf{X}}^{(n)} &\leftarrow& \mathcal{L}(\underline{\mathbf{C}}^{(n)}),\\
\label{Step2}
\underline{\mathbf{C}}^{(n+1)} &\leftarrow& \underline{\mathbf{C}}^{(0)} + (\underline{\mathbf{1}} - \underline{\Omega}) \oast \underline{\mathbf{X}}^{(n)},
\end{eqnarray}
where $\underline{\mathbf{C}}^{(0)}$ is the initial image with missing elements $\underline{\mathbf{M}}$, and $\mathcal{L}$ is an operator that computes a low-rank tensor approximation of the tensor $\underline{\mathbf{C}}^{(0)}$. Here, $\oast$ denotes the Hadamard (element-wise) product of two tensors, $\underline{\mathbf{1}}$ is a tensor of the same size as the original image with all elements equal to one, and $\underline{\Omega}$ is a binary mask tensor that indicates the locations of known elements: if an element in the tensor is known, the corresponding entry in $\underline{\Omega}$ is 1; otherwise, it is 0. 

The stages \eqref{Step1} and \eqref{Step2} are repeated until a specified approximation error bound is satisfied or a maximum number of iterations (80) is reached. The first stage is computationally expensive, so we employ our proposed single-pass approaches to compute a low-rank tensor approximation. Specifically, we use the proposed randomized single-pass algorithms. To further improve image recovery after the second stage \eqref{Step2}, we apply a Gaussian filter (using the MATLAB function \texttt{imgaussfilt}). 

Five images are considered: ``Peppers'', ``Airplane'', ``Kodim01'', ``Kodim02'', and ``Kodim03''. The first two images are of size $256\times 256\times 3$, while the remaining three are of size $512 \times 768 \times 3$. The images are up-sampled four times along both the $x$- and $y$-axes, and a tubal rank of $R=60$ is used for all images.

The simulation results are reported in Table~\ref{tab:images}, and two reconstructed images are shown in Figure~\ref{test}. From Table~\ref{tab:images}, we observe that the proposed randomized fixed-precision algorithms produce recovered images comparable to those of the deterministic method (truncated T-SVD) but with significantly less computational time. This demonstrates the advantage of the proposed randomized single-pass algorithms for image super-resolution.

To further benchmark the efficiency of our proposed single-pass algorithms, we compared them against the single-pass algorithm from \cite{sun2020low} used as the operator $\mathcal{L}$ in the first step of the completion process \eqref{Step1}. That baseline, which is based on the Tucker decomposition with a multilinear rank of $(60,60,3)$, was approximately twice as slow as our method. Our algorithm achieved comparable PSNR values for all test images while requiring half the computation time. A key limitation of this baseline is its large parameter count, which makes fine-tuning challenging. By contrast, our approaches exhibit lower sensitivity to parameter selection and are consequently easier to deploy.

\begin{table}[htbp]
\centering
\caption{Running time (seconds) and PSNR (dB) for image super-resolution (Example~4) using single-pass tensor completion vs. deterministic T-SVD completion.}
\label{tab:images}
\begin{tabular}{|l|c|c|}
\hline
Image & Single-pass completion &  Deterministic completion\\\hline
Peppers & (26.99, 22.01) & (42.70, 22.01) \\\hline
Airplane & (27.69, 22.13) & (44.58 ,22.11) \\\hline
Kodim01 & (50.34, 20.44) & (142.42, 20.56) \\\hline
Kodim02 & (44.58, 26.70) & (137.30, 26.91) \\\hline
Kodim03 & (48.05, 26.96) & (146.76, 27.42) \\\hline
\end{tabular}
\end{table}

\begin{figure}[htbp]
\begin{center}
\includegraphics[width=1\linewidth]{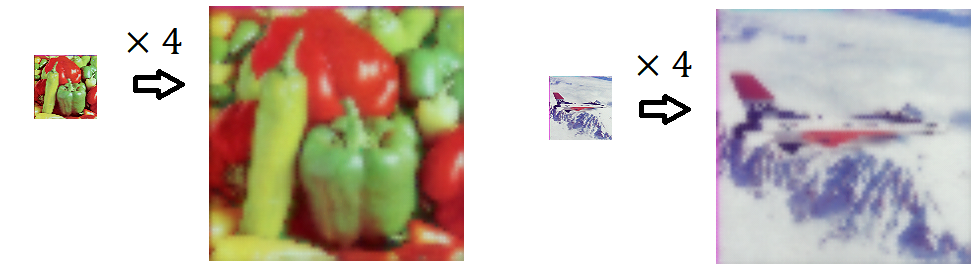}
\caption{Image super-resolution results using tensor completion with the proposed randomized single-pass algorithms for the low-tubal-rank operator $\mathcal{L}$ in Eq.~(50). Input low-resolution images (left column) are upsampled by a factor of $4$, and missing pixels are recovered via iterative completion. Output high-resolution images (right column) show restored details (e.g., edges, textures) comparable to deterministic T-SVD completion but with much lower computational time (see Table~\ref{tab:images}).}
\label{test}
\end{center}
\end{figure}
\end{Example}
 
\begin{Example}
({\bf Deep Learning Application})
In this example, we apply the proposed completion method to improve the accuracy of object detection, an important computer vision task in deep learning. Specifically, we consider two images (which we refer to as ``dog'' and ``horses'') shown in Figure~\ref{test_obdet} (first column) and manually remove some parts of the images in a structural manner, as depicted in Figure~\ref{test_obdet} (second column). Both images are of size $512\times 768 \times 3$.

We employ the pretrained YOLOv3 (\textbf{You Only Look Once}, Version 3), an efficient deep neural network for object detection \cite{redmon2018yolov3}. YOLOv3 improves upon its predecessor, YOLOv2 \cite{redmon2017yolo9000}, and aims to address several limitations of earlier versions. It is known for its speed and efficiency in real-time object detection tasks. YOLOv3 divides the input image into a grid and predicts bounding boxes and class probabilities for each grid cell, enabling faster processing than other object detection algorithms that scan the entire image multiple times. YOLOv3 has become a popular choice for real-time object detection tasks due to its speed, accuracy, and efficiency.

We then applied YOLOv3 to detect objects in the degraded images. The results are shown in Figure~\ref{test_obdet_bench}. The first column presents detection results for degraded images, whereas the second column presents those for the recovered images. As can be observed, for the dog image, the network detected only one object (a cat), which is also incorrect. For the horse image, the network detected two objects: one horse and one giraffe. It failed to detect two additional horses and incorrectly classified one as a giraffe.

In contrast, when our completion technique is applied as a preprocessing step, the degraded image is first processed quickly, after which YOLOv3 detects objects in the restored image. Here, we applied our single-pass algorithms with the tubal rank $R=150$, using 100 iterations, the Gaussian filtering parameter $0.6$, and $\epsilon=10^{-4}$ as a stopping tolerance. After this fast preprocessing with our techniques, the network detected three objects in the dog image: a bicycle, a dog, and a truck, all with very precise bounding boxes. Furthermore, for the horse image, the network successfully detected all four horses. 

These results clearly demonstrate the efficiency and feasibility of the proposed completion procedure for stabilizing YOLOv3 against pixel removal.

To further validate the robustness of our approach, we extend the experiment to YOLO11, one of the latest iterations in the YOLO series \cite{redmon2025yolo11}. YOLO11 introduces several architectural improvements over its predecessors, including:

\begin{itemize}
    \item \textbf{Enhanced backbone network:} A modified CSPNet (Cross Stage Partial Network) with improved gradient flow;
    \item \textbf{Attention mechanisms:} Integrated multi-head self-attention for better feature representation;
    \item \textbf{Adaptive anchor boxes:} Dynamic anchor box generation for improved localization;
    \item \textbf{Optimized inference:} Reduced latency with comparable or higher mean Average Precision (mAP).
\end{itemize}

We evaluate YOLO11 on the same degraded images under two configurations:

\begin{enumerate}
    \item \textbf{Baseline:} YOLO11 applied directly to degraded images (no preprocessing);
    \item \textbf{Proposed preprocessing:} Our tensor completion method followed by YOLO11;
\end{enumerate}

For the dog image, baseline YOLO11 partially fails (detects a cat instead of a dog). With our preprocessing, YOLO11 correctly identifies the dog, bicycle, and truck with high confidence scores ($>0.85$). For the horse image, our method enables YOLO11 to detect all four horses, while the baseline detects only two (one misclassified as a giraffe).

\begin{figure*}[htbp]
\begin{center}
\includegraphics[width=0.6\linewidth]{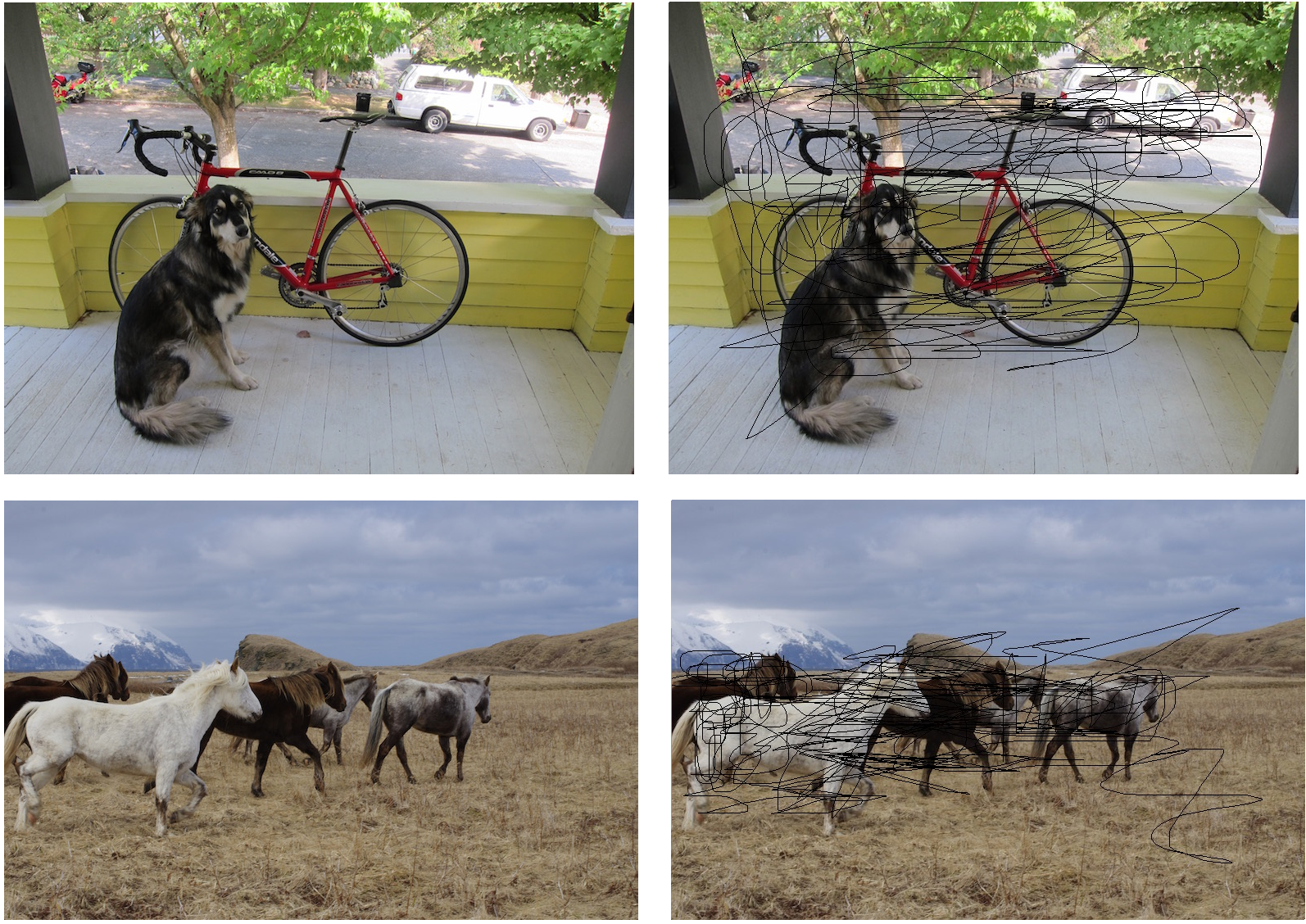}
\caption{Original (first column) and degraded (second column) images used in the object detection experiment (Example~5). Degraded images have manually removed structural parts (rectangular occlusions).}
\label{test_obdet}
\end{center}
\end{figure*}

\begin{figure*}[htbp]
\begin{center}
\includegraphics[width=0.6\linewidth]{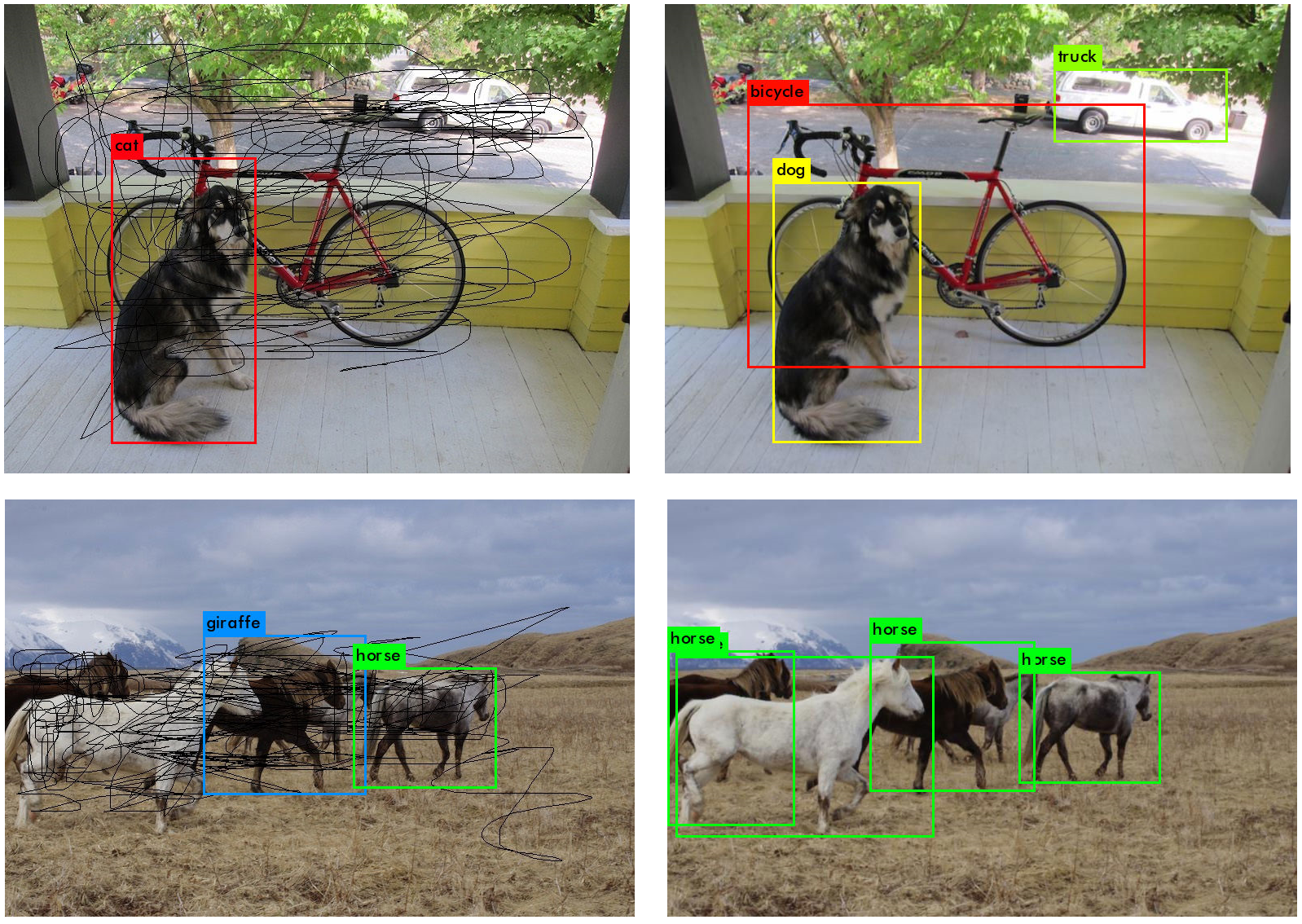}
\caption{Object detection results on degraded images (first column) vs. after recovery using the proposed tensor completion with randomized single-pass algorithms (second column).}
\label{test_obdet_bench}
\end{center}
\end{figure*}
\end{Example}

\begin{Example}
({\bf COIL-100 Data Set}) The Columbia Object Image Library (COIL-100) is a benchmark dataset widely used for object recognition and image classification tasks \cite{nene1996columbia}. It was curated by the Center for Research on Intelligent Systems at the Department of Computer Science, Columbia University.

The dataset contains \textbf{7,200} color images of \textbf{100 distinct objects}. For each object, images are captured at 5-degree pose intervals as the object rotates through a full 360 degrees on a motorized turntable. Consequently, each object has \textbf{72 images} representing different viewing angles \cite{huggingface_coil100, tensorflow_coil100}. The technical Specifications are as follows:

\begin{itemize}
    \item \textbf{Image Size}: Each image is resized and normalized to $128 \times 128$ pixels \cite{tensorflow_coil100}.
    \item \textbf{Color Format}: RGB color images (3 channels).
    \item \textbf{Background}: Uniform black background to simplify object segmentation.
    \item \textbf{File Format}: Typically stored as PNG or JPEG images.
    \item \textbf{Total Size}: Approximately 124.63 MB \cite{tensorflow_coil100}.
\end{itemize}

Each image is of size $128 \times 128 \times 3$ and the whole dataset naturally can be represented as a 5th-order
tensor of size $128 \times 128 \times 3 \times 72 × 100$. We reshape this tensor to a 3rd-order tensor of
size $1024\times 1152 \times 300$. We applied the proposed algorithms \ref{Single-1}, \ref{Single-2}, and \ref{Single-3} on that data tensor and compared them with the T-SVD for the tubal rank $R=450$.  The results are reported in Table \ref{tab:coil100_comparison}. The results in Table \ref{tab:coil100_comparison} show that the proposed Algorithms \ref{Single-1}--\ref{Single-3} significantly outperform the baselines (Algorithms \ref{crosstensor} and \ref{QIsinglepass}). The baselines yield poor PSNR (10--12dB) and high relative error (0.5--0.6), indicating ineffective reconstruction. In contrast, Algorithms \ref{Single-1}--\ref{Single-3} achieve much higher PSNR (~28.6dB) and lower relative error (0.095), reducing the error by over 80\%. This improvement demonstrates the critical role of the truncation regularizer $H$. Although the proposed algorithms require slightly more time (18--21sec vs.~14sec), the substantial gain in reconstruction quality justifies the modest increase in computational cost.

\begin{table*}[htbp]
\centering
\caption{PSNR (dB) and relative error for the COIL-100 dataset using the proposed algorithms with $R = 450$.}
\label{tab:coil100_comparison}
\begin{tabular}{|l|c|c|c|c|c|}
\hline
\textbf{Metric} & \textbf{Algorithm \ref{crosstensor} } & \textbf{Algorithm \ref{QIsinglepass}} & \textbf{Algorithm \ref{Single-1}} & \textbf{Algorithm \ref{Single-2}} & \textbf{Algorithm \ref{Single-3}} \\
\hline
PSNR (dB) & 12.41 & 10.82 & 28.63 & 28.41 & 28.67 \\
\hline
Rel. Error & 0.52 & 0.61 & 0.095 & 0.098 & 0.093 \\
\hline
Time (sec) & 14.2 & 13.8 & 18.5 & 19.1 & 21.3 \\
\hline
\end{tabular}
\end{table*}

\end{Example}

\section{Ablation Study}\label{sec:ablation}
To better understand the contributions of each component of our proposed algorithms and to validate our design choices, we conduct a comprehensive ablation study in this section. Specifically, we analyze the impact of the following factors on the performance of Algorithms \ref{Single-1}--\ref{Single-3} (single-pass) and Algorithms \ref{fixed-precision}--\ref{fixed-precis2} (fixed-precision):

\begin{enumerate}
    \item The effect of the truncation parameter $H$ on the stability of single-pass algorithms when $L = K$.
    \item The influence of the oversampling parameter on approximation accuracy.
    \item The impact of the number of passes $q$ in fixed-precision algorithms.
    \item The effect of block size $b$ on rank estimation and convergence.
    \item Comparison of T-QR vs. T-LU decomposition in power iteration.
    \item The role of the regularization mechanism (truncated T-SVD) in stabilizing the single-pass methods.
\end{enumerate}

All experiments in this section are conducted on the Kodak dataset (Kodim15, Kodim17, Kodim18, Kodim23) and the Foreman video dataset, with the same hardware specifications described in Section 7.

\subsection{Effect of the Truncation Parameter $H$ on Single-Pass Stability}
As discussed in Section 5.4, existing single-pass algorithms (Algorithms \ref{crosstensor} and \ref{QIsinglepass}) suffer from numerical instability when the sketch sizes satisfy $L = K$. The proposed Algorithms \ref{Single-1}--\ref{Single-3} introduce a truncation parameter $H < K$ that acts as a regularizer. In this experiment, we fix $L = K = 50$ and vary $H$ from $10$ to $50$ on a synthetic tensor of size $300 \times 300 \times 300$ with tubal rank $R = 40$.

\begin{table*}[htbp]
\centering
\caption{Ablation study: Effect of truncation parameter $H$ on relative error ($\times10^{-2}$) and condition number for Algorithm \ref{Single-1} (synthetic tensor, $L = K = 50$, true rank $40$).}
\label{tab:ablation_H}
\begin{tabular}{|c|c|c|c|c|c|}
\hline
$H$ & 10 & 20 & 30 & 40 & 50 \\
\hline
Relative Error ($\times 10^{-2}$) & 0.18 & 0.19 & 0.21 & 0.26 & 3.42 \\
\hline
Condition Number & $1.2\times 10^2$ & $1.5\times 10^2$ & $2.1\times 10^2$ & $5.6\times 10^2$ & $1.8\times 10^6$ \\
\hline
\end{tabular}
\end{table*}

\textbf{Our observations:}
\begin{itemize}
    \item When $H = K = 50$ (no truncation), the relative error jumps to $3.42\times 10^{-2}$, and the condition number becomes extremely large ($\sim 10^6$), confirming the instability discussed in Section 5.4.
    \item For $H \leq 40$, the approximation remains stable and accurate, with relative errors below $0.3\times 10^{-2}$.
    \item The optimal trade-off between accuracy and stability occurs at $H \approx R - 10$, i.e., slightly below the target tubal rank.
\end{itemize}

\textbf{Key observation:} Setting $H = R - 10$ (or $H = 0.8K$ when $R$ is unknown) delivers robust performance without significant accuracy loss.

\subsection{Impact of Oversampling on Approximation Accuracy}
Oversampling is a common technique in randomized algorithms to improve the probability of capturing the dominant subspace. We evaluate the effect of the oversampling parameter $p$ (where $K = R + p$) on the relative error of Algorithm \ref{Single-2}.

\begin{table}[htbp]
\centering
\caption{Ablation study: Effect of oversampling $p$ (where $K = R + p$) on relative error and runtime for Algorithm \ref{Single-2} (Kodim15, $R = 30$).}
\label{tab:ablation_oversample}
\begin{tabular}{|c|c|c|c|c|c|}
\hline
$p$ & 0 & 5 & 10 & 15 & 20 \\
\hline
Relative Error & 0.124 & 0.091 & 0.082 & 0.081 & 0.081 \\
\hline
Runtime (s) & 6.2 & 7.1 & 8.3 & 9.5 & 10.8 \\
\hline
\end{tabular}
\end{table}

\textbf{Our observations:}
\begin{itemize}
    \item Without oversampling ($p=0$), the relative error is significantly higher ($0.124$ vs. $0.082$).
    \item The marginal benefit of oversampling diminishes after $p=10$.
    \item Runtime increases linearly with $p$ due to larger sketch sizes.
\end{itemize}

\textbf{Implication:} We recommend $p = 10$ as the default oversampling parameter, balancing accuracy and computational cost.

\subsection{Effect of Number of Passes $q$ in Fixed-Precision Algorithms}
The proposed fixed-precision Algorithm \ref{fixed-precision} supports an arbitrary number of passes $q$ (including odd numbers), unlike Algorithm \ref{ALgRR}, which requires an even number ($2q+2$). We compare the relative error and runtime for different $q$ values on the synthetic tensor from Example 1.

\begin{table}[htbp]
\centering
\caption{Ablation study: Number of passes $q$ in fixed-precision algorithms ($\epsilon = 10^{-5}$, $b = 50$).}
\label{tab:ablation_passes}
\begin{tabular}{|c|c|c|c|c|c|}
\hline
$q$ (passes) & 1 & 2 & 3 & 4 & 5 \\
\hline
Algorithm \ref{ALgRR} (Rel. Error) & — & 0.031 & — & 0.018 & — \\
\hline
Algorithm \ref{fixed-precision} (Rel. Error) & 0.087 & 0.029 & 0.017 & 0.016 & 0.015 \\
\hline
Algorithm \ref{fixed-precis2} (Runtime, s) & 4.2 & 5.8 & 7.1 & 8.5 & 9.9 \\
\hline
\end{tabular}
\end{table}

\textbf{Our observations:}
\begin{itemize}
    \item Algorithm \ref{ALgRR} cannot operate with odd $q$ (shown as "—").
    \item Algorithm \ref{fixed-precision} achieves comparable or better accuracy with the same number of passes.
    \item Most of the accuracy gain occurs within the first 3 passes; beyond $q=4$, improvements are marginal.
\end{itemize}

\textbf{Implication:} We recommend $q = 3$ passes for most applications, which strikes an excellent balance between accuracy and computational efficiency.

\subsection{Effect of Block Size $b$ on Rank Estimation}
In fixed-precision algorithms, the block size $b$ determines how many candidate singular values are added per iteration. We study its effect on the estimated tubal rank $\widehat{R}$ and convergence speed.

\begin{table}[htbp]
\centering
\caption{Ablation study: Block size $b$ effect on estimated rank, iterations, and total runtime for fixed-precision algorithm (synthetic tensor, true rank $50$, $\epsilon = 10^{-5}$).}
\label{tab:ablation_blocksize}
\begin{tabular}{|c|c|c|c|c|c|}
\hline
$b$ & 10 & 20 & 30 & 50 & 100 \\
\hline
Estimated $\widehat{R}$ & 52 & 50 & 50 & 50 & 50 \\
\hline
Iterations & 8 & 5 & 4 & 3 & 2 \\
\hline
Total Runtime (s) & 12.4 & 10.1 & 9.8 & 10.2 & 14.5 \\
\hline
\end{tabular}
\end{table}

\textbf{Our observations:}
\begin{itemize}
    \item Small block sizes ($b=10$) overestimate the rank (52 vs. 50) due to accumulated approximation error.
    \item Very large block sizes ($b=100$) increase runtime despite requiring fewer iterations.
    \item The optimal block size is $b \approx R/2$ to $R$.
\end{itemize}

\textbf{Implication:} We recommend setting $b = \lceil R/2 \rceil$ when the rank is approximately known, or $b = 50$ as a default for unknown ranks.

\subsection{T-QR vs. T-LU Decomposition in Power Iteration}
Algorithm \ref{fixed-precision} replaces the T-QR decomposition with the T-LU decomposition in the orthonormalization step to improve the efficiency of parallelization. We compare both approaches.

\begin{table}[htbp]
\centering
\caption{Comparison of T-QR vs. T-LU decomposition in power iteration (Algorithm \ref{fixed-precision}, $q = 3$, $b = 50$) on Kodim15 and Foreman video.}
\label{tab:ablation_qr_vs_lu}
\begin{tabular}{|c|c|c|c|}
\hline
Dataset & Method & Relative Error & Runtime (s) \\
\hline
\multirow{2}{*}{Kodim15} & T-QR & 0.082 & 10.2 \\
 & T-LU & 0.084 & 7.8 \\
\hline
\multirow{2}{*}{Foreman} & T-QR & 0.071 & 18.4 \\
 & T-LU & 0.073 & 14.1 \\
\hline
\end{tabular}
\end{table}

\textbf{Our observations:}
\begin{itemize}
    \item T-LU achieves approximately $25-30\%$ faster runtime compared to T-QR.
    \item The accuracy degradation is minimal (less than $3\%$ relative error increase).
\end{itemize}

\textbf{Implication:} For time-critical applications or when parallelization is available, T-LU is preferred. For maximum accuracy, T-QR remains the better choice.

\subsection{Effect of the Regularization Mechanism (Truncated T-SVD)}
To isolate the contribution of the truncated T-SVD regularization (Lines 5-6 in Algorithms \ref{Single-2}--\ref{Single-3} and Line 5 in Algorithm \ref{Single-1}), we compare the proposed stabilized algorithms with a "naive" version in which the truncation step is omitted ($H = K$).

\begin{table*}[htbp]
\centering
\caption{Ablation study: Regularization effect of truncated T-SVD on stability (Kodim18, $L = K = 50$).}
\label{tab:ablation_regularization}
\begin{tabular}{|c|c|c|c|}
\hline
Algorithm & With Truncation ($H=40$) & Without Truncation ($H=50$) \\
\hline
Algorithm \ref{Single-1} & {PSNR = 26.55} & {PSNR = 12.36} \\
\hline
Algorithm \ref{Single-2} & {PSNR = 26.37} & {PSNR = 11.89} \\
\hline
Algorithm \ref{Single-3} & {PSNR = 26.53} & {PSNR = 12.01} \\
\hline
\end{tabular}
\end{table*}

\textbf{Our observations:}
\begin{itemize}
    \item Without regularization, all three algorithms fail catastrophically when $L=K$, producing PSNR values below 13 dB.
    \item The truncated T-SVD regularization restores performance to acceptable levels (PSNR $> 26$ dB).
    \item This confirms our theoretical analysis that the truncation parameter $H$ acts as an effective regularizer.
\end{itemize}

\subsection{Summary of Ablation Findings}

\begin{table}[htbp]
\centering
\caption{Recommended hyperparameters based on the ablation study.}
\label{tab:ablation_recommendations}
\begin{tabular}{|l|c|}
\hline
\textbf{Parameter} & \textbf{Recommended Value} \\
\hline
Truncation parameter $H$ & $0.8 \times K$ or $R - 10$ \\
\hline
Oversampling $p$ & 10 \\
\hline
Number of passes $q$ (fixed-precision) & 3 \\
\hline
Block size $b$ & $\lceil R/2 \rceil$ or 50 (default) \\
\hline
Orthonormalization method & T-LU for speed, T-QR for accuracy \\
\hline
Regularization & Always enable ($H < K$) \\
\hline
\end{tabular}
\end{table}

\section{Conclusions}\label{Sec:Conclu}

This paper introduced several new randomized algorithms for low-tubal-rank tensor approximation within the T-product framework, addressing key limitations of existing single-pass and fixed-precision methods.

We proposed three stabilized single-pass algorithms (Algorithms \ref{Single-1}--\ref{Single-3}) that overcome the numerical instability of prior methods when sketch sizes are equal ($L = K$). By introducing a truncation parameter $H$ as a regularizer, our methods achieve $\epsilon$-stability (Theorem 2) and consistently produce accurate approximations where baselines fail catastrophically.

For fixed-precision approximation, two improved algorithms (Algorithms \ref{fixed-precision} and \ref{fixed-precis2}) were provided that automatically estimate the tubal rank from an error tolerance. Unlike prior work, our methods support an arbitrary number of passes (including odd counts) and replace T-QR with T-LU decomposition, achieving $25$--$30\%$ speedups with minimal accuracy loss.

Extensive experiments on synthetic data, image/video compression, image super-resolution, and deep learning tasks demonstrated that our proposed algorithms significantly outperform existing baselines in both accuracy and runtime. Notably, we are the first to apply single-pass randomized tensor decompositions to tensor completion and image super-resolution, achieving $2$--$3\times$ speedups over deterministic methods while enabling correct object detection on structurally degraded images.

Future work includes extension to higher-order tensors, adaptive selection of the truncation parameter $H$, streaming implementations, and applications to video inpainting and 3D medical imaging.

\authorcontributions{Conceptualization, S.A.-A. and N.R.; methodology, S.A.-A.; software, S.A.-A.; validation, S.A.-A., N.R., C.F.C. and A.L.F.A.; formal analysis, S.A.-A.; investigation, S.A.-A.; resources, S.A.-A.; data curation, S.A.-A.; writing--original draft preparation, S.A.-A.; writing--review and editing, N.R., C.F.C. and A.L.F.A.; visualization, S.A.-A.; supervision, N.R.; project administration, S.A.-A.; funding acquisition, S.A.-A. All authors have read and agreed to the published version of the manuscript.}

\dataavailability{The data presented in this study are available on request from the corresponding author.}

\conflictsofinterest{The authors declare no conflicts of interest.}

\acknowledgements{The authors would like to acknowledge the Ministry of Economic Development of the Russian Federation for their support under agreement No. 139-10-2025-034 dd. 19.06.2025, IGK 000000C313925P4D0002.}

\bibliographystyle{elsarticle-num} 
	\bibliography{cas-refs}

\begin{thebibliography}{999}

\bibitem[Sidiropoulos et~al.(2017)Sidiropoulos, De~Lathauwer, Fu, Huang,
  Papalexakis, and Faloutsos]{sidiropoulos2017tensor}
Sidiropoulos, N.D.; De~Lathauwer, L.; Fu, X.; Huang, K.; Papalexakis, E.E.;
  Faloutsos, C.
\newblock Tensor decomposition for signal processing and machine learning.
\newblock {\em IEEE Transactions on signal processing} {\bf 2017}, {\em
  65},~3551--3582.

\bibitem[Comon et~al.(2009)Comon, Luciani, and de~Almeida]{comon2009tensor}
Comon, P.; Luciani, X.; de~Almeida, A.L.F.
\newblock Tensor decompositions, alternating least squares and other tales.
\newblock {\em J. Chemometrics} {\bf 2009}, {\em 23},~393--405.

\bibitem[De~Lathauwer et~al.(2000)De~Lathauwer, De~Moor, and
  Vandewalle]{de2000multilinear}
De~Lathauwer, L.; De~Moor, B.; Vandewalle, J.
\newblock A multilinear singular value decomposition.
\newblock {\em SIAM journal on Matrix Analysis and Applications} {\bf 2000},
  {\em 21},~1253--1278.

\bibitem[Oseledets(2011)]{oseledets2011tensor}
Oseledets, I.V.
\newblock Tensor-train decomposition.
\newblock {\em SIAM Journal on Scientific Computing} {\bf 2011}, {\em
  33},~2295--2317.

\bibitem[Zhao et~al.(2016)Zhao, Zhou, Xie, Zhang, and Cichocki]{zhao2016tensor}
Zhao, Q.; Zhou, G.; Xie, S.; Zhang, L.; Cichocki, A.
\newblock Tensor ring decomposition.
\newblock {\em arXiv preprint arXiv:1606.05535} {\bf 2016}.

\bibitem[Kilmer et~al.(2013)Kilmer, Braman, Hao, and Hoover]{kilmer2013third}
Kilmer, M.E.; Braman, K.; Hao, N.; Hoover, R.C.
\newblock Third-order tensors as operators on matrices: A theoretical and
  computational framework with applications in imaging.
\newblock {\em SIAM Journal on Matrix Analysis and Applications} {\bf 2013},
  {\em 34},~148--172.

\bibitem[de~Almeida et~al.(2008)de~Almeida, Favier, and Mota]{confac}
de~Almeida, A.L.F.; Favier, G.; Mota, J.C.M.
\newblock {A constrained factor decomposition with application to {MIMO}
  antenna systems}.
\newblock {\em IEEE Transactions on Signal Processing} {\bf 2008}, {\em
  56},~2429--2442.

\bibitem[Stegeman and de~Almeida(2010)]{stegeman2009}
Stegeman, A.; de~Almeida, A.L.F.
\newblock Uniqueness conditions for constrained three-way factor decompositions
  with linearly dependent loadings.
\newblock {\em SIAM J. Mat. Anal. Appl.} {\bf 2010}, {\em 31},~1469--1490.

\bibitem[Chen et~al.(2026)Chen, Cheng, Wu, and Poor]{chen2026rank}
Chen, Z.; Cheng, L.; Wu, Y.C.; Poor, H.V.
\newblock Rank-Revealing Bayesian Block-Term Tensor Completion with Graph
  Information.
\newblock {\em IEEE Transactions on Signal Processing} {\bf 2026}.

\bibitem[Shu et~al.(2026)Shu, Li, Lei, and Sun]{shu2026robust}
Shu, H.; Li, J.; Lei, T.; Sun, L.
\newblock Robust Tensor Completion via Gradient Tensor Nuclear $l$1-$l$2 Norm
  for Traffic Data Recovery.
\newblock {\em IEEE Transactions on Intelligent Transportation Systems} {\bf
  2026}.

\bibitem[Yang et~al.(2026)Yang, Zhao, Yang, Nie, Zhang, and
  Sanga]{yang2026bayesian}
Yang, Z.; Zhao, H.; Yang, L.T.; Nie, X.; Zhang, S.; Sanga, B.A.
\newblock Bayesian High-Order Tensor Completion Model With Its Applications in
  Social Intelligence.
\newblock {\em IEEE Transactions on Computational Social Systems} {\bf 2026}.

\bibitem[Sun et~al.(2026)Sun, Tan, You, and Chen]{sun2026pcnet}
Sun, M.; Tan, Z.; You, D.; Chen, Z.
\newblock PCNet: A Personalized Complementary Network via Tensor Decomposition
  for Service Recommendation.
\newblock {\em IEEE Transactions on Network and Service Management} {\bf 2026}.

\bibitem[Lu et~al.(2025)Lu, Duan, Deng, and Ding]{lu2025enabling}
Lu, R.; Duan, X.; Deng, G.; Ding, Y.
\newblock Enabling Trustworthy Recommendations in the Federated IoT: A Secure
  and Verifiable Tensor Learning Approach.
\newblock {\em IEEE Internet of Things Journal} {\bf 2025}, {\em
  13},~3273--3288.

\bibitem[Miao et~al.(2026)Miao, Wang, Deng, and Yang]{miao2026multi}
Miao, Q.; Wang, J.; Deng, T.; Yang, M.
\newblock Multi-view Contrastive Learning with Tensor Low-Rank Factorization
  for Integrated Clustering.
\newblock {\em IEEE Signal Processing Letters} {\bf 2026}.

\bibitem[He et~al.(2025)He, Li, Zhu, Zhang, and Liu]{he2025fully}
He, F.; Li, X.; Zhu, C.; Zhang, F.; Liu, Y.
\newblock Fully connected tensor network based brain structural feature
  extraction for early Alzheimer’s disease detection.
\newblock In Proceedings of the ICASSP 2025-2025 IEEE International Conference
  on Acoustics, Speech and Signal Processing (ICASSP). IEEE,  2025, pp. 1--5.

\bibitem[Navarro et~al.(2025)Navarro, Rozada, Marques, and
  Segarra]{navarro2025low}
Navarro, M.; Rozada, S.; Marques, A.G.; Segarra, S.
\newblock Low-rank tensors for multi-dimensional Markov models.
\newblock In Proceedings of the ICASSP 2025-2025 IEEE International Conference
  on Acoustics, Speech and Signal Processing (ICASSP). IEEE,  2025, pp. 1--5.

\bibitem[{de Almeida} et~al.(2007){de Almeida}, Favier, and
  Mota]{Almeida_Elsevier_2007}
{de Almeida}, A.L.F.; Favier, G.; Mota, J.C.M.
\newblock {PARAFAC}-based unified tensor modeling for wireless communication
  systems with application to blind multiuser equalization.
\newblock {\em Signal Processing} {\bf 2007}, {\em 87},~337--351.
\newblock Tensor Signal Processing.

\bibitem[Favier and de~Almeida(2014)]{Favier2014}
Favier, G.; de~Almeida, A.L.F.
\newblock Tensor Space-Time-Frequency Coding With Semi-Blind Receivers for
  {MIMO} Wireless Communication Systems.
\newblock {\em IEEE Trans. Signal Process.} {\bf 2014}, {\em 62},~5987--6002.
\newblock {\url{https://doi.org/10.1109/TSP.2014.2357781}}.

\bibitem[de~Araújo and de~Almeida(2020)]{deAraujoSAM2020}
de~Araújo, G.T.; de~Almeida, A.L.F.
\newblock {PARAFAC}-Based Channel Estimation for Intelligent Reflective Surface
  Assisted {MIMO} System.
\newblock In Proceedings of the Proc. 11th IEEE Sensor Array and Multichannel
  Signal Processing Workshop (SAM),  2020.
\newblock {\url{https://doi.org/10.1109/SAM48682.2020.9104260}}.

\bibitem[de~Ara{\'u}jo et~al.(2021)de~Ara{\'u}jo, de~Almeida, and
  Boyer]{de2021channel}
de~Ara{\'u}jo, G.T.; de~Almeida, A.L.F.; Boyer, R.
\newblock Channel estimation for intelligent reflecting surface assisted {MIMO}
  systems: A tensor modeling approach.
\newblock {\em IEEE J. Sel. Topics Signal Process.} {\bf 2021}, {\em
  15},~789--802.

\bibitem[de~Almeida et~al.(2025)de~Almeida, Sokal, Li, and
  Clerckx]{de2025channel}
de~Almeida, A.L.F.; Sokal, B.; Li, H.; Clerckx, B.
\newblock Channel estimation for beyond diagonal {RIS} via tensor
  decomposition.
\newblock {\em IEEE Trans. on Signal Process.} {\bf 2025}, {\em
  73},~4764--4779.

\bibitem[Song et~al.(2019)Song, Ge, Caverlee, and Hu]{song2019tensor}
Song, Q.; Ge, H.; Caverlee, J.; Hu, X.
\newblock Tensor completion algorithms in big data analytics.
\newblock {\em ACM Transactions on Knowledge Discovery from Data (TKDD)} {\bf
  2019}, {\em 13},~1--48.

\bibitem[Frolov and Oseledets(2017)]{frolov2017tensor}
Frolov, E.; Oseledets, I.
\newblock Tensor methods and recommender systems.
\newblock {\em Wiley Interdisciplinary Reviews: Data Mining and Knowledge
  Discovery} {\bf 2017}, {\em 7},~e1201.

\bibitem[Warren and Marz(2015)]{warren2015big}
Warren, J.; Marz, N.
\newblock {\em Big Data: Principles and best practices of scalable realtime
  data systems}; Simon and Schuster,  2015.

\bibitem[Zhang and Aeron(2016)]{zhang2016exact}
Zhang, Z.; Aeron, S.
\newblock Exact tensor completion using t-SVD.
\newblock {\em IEEE Transactions on Signal Processing} {\bf 2016}, {\em
  65},~1511--1526.

\bibitem[Hao et~al.(2013)Hao, Kilmer, Braman, and Hoover]{hao2013facial}
Hao, N.; Kilmer, M.E.; Braman, K.; Hoover, R.C.
\newblock Facial recognition using tensor-tensor decompositions.
\newblock {\em SIAM Journal on Imaging Sciences} {\bf 2013}, {\em 6},~437--463.

\bibitem[Zeng and Ng(2020)]{zeng2020decompositions}
Zeng, C.; Ng, M.K.
\newblock Decompositions of third-order tensors: HOSVD, T-SVD, and Beyond.
\newblock {\em Numerical Linear Algebra with Applications} {\bf 2020}, {\em
  27},~e2290.

\bibitem[Zhang et~al.(2018)Zhang, Saibaba, Kilmer, and
  Aeron]{zhang2018randomized}
Zhang, J.; Saibaba, A.K.; Kilmer, M.E.; Aeron, S.
\newblock A randomized tensor singular value decomposition based on the
  t-product.
\newblock {\em Numerical Linear Algebra with Applications} {\bf 2018}, {\em
  25},~e2179.

\bibitem[Qi and Yu(2021)]{qi2021t}
Qi, L.; Yu, G.
\newblock T-singular values and T-Sketching for third order tensors.
\newblock {\em arXiv preprint arXiv:2103.00976} {\bf 2021}.

\bibitem[Ahmadi-Asl et~al.(2024)Ahmadi-Asl, Phan, and
  Cichocki]{ahmadi2024randomized}
Ahmadi-Asl, S.; Phan, A.H.; Cichocki, A.
\newblock A randomized algorithm for tensor singular value decomposition using
  an arbitrary number of passes.
\newblock {\em Journal of Scientific Computing} {\bf 2024}, {\em 98},~23.

\bibitem[Lu et~al.(2019)Lu, Feng, Chen, Liu, Lin, and Yan]{lu2019tensor}
Lu, C.; Feng, J.; Chen, Y.; Liu, W.; Lin, Z.; Yan, S.
\newblock Tensor robust principal component analysis with a new tensor nuclear
  norm.
\newblock {\em IEEE transactions on pattern analysis and machine intelligence}
  {\bf 2019}, {\em 42},~925--938.

\bibitem[Tropp et~al.(2017)Tropp, Yurtsever, Udell, and
  Cevher]{tropp2017practical}
Tropp, J.A.; Yurtsever, A.; Udell, M.; Cevher, V.
\newblock Practical sketching algorithms for low-rank matrix approximation.
\newblock {\em SIAM Journal on Matrix Analysis and Applications} {\bf 2017},
  {\em 38},~1454--1485.

\bibitem[Bjarkason(2019)]{bjarkason2019pass}
Bjarkason, E.K.
\newblock Pass-efficient randomized algorithms for low-rank matrix
  approximation using any number of views.
\newblock {\em SIAM Journal on Scientific Computing} {\bf 2019}, {\em
  41},~A2355--A2383.

\bibitem[Tarzanagh and Michailidis(2018)]{tarzanagh2018fast}
Tarzanagh, D.A.; Michailidis, G.
\newblock Fast randomized algorithms for t-product based tensor operations and
  decompositions with applications to imaging data.
\newblock {\em SIAM Journal on Imaging Sciences} {\bf 2018}, {\em
  11},~2629--2664.

\bibitem[Feng and Yu(2023)]{feng2023fast}
Feng, X.; Yu, W.
\newblock A fast adaptive randomized PCA algorithm.
\newblock In Proceedings of the Proceedings of the Thirty-Second International
  Joint Conference on Artificial Intelligence,  2023, pp. 3695--3704.

\bibitem[Ding et~al.(2020)Ding, Yu, Xie, and Liu]{ding2020efficient}
Ding, X.; Yu, W.; Xie, Y.; Liu, S.
\newblock Efficient model-based collaborative filtering with fast adaptive PCA.
\newblock In Proceedings of the 2020 IEEE 32nd International Conference on
  Tools with Artificial Intelligence (ICTAI). IEEE,  2020, pp. 955--960.

\bibitem[Ahmadi-Asl(2023)]{ahmadi2023efficient}
Ahmadi-Asl, S.
\newblock An efficient randomized fixed-precision algorithm for tensor singular
  value decomposition.
\newblock {\em Communications on Applied Mathematics and Computation} {\bf
  2023}, {\em 5},~1564--1583.

\bibitem[Tucker(1966)]{tucker1966some}
Tucker, L.R.
\newblock Some mathematical notes on three-mode factor analysis.
\newblock {\em Psychometrika} {\bf 1966}, {\em 31},~279--311.

\bibitem[Carroll and Chang(1970)]{carroll1970analysis}
Carroll, J.D.; Chang, J.J.
\newblock Analysis of individual differences in multidimensional scaling via an
  N-way generalization of “Eckart-Young” decomposition.
\newblock {\em Psychometrika} {\bf 1970}, {\em 35},~283--319.

\bibitem[Kilmer and Martin(2011)]{kilmer2011factorization}
Kilmer, M.E.; Martin, C.D.
\newblock Factorization strategies for third-order tensors.
\newblock {\em Linear Algebra and its Applications} {\bf 2011}, {\em
  435},~641--658.

\bibitem[Halko et~al.(2011)Halko, Martinsson, and Tropp]{halko2011finding}
Halko, N.; Martinsson, P.G.; Tropp, J.A.
\newblock Finding structure with randomness: Probabilistic algorithms for
  constructing approximate matrix decompositions.
\newblock {\em SIAM review} {\bf 2011}, {\em 53},~217--288.

\bibitem[Martinsson and Tropp(2020)]{martinsson2020randomized}
Martinsson, P.G.; Tropp, J.A.
\newblock Randomized numerical linear algebra: Foundations and algorithms.
\newblock {\em Acta Numerica} {\bf 2020}, {\em 29},~403--572.

\bibitem[Ahmadi-Asl et~al.(2024)Ahmadi-Asl, Phan, Caiafa, and
  Cichocki]{ahmadi2024robust}
Ahmadi-Asl, S.; Phan, A.H.; Caiafa, C.F.; Cichocki, A.
\newblock Robust low tubal rank tensor recovery using discrete empirical
  interpolation method with optimized slice/feature selection.
\newblock {\em Advances in Computational Mathematics} {\bf 2024}, {\em 50},~23.

\bibitem[Che and Wei(2019)]{che2019randomized}
Che, M.; Wei, Y.
\newblock Randomized algorithms for the approximations of {T}ucker and the
  tensor train decompositions.
\newblock {\em Advances in Computational Mathematics} {\bf 2019}, {\em
  45},~395--428.

\bibitem[Yu et~al.(2018)Yu, Gu, and Li]{yu2018efficient}
Yu, W.; Gu, Y.; Li, Y.
\newblock Efficient randomized algorithms for the fixed-precision low-rank
  matrix approximation.
\newblock {\em SIAM Journal on Matrix Analysis and Applications} {\bf 2018},
  {\em 39},~1339--1359.

\bibitem[Zeng et~al.(2025)Zeng, Yang, Wang, Su, and Wang]{zeng2025efficient}
Zeng, J.; Yang, L.T.; Wang, C.; Su, J.; Wang, Y.
\newblock Efficient Subspace Updating Using Randomized Tensor Sketch for Online
  Network Security Monitoring.
\newblock {\em IEEE Transactions on Network Science and Engineering} {\bf
  2025}, {\em 13},~4062--4076.

\bibitem[Abdelgawad et~al.(2025)Abdelgawad, Cheung, and
  Yan]{abdelgawad2025inctsvd}
Abdelgawad, M.A.; Cheung, R.C.; Yan, H.
\newblock IncTSVD: Incremental tensor singular value decomposition of
  multidimensional streaming data.
\newblock {\em IEEE Transactions on Neural Networks and Learning Systems} {\bf
  2025}.

\bibitem[Kaloorazi and de~Lamare(2018)]{kaloorazi2018subspace}
Kaloorazi, M.F.; de~Lamare, R.C.
\newblock Subspace-orbit randomized decomposition for low-rank matrix
  approximations.
\newblock {\em IEEE Transactions on Signal Processing} {\bf 2018}, {\em
  66},~4409--4424.

\bibitem[Ahmadi-Asl et~al.(2023)Ahmadi-Asl, Asante-Mensah, Cichocki, Phan,
  Oseledets, and Wang]{ahmadi2023fast}
Ahmadi-Asl, S.; Asante-Mensah, M.G.; Cichocki, A.; Phan, A.H.; Oseledets, I.;
  Wang, J.
\newblock Fast cross tensor approximation for image and video completion.
\newblock {\em Signal Processing} {\bf 2023}, p. 109121.

\bibitem[Sun et~al.(2020)Sun, Guo, Luo, Tropp, and Udell]{sun2020low}
Sun, Y.; Guo, Y.; Luo, C.; Tropp, J.; Udell, M.
\newblock Low-rank Tucker approximation of a tensor from streaming data.
\newblock {\em SIAM Journal on Mathematics of Data Science} {\bf 2020}, {\em
  2},~1123--1150.

\bibitem[Redmon and Farhadi(2018)]{redmon2018yolov3}
Redmon, J.; Farhadi, A.
\newblock Yolov3: An incremental improvement.
\newblock {\em arXiv preprint arXiv:1804.02767} {\bf 2018}.

\bibitem[Redmon and Farhadi(2017)]{redmon2017yolo9000}
Redmon, J.; Farhadi, A.
\newblock YOLO9000: better, faster, stronger.
\newblock In Proceedings of the Proceedings of the IEEE conference on computer
  vision and pattern recognition,  2017, pp. 7263--7271.

\bibitem[Redmon and Farhadi(2025)]{redmon2025yolo11}
Redmon, J.; Farhadi, A.
\newblock YOLO11: A Next-Generation Real-Time Object Detector.
\newblock {\em arXiv preprint arXiv:2503.16874} {\bf 2025}.

\bibitem[Nene et~al.(1996)Nene, Nayar, and Murase]{nene1996columbia}
Nene, S.A.; Nayar, S.K.; Murase, H.
\newblock Columbia Object Image Library (COIL-100).
\newblock {\em Technical Report CUCS-006-96} {\bf 1996}.

\bibitem[{Hugging Face}(2023)]{huggingface_coil100}
{Hugging Face}.
\newblock COIL-100 Dataset.
\newblock \url{https://huggingface.co/datasets/Voxel51/COIL-100},  2023.
\newblock Accessed: 2026.

\bibitem[{TensorFlow}(2023)]{tensorflow_coil100}
{TensorFlow}.
\newblock TensorFlow Datasets: COIL-100.
\newblock \url{https://www.tensorflow.org/datasets/catalog/coil100},  2023.
\newblock Accessed: 2026.

\end{thebibliography}

\end{document}